\documentclass[11pt]{article}

\usepackage[T1]{fontenc}
\usepackage[utf8]{inputenc}
\usepackage{lmodern}
\usepackage{amsmath,amssymb,amsthm,mathtools}
\usepackage{mathrsfs}
\usepackage{bm}
\usepackage{geometry}
\usepackage{microtype}
\usepackage{enumitem}
\usepackage{booktabs}
\usepackage{tikz}
\usetikzlibrary{positioning, fit, backgrounds, arrows.meta}
\usepackage{xcolor}
\usepackage{hyperref}

\usepackage{hyperxmp}

\hypersetup{%
  colorlinks       = true,
  linkcolor        = blue!55!black,
  citecolor        = blue!55!black,
  urlcolor         = blue!55!black,
  pdftitle    = {Autonomous Neural Adaptive Control as a Self-Consistent
                 Transfer Operator: A Barycentric State-Parameter Formulation},
  pdfauthor   = {Stefano Galatolo, Vittorio Lippi},
  pdfsubject  = {%
    A mathematical framework unifying neural adaptive control and transfer
    operator theory. Introduces the Self-Consistent Transfer Operator (STO)
    on the extended state-parameter space, proves existence and uniqueness
    of adaptive equilibria via Brouwer degree theory, establishes local
    exponential stability under time-scale separation, and derives the
    linear response identity linking the reduced learning field Jacobian
    to the slow adaptive dynamics. Validated on a reflected noisy
    Duffing-type oscillator with neural residual compensator.},
  pdfkeywords = {%
    adaptive control,
    neural network control,
    self-consistent transfer operator,
    closed-loop learning,
    barycentric formulation,
    state-parameter space,
    Brouwer degree theory,
    invariant measure,
    Wasserstein contraction,
    time-scale separation,
    learning rate optimization,
    action error,
    distributional shift,
    performative prediction,
    Duffing oscillator,
    nonlinear control,
    Lyapunov stability,
    reinforcement learning,
    online adaptation,
    persistent excitation,
    identifiability trap,
    monokinetic measure,
    frozen dynamics,
    spectral gap,
    linear response,
    reduced learning field,
    neural residual compensator,
    stochastic dynamical systems,
    ergodic theory,
    measure-valued dynamics,
    control theory,
    system identification,
    robotics,
    exoskeleton control,
    Vittorio Lippi,
    Stefano Galatolo,
    Calejo},
  pdfcreator  = {pdfLaTeX},
  pdfproducer = {pdfLaTeX},
  pdfinfo     = {%
    ORCID       = {0000-0001-5520-8974},
    Affiliation = {Calejo Hybrid Intelligence AB},
    Cite        = {Galatolo, S., Lippi, V. (2026). Autonomous Neural Adaptive
                   Control as a Self-Consistent Transfer Operator: A Barycentric
                   State-Parameter Formulation.},
    RelatedWork = {Lippi 2026 ICINCO; Galatolo CMP 2022;
                   Hauswirth et al. 2024; Vernade et al. 2026},
    License     = {CC-BY-4.0},
    CodeRepository = {https://github.com/VittorioFreiburg/sto-adaptive-control},
    SoftwareDOI = {10.5281/zenodo.24112629}%
  },
  pdfcopyright     = {Copyright 2026 Stefano Galatolo, Vittorio Lippi. CC-BY-4.0.},
  pdflicenseurl    = {https://creativecommons.org/licenses/by/4.0/},
  pdfcontacturl    = {https://orcid.org/0000-0001-5520-8974},
  pdfcontactemail  = {vittorio.lippi@calejo.ai},
  pdftype          = {Text},
  pdflang          = {en},
  pdfmetalang      = {en}
}

\newtheorem{theorem}{Theorem}[section]
\newtheorem{proposition}[theorem]{Proposition}

\newtheorem{corollary}[theorem]{Corollary}
\theoremstyle{definition}
\newtheorem{definition}[theorem]{Definition}
\newtheorem{assumption}[theorem]{Assumption}
\theoremstyle{remark}
\newtheorem{remark}[theorem]{Remark}

\newcommand{\SM}{\mathrm{SM}}
\newcommand{\Prob}{\mathcal{P}}
\newcommand{\ProbThetaOne}{\mathcal{P}_{\Theta}^{1}}
\newcommand{\STO}{\mathcal{L}}
\newcommand{\Lin}{L}
\newcommand{\Id}{\mathrm{Id}}
\newcommand{\dd}{\,\mathrm{d}}
\newcommand{\norm}[1]{\left\lVert #1\right\rVert}
\newcommand{\R}{\mathbb{R}}
\newcommand{\ThetaSet}{\Theta}
\newcommand{\Gammahat}{\widehat{\Gamma}}
\newcommand{\Z}{\mathsf{Z}}
\newcommand{\Mmono}{\mathscr{M}_{\mathrm{mono}}}

\title{Autonomous Neural Adaptive Control as a Self-Consistent Transfer Operator}
\author{Stefano Galatolo 
\footnote{Dipartimento di Matematica, Università di Pisa,
Largo Pontecorvo 5, 56127 Pisa, Italy email: stefano.galatolo@unipi.it}
\and Vittorio Lippi
\footnote{Calejo Hybrid Intelligence AB, Mellansj\"{o} 109, 92592 T\"{a}rnaby, Sweden email: vittorio.lippi@calejo.ai}}

\date{\today}

\begin{document}
\maketitle

\begin{abstract}
We study autonomous neural adaptive control, where the controller modifies the physical dynamics while the data generated by the physical dynamics itself drive subsequent parameter updates. We formulate the coupled evolution of the physical state distribution and controller parameters through a self-consistent transfer operator.

Assuming that the physical dynamics under each fixed controller admits a unique invariant distribution, we characterize adaptive equilibria through a finite-dimensional self-consistency equation for the stationary mean learning update. We establish existence using Brouwer degree, without requiring global convexity, and reduce uniqueness to uniqueness of the corresponding zero.

Under suitable differentiability and spectral-gap assumptions, we derive a finite-dimensional criterion for local exponential convergence of both the physical distribution and controller parameters at sufficiently small positive learning rates. The criterion involves the Jacobian of the reduced equation, computed using linear response to include the effect of controller-induced changes in the invariant distribution. The asymptotic slow relaxation time scales inversely with the learning rate.

We rigorously verify the assumptions and apply the theory to a noisy Duffing-type system with a neural residual compensator, providing a nontrivial example of adaptive control with provable existence, uniqueness, and local stability of an adaptive equilibrium. Numerical experiments illustrate statistical feedback and the role of time-scale separation.

\end{abstract}


\tableofcontents

\section{Introduction}

Learning-based control, including reinforcement learning and neural adaptive control, is becoming an increasingly important approach to control the behavior of complex dynamical systems, with applications ranging from robotic locomotion and manipulation to flow and fluids control \cite{TangEtAl2025}, \cite{BruntonEtAl2020}. By learning from interactions with the controlled system, these methods offer ways to design feedback policies when accurate models are unavailable or difficult to obtain \cite{Recht2019} \cite{Annaswamy2023}. Their growing use motivates fundamental mathematical questions about the dynamics generated when learning and physical evolution remain coupled.

In this setting, adaptation modifies the dynamics that produces the data, while the resulting statistical behavior drives subsequent adaptation.
Online adaptation therefore closes a feedback loop,
\[
  \text{controller parameters}
  \longrightarrow
  \text{state distribution}
  \longrightarrow
  \text{learning update}
  \longrightarrow
  \text{controller parameters}.
\]

Understanding this feedback requires a theory of the joint evolution of physical statistics and controller parameters: whether stationary configurations exist, under what conditions they are stable, and how the time scales of physical relaxation and learning determine convergence. These questions place continuously adapting control systems within the broader study of nonlinear dynamics and its statistical properties.

An important tool in this context is linear response theory,
which describes how statistical averages change to first order
under small perturbations of the dynamics; see \cite{Baladi2014} for a survey in
deterministic dynamics and
\cite{HairerMajda2010,GalatoloGiulietti2019,GalatoloLucarini2026}
for stochastic systems.
Here, we use the response of the invariant distribution associated
with a fixed controller to quantify how changes in its parameters
affect physical statistics and adaptive stability.

We study these questions through a self-consistent transfer
operator formulation of the coupled evolution of the physical
state distribution and the controller parameters.
The controller determines how the physical distribution evolves,
while the parameters are updated using the loss gradient averaged
over the current distribution.
In this description, an adaptive equilibrium consists of
a controller and an invariant state distribution for which
this averaged  loss gradient vanishes.

When each fixed controller admits a unique invariant distribution,
the equilibrium condition becomes a finite-dimensional equation:
the loss gradient, averaged over the invariant distribution
generated by that controller, must vanish.
This is a self-consistency equation for the controller parameters,
analogous to those used in the theory of self-consistent transfer
operators (\cite{APK,GNT,BKZ,BL,Selley2021Multiple,Tanzi2023}).
Under suitable assumptions, this reduction allows us to establish
existence and uniqueness of adaptive equilibria and, through linear
response, to give explicit conditions for their local exponential
stability and convergence from nearby initial states.

We focus on controllers parameterized by neural networks. Neural networks have long been used as flexible models and
controllers for nonlinear systems
\cite{NarendraParthasarathy1990,hunt1992neural,
sjoberg1995nonlinear,levine2016end}.
Because of the feedback between control, physical
dynamics, and state statistics, a neural controller may
perform differently from what its training performance
would suggest once it is deployed to control the physical
system using observations of that system's state. This deployment-induced mismatch, or \emph{Action Error}, is a critical bottleneck for dynamic physical systems. In human-inspired posture control, sensory feedback inevitably correlates the system's state with controlled inputs,  complicating the offline identification of learned predictors and necessitating online adaptation \cite{lippi2018prediction}. Similarly, in the control of lower-limb exoskeletons, static offline models struggle with phase transitions, requiring continuous neural learning to maintain smooth, transparent tracking \cite{lippi2021icinco}. 
To formalize these sequential failure modes, \cite{LippiActionError} introduces a three-part taxonomy: the initial \emph{Closed-Loop Identification Bias} caused by unmeasured process noise, the deployment-induced \emph{Action Error} which shifts the physical state distribution, and the continuous \emph{Tracking Error} encountered during online adaptation. 
The present self-consistent transfer operator framework
provides a mathematical description of this interaction
and conditions under which the physical distribution and
the learning parameters converge jointly to an adaptive
equilibrium.
Online adaptation may also face an identifiability problem:
successful regulation can reduce the excitation needed
to distinguish different controller parameters
\cite{LippiActionError,Annaswamy2023}.
Here we focus on the coupled stationarity problem underlying the  Action Error.



Related feedback mechanisms between learning and data distributions
arise in several settings. In imitation learning, the learned policy
determines the states encountered during deployment
\cite{RossGordonBagnell2011}. In performative prediction, a stable
model is optimal for the distribution it induces
\cite{PerdomoEtAl2020}, with stateful extensions retaining dependence
on the current distribution \cite{BrownHodKalemaj2022}.
Stochastic approximation with parameter-dependent Markovian data
also accounts for evolving sampling statistics
\cite{KondaTsitsiklis2003,LiWai2022,LiuLiWai2024},
while performative reinforcement learning allows policies to modify
the environment's transition and reward mechanisms
\cite{MandalEtAl2023}.
Feedback-based optimization similarly treats the optimizer and
physical process as a coupled dynamical system
\cite{HauswirthEtAl2024}.
These approaches draw on a range of mathematical tools, including
online learning and regret bounds, contraction arguments,
stochastic approximation, and Lyapunov stability theory.

Our approach uses the physical transfer operators to derive
properties of the distributional response from the underlying
dynamics. This makes topological and spectral tools available:
degree theory gives existence criteria without global convexity,
while linear response quantifies how changes in physical
statistics contribute to adaptive stability and relaxation rates.
The resulting guarantees concern local convergence of both the
physical distribution and the controller parameters, under
regularity and mixing assumptions on the physical dynamics.

We now introduce informally the mathematical formulation of the adaptive
control systems considered in this paper.
Let $\mu_n\in\Prob(X)$ denote the physical state distribution
on the state space $X$ at time $n$, and let
$\theta_n\in\ThetaSet\subset\R^p$ denote the controller parameters.
For each $\theta\in\ThetaSet$, a Markov transfer operator $P_\theta$
describes the evolution of the physical distribution under
that controller.

The controller parameters are updated using
$G(\mu,\theta)\in\R^p$, the gradient with
respect to $\theta$, of a loss function  averaged over the current physical
distribution $\mu$; see \eqref{eq:mean-gradient}.
With learning rate $\delta>0$ and preconditioning matrix
$\Gammahat$, the coupled autonomous dynamics is
\begin{equation}
  \mu_{n+1}=P_{\theta_n}\mu_n,
  \qquad
  \theta_{n+1}=\theta_n-\delta\Gammahat G(\mu_n,\theta_n).
  \label{eq:intro-pair-dynamics}
\end{equation}
Thus the controller determines the evolution of the physical
distribution, while that distribution determines the averaged
gradient used to update the controller.

The purpose of this paper is to place \eqref{eq:intro-pair-dynamics} in a
transfer-operator framework without introducing artificial stationary states
corresponding to simultaneously deployed ensembles of controllers.  We use the
extended space
\[
  \Z=X\times\ThetaSet
\]
and write \(\ProbThetaOne(\Z)\) for the ambient probabilities having finite
first moment in the parameter coordinate; the formal definition is given in
Section~2.  We embed a physical distribution--controller pair  in \(\ProbThetaOne(\Z)\) as
\[
  \iota(\mu,\theta)=\mu\otimes\delta_\theta.
\]

We call measures of the form $\mu\otimes\delta_\theta$
{\em monokinetic} and denote their collection by $\Mmono$.

For an arbitrary joint distribution \(\rho\in\ProbThetaOne(\Z)\), we extract its physical marginal \(\mu_\rho\) alongside a single effective controller from its parameter marginal. A natural and convenient choice for this controller is the barycenter
\[
  \bar\theta_\rho
  =\int_{\Z}\vartheta\,\dd\rho(x,\vartheta)
\]
which defines the single-controller reduction
\[
  \mathcal R(\rho)=(\mu_\rho,\bar\theta_\rho).
\]

With these notations we can define the joint evolution of the statistical state and the control parameters by the following operator
  
\begin{equation}
  \boxed{
  \STO_\delta
  =\iota\circ\Phi_\delta\circ\mathcal R,}
  \label{eq:intro-factorization}
\end{equation}
where
\[
  \Phi_\delta(\mu,\theta)
  =\bigl(P_\theta\mu,
         \theta-\delta\Gammahat G(\mu,\theta)\bigr).
\]

  This is a nonlinear Markov operator and has the standard
self-consistent transfer operator form
\[
  \STO_\delta(\rho)=\Lin_{\delta,\rho}\rho.
\]
Here the current statistical state $\rho$ selects a linear
Markov operator $\Lin_{\delta,\rho}$, which is then applied
to $\rho$ itself.
Thus, each selected operator is linear, while the dependence
of its selection on the evolving distribution makes the
resulting dynamics nonlinear.
\footnote{For fixed $\rho$, define
\[
  \Lin_{\delta,\rho}\nu
  :=
  \bigl(P_{\bar\theta_\rho}(\pi_X)_\#\nu\bigr)
  \otimes
  \delta_{\bar\theta_\rho
  -\delta\widehat\Gamma G(\mu_\rho,\bar\theta_\rho)}.
\]
Both the physical transfer operator and the updated controller
parameter are fixed by $\rho$. Thus this expression is linear
in $\nu$, preserves positivity and total mass, and satisfies
$\STO_\delta(\rho)=\Lin_{\delta,\rho}\rho$.
See Section~\ref{2.3} for the corresponding Markov kernel construction} 
(see \cite{GalatoloCMP} for a general introduction to this
approach).  In particular,
\[
  \STO_\delta\bigl(\ProbThetaOne(\Z)\bigr)
  \subseteq\Mmono,
\]
so every invariant probability measure is automatically monokinetic.

For a positive learning rate and an injective preconditioning
matrix $\widehat\Gamma$, adaptive equilibria are characterized by

\begin{equation}
 P_{\theta_*}\mu_*=\mu_* ,
  \qquad
  G(\mu_*,\theta_*)=0.
  \label{eq:intro-equilibrium-system}
\end{equation}

\paragraph{Main contributions.}
The Self Consistent Transfer Operator formulation allows us to apply topological, spectral,
and linear-response methods to study the existence and uniqueness
of adaptive equilibria and the convergence of nearby states
toward them.
Our main results are as follows.

\begin{enumerate}[label=(\roman*)]

  \item \emph{Stationary reduction and uniqueness.}
  For positive learning rate and injective $\Gammahat$,
  suppose that each fixed controller $\theta$ admits a unique
  invariant physical distribution $\mu_\theta$ in a prescribed
  Banach space of measures $B_s(X)$
  (Section~\ref{4}).
  Adaptive equilibria in this class then correspond bijectively
  to solutions of
  \[
    \overline G(\theta):=G(\mu_\theta,\theta)=0,
  \]
  where $\overline G$ is the loss gradient averaged over the
  invariant distribution generated by the controller itself
  (Theorem~\ref{thm:equilibrium-correspondence}).
  Uniqueness reduces to uniqueness of a zero of this
  finite-dimensional field
  (Corollary~\ref{cor:reduced-global-uniqueness});
  under the differentiability assumptions of Section~\ref{LR},
  invertibility of $D\overline G(\theta_*)$ gives local uniqueness
  (Remark~\ref{cor:local-uniqueness}).

  \item \emph{Existence by topological methods.}
  Under continuity assumptions, a nonzero Brouwer degree of
  $\overline G$ on a bounded parameter region guarantees
  an adaptive equilibrium
  (Theorem~\ref{thm:degree-existence}).
  Section~\ref{4.1} gives sufficient conditions based on
  outward-pointing behavior, comparison with simpler fields,
  or validated numerical approximations.
  The outward-pointing condition can be checked using the
  growth and regularization properties of the loss,
  without requiring global convexity or global monotonicity.

  \item \emph{Linear response and local exponential stability.}
  Linear response differentiates the stationary distribution
  $\mu_\theta$ and yields a formula for $D\overline G(\theta)$
  that separates the direct effect of changing the controller
  from the contribution of the changing physical statistics.
  This derivative also determines the slow part of the
  linearized adaptive dynamics
  (Proposition~\ref{prop:effective-jacobian-identification}).
  Under suitable differentiability and a spectral gap for the
  frozen physical operator, positivity of the real parts of all
  eigenvalues of $\Gammahat D\overline G(\theta_*)$ implies
  local exponential convergence of both the physical
  distribution and the controller parameters for every
  sufficiently small positive learning rate $\delta$
  (Theorem~\ref{thm:local-stability} and
  Remark~\ref{cor:LR-stability}).
  The slow relaxation time is of order $\delta^{-1}$.

  \item \emph{A concrete example and numerical experiments.}
  For a reflected noisy cubic (Duffing-type) control system 
  with a residual compensator whose readout weights are adapted
  while its features remain fixed, we verify the mixing and 
  differentiability assumptions and obtain explicit conditions 
  for existence, uniqueness in a parameter region, and local 
  exponential stability (Section~\ref{sec:duffing-example}).
  The experiments in Section~\ref{sec:empirical-validation}
  compare the predicted regularization-induced equilibrium
  shift and relaxation rates with numerical results, and
  isolate the statistical contribution to adaptive stability.
  Further experiments explore weak excitation, deployment
  mismatch, larger learning rates, and adaptation of internal
  neural-network weights beyond the explicitly verified setting.

\end{enumerate}

The theoretical results concern autonomous dynamics driven by
gradients averaged over the current physical distribution.
Finite-sample updates, additional learning variables,
time-dependent references, and active safety constraints
require further analysis or extensions of the model.
The barycentric selector is a modeling choice; substantial
parameter symmetries may motivate a formulation on a policy
or quotient space.

\paragraph{Relation with the literature and novelty.}
Our contribution connects autonomous adaptive control with the statistical theory of nonlinear dynamical systems
through an STO formulation.
This connection yields verifiable criteria for the existence,
uniqueness, and local stability of adaptive equilibria.

The reduction of stationary self-consistency to a
finite-dimensional equation has precedents in the study
of mean-field coupled systems. The work
\cite{APK} studies collective behavior
through a self-consistent nonlinear Frobenius-Perron equation.
A particularly explicit stationary reduction appears in \cite{BKZ}, where a scalar
self-consistency equation selects invariant densities and
enters the analysis of stability and bistability.
Related reductions also appear 
\cite{Selley2021Multiple},  characterizing stationary
absolutely continuous measures of self-consistent $\beta$-maps
through a scalar fixed-point equation.
The article \cite{BL} exploits the dependence of the coupling on a scalar observable, reduced equations and
linear-response formulas to reduce invertibility conditions
to scalar criteria, thereby developing a bifurcation theory
for stationary measures in this setting.

In \cite{GNT}, finite-dimensional reduction is combined with
validated numerics to obtain computer-assisted proofs of
stationary measures in mean-field coupled systems.
A close connection with our approach is the use of linear
response to determine the derivatives of reduced self-consistency
equations.
Here these equations take a different form, reflecting the
adaptive coupling designed to decrease a loss function.

Our stability analysis builds on the operator and spectral viewpoint
developed in
\cite{GalatoloCMP,CastorriniGalatoloTanziJNS,BL}.
In our setting, the reduced self-consistency equation simplifies
this approach: its derivative is a finite-dimensional matrix,
explicitly expressed through linear response, whose spectrum
governs the stability of slow adaptation.
This avoids the difficulties associated with differentiating
the ambient STO on infinite-dimensional spaces and makes the
resulting stability criterion more accessible to analysis
and computation.

%

As discussed above, the problem considered here shares
with performative learning and feedback-based optimization
a feedback mechanism through which learning influences
the system generating its data.
These connections deserve further investigation.
Our formulation places this adaptive-control problem
within the developed theory of statistical properties
of dynamical systems, making its tools available for
the analysis of the coupled evolution.
In particular, transfer operators describe the physical
statistics, while linear response quantifies how changes
in the controller affect the stationary distribution
and thereby contribute to adaptive stability.

{\noindent \bf Plan of the paper.}
The paper is organized as follows.
Section \ref{2} constructs the extended self-consistent transfer
operator and its reduction to the physical distribution--controller
dynamics.
Section \ref{3} characterizes adaptive equilibria and introduces an
invariance-residual interpretation of Action Error.
Section \ref{4} reduces the stationary problem to the field
\(\overline G\) and establishes existence and uniqueness criteria.
Section \ref{5} proves local exponential stability for sufficiently
small learning rates, while Section \ref{LR} uses linear response
to connect the reduced equilibrium equation with the
linearized adaptive dynamics.
Section \ref{sec:duffing-example} develops the reflected noisy Duffing-type example,
Section \ref{sec:empirical-validation} presents numerical experiments, and Section \ref{9}
discusses the conclusions and further directions.

\section{Autonomous adaptive control and the extended STO}
\label{2}


A self-consistent transfer operator (STO) is constructed from a family of linear Markov operators whose transition mechanisms depend on the state distribution itself. Specifically, an STO is a nonlinear Markov operator of the form
\begin{equation}
  \STO_\delta(\rho):=\Lin_{\delta,\rho}\rho,
  \label{eq:abstract-sto}
\end{equation}
where, for fixed \(\rho\), the frozen operator \(\Lin_{\delta,\rho}\) is
linear and Markov. This formulation typically requires the mapping \(\rho \mapsto \Lin_{\delta,\rho}\) to possess sufficient regularity, such as continuity or Lipschitz bounds, to ensure the stability of the resulting nonlinear dynamics (see \cite{GalatoloCMP}). Its nonlinearity comes from the dependence of the frozen
operator on the measure that determines it.  We now construct an extended STO
for autonomous adaptive control.  The extension is chosen so that an arbitrary
probability measure on the state--parameter space represents one physical
state distribution together with one current controller parameter, extracted
from the parameter marginal.  Consequently, after one application of the STO
the state is monokinetic in the parameter coordinate.

\subsection{Physical dynamics and learning rule}

Let $X$ be the physical state space, a Polish space equipped with its Borel
$\sigma$-algebra, and let $\ThetaSet\subset\mathbb R^p$
be the set of neural-network parameters, we suppose this is nonempty closed convex parameter set.
  For a fixed
\(\theta\in\ThetaSet\), let
\begin{equation}
  K_\theta(x,\dd y)
  \label{eq:kernel}
\end{equation}
be the Markov transition kernel of the controlled plant.  Its transfer operator
is
\begin{equation}
  (P_\theta\mu)(A)
  :=\int_X K_\theta(x,A)\,\dd\mu(x).
  \label{eq:plant-transfer}
\end{equation}
In the deterministic case,
\[
  x_{n+1}=F_{\theta}(x_n)
  =F\bigl(x_n,g(x_n,\theta)\bigr),
\]
and \(P_\theta=(F_\theta)_\#\).
Here $g(x,\theta)$ is the control input produced by the
neural controller, and $F(x,u)$ describes the physical
state update under input $u$.

Let \(\ell(x,\theta)\) be the instantaneous loss.  For the algebraic
construction we assume, for simplicity, that
\(x\mapsto\nabla_\theta\ell(x,\theta)\) is bounded and measurable for each
\(\theta\), so that the following integral is defined for every
\(\mu\in\Prob(X)\).  The same construction works on any forward-invariant
class of physical probabilities on which the integral is finite.  The
distributional or infinite-batch gradient is
\begin{equation}
  G(\mu,\theta)
  :=\int_X \nabla_\theta\ell(x,\theta)\,\dd\mu(x).
  \label{eq:mean-gradient}
\end{equation}
This is the deterministic mean-field counterpart of the single-sample SGD
update used in the numerical implementation of \cite{LippiActionError}.

We factor the learning-rate matrix as
\begin{equation}
  \Gamma=\delta\Gammahat,
  \qquad \delta\geq0,
  \label{eq:learning-rate}
\end{equation}
where \(\delta\) is the scalar learning-rate scale and \(\Gammahat\) is fixed.
The autonomous measure--parameter dynamics is
\begin{equation}
  \boxed{
  \begin{aligned}
    \mu_{n+1}&=P_{\theta_n}\mu_n,\\
    \theta_{n+1}&=\theta_n-\delta\Gammahat
                    G(\mu_n,\theta_n).
  \end{aligned}}
  \label{eq:coupled-system}
\end{equation}
Equation \eqref{eq:coupled-system} is the autonomous version of the two-time-scale mechanism in \cite{LippiActionError}: the physical distribution evolves
under the current controller and the controller parameters evolve according to
the gradient generated by the current distribution.

\subsection{Extended state space and the single-controller reduction}

We now formalize the coupled dynamics of Equation~\eqref{eq:coupled-system} as a self-consistent transfer operator of the kind introduced above. To achieve this without introducing artificial macroscopic states representing arbitrary ensembles of controllers, we represent the pair dynamics using a special class of monokinetic measures.

Set
\begin{equation}
  \Z:=X\times\ThetaSet.
  \label{eq:extended-space}
\end{equation}
A pair \((\mu,\theta)\in\Prob(X)\times\ThetaSet\) is represented by
\begin{equation}
  \iota(\mu,\theta)
  :=\mu\otimes\delta_\theta
  \in\ProbThetaOne(\Z),
  \label{eq:embedding}
\end{equation}
where \(\ProbThetaOne(\Z)\) is the ambient domain of probability measures defined in Equation~\eqref{eq:ambient-domain} below.


We call
\begin{equation}
  \Mmono
  :=\{\mu\otimes\delta_\theta:
       \mu\in\Prob(X),\ \theta\in\ThetaSet\}
  \subset\ProbThetaOne(\Z)
  \label{eq:monokinetic-set}
\end{equation}
the monokinetic set.  The terminology emphasizes that, at a given iteration,
all the physical mass carries the same controller parameter vector \(\theta\).

The physical adaptive-control system has one controller parameter at each
iteration.  Therefore, when extending the dynamics from \(\Mmono\) to an
ambient space of probability measures on \(\Z\), we associate an arbitrary
extended measure with a single effective parameter.  We assume that \(\ThetaSet\subset\R^p\) is closed and convex.  The natural
ambient domain is
\begin{equation}
  \ProbThetaOne(\Z)
  :=\left\{\rho\in\Prob(\Z):
       \int_{\Z}\norm{\vartheta}\,\dd\rho(x,\vartheta)<\infty\right\}.
  \label{eq:ambient-domain}
\end{equation}
If \(\ThetaSet\) is bounded this is simply \(\Prob(\Z)\).  Denote the two marginals
of \(\rho\in\ProbThetaOne(\Z)\) by
\begin{equation}
  \mu_\rho:=(\pi_X)_\#\rho,
  \qquad
  m_\rho:=(\pi_\Theta)_\#\rho.
  \label{eq:extended-marginals}
\end{equation}
Because \(\ThetaSet\) is closed and convex, the barycenter of every measure in
\eqref{eq:ambient-domain} belongs to \(\ThetaSet\).  Define the barycentric controller
\begin{equation}
  \bar\theta_\rho
  :=\int_{\ThetaSet}\vartheta\,\dd m_\rho(\vartheta)
  =\int_{\Z}\vartheta\,\dd\rho(x,\vartheta)
  \in\ThetaSet.
  \label{eq:barycentric-controller}
\end{equation}
The associated single-controller reduction map is
\begin{equation}
  \mathcal R(\rho):=(\mu_\rho,\bar\theta_\rho).
  \label{eq:single-controller-retraction}
\end{equation}
By construction,
\begin{equation}
  \boxed{\mathcal R\circ\iota=\Id_{\Prob(X)\times\ThetaSet}.}
  \label{eq:retraction-identity}
\end{equation}
Thus the reduction map is a left inverse of the monokinetic embedding.  In the
set-theoretic sense, the associated retraction onto the monokinetic subset is
\begin{equation}
  \Pi_{\mathrm{mono}}:=\iota\circ\mathcal R:
  \ProbThetaOne(\Z)\longrightarrow\Mmono,
  \qquad
  \Pi_{\mathrm{mono}}|_{\Mmono}=\Id_{\Mmono},
  \qquad \Pi_{\mathrm{mono}}^2=\Pi_{\mathrm{mono}}.
  \label{eq:monokinetic-retraction}
\end{equation}
No continuity of \(\Pi_{\mathrm{mono}}\) is asserted until a topology on the
ambient measure space is specified.

\begin{remark}[Choice of the single-controller selector]
The barycenter is a convenient choice, but the algebraic
construction only requires a selector satisfying
$S(\delta_\theta)=\theta$.
All such selectors induce the same dynamics on monokinetic
states, whereas ambient stability additionally requires
appropriate continuity or Lipschitz bounds.
We use the barycentric selector throughout.
\end{remark}

\subsection{Frozen and self-consistent transfer operators}\label{2.3}

Fix \(\rho\in\ProbThetaOne(\Z)\).  The measure \(\rho\) determines the current
physical distribution \(\mu_\rho\) and the current single controller
\(\bar\theta_\rho\).  Define the updated controller
\begin{equation}
  \theta^+_{\delta,\rho}
  :=\bar\theta_\rho
    -\delta\Gammahat G(\mu_\rho,\bar\theta_\rho).
  \label{eq:barycentric-update}
\end{equation}
Throughout the paper we assume that this update remains in \(\ThetaSet\). This condition is automatic when $\ThetaSet=\mathbb R^p$.
Projected updates can also be considered, but active projection constraints
change the fixed-point and differentiability formulas and are not covered by
the results below.

The frozen Markov kernel is
\begin{equation}
  \boxed{
  \mathsf K_{\delta,\rho}
  \bigl((x,\vartheta),\dd y\,\dd\eta\bigr)
  :=K_{\bar\theta_\rho}(x,\dd y)\,
    \delta_{\theta^+_{\delta,\rho}}(\dd\eta).}
  \label{eq:extended-kernel}
\end{equation}

The kernel is independent of the incoming parameter coordinate:
all physical states evolve under the selected controller
$\bar\theta_\rho$ and receive the same updated parameter
$\theta_{\delta,\rho}^{+}$.

The kernel \eqref{eq:extended-kernel} induces the frozen linear Markov operator
\begin{equation}
  \bigl(\Lin_{\delta,\rho}\nu\bigr)(E)
  :=\int_{\Z}
       \mathsf K_{\delta,\rho}(z,E)\,\dd\nu(z),
  \qquad
  E\in\mathcal B(\Z),\quad \nu\in\SM(\Z).
  \label{eq:frozen-extended-operator}
\end{equation}

For fixed $\rho$, this operator is linear, preserves positivity
and total mass, and acts on a probability measure $\nu$ with
physical marginal $\mu_\nu=(\pi_X)_\#\nu$ as

If
\(\mu_\nu=(\pi_X)_\#\nu\), then for every probability measure \(\nu\)
\begin{equation}
  \boxed{
  \Lin_{\delta,\rho}\nu
  =P_{\bar\theta_\rho}\mu_\nu
   \otimes\delta_{\theta^+_{\delta,\rho}}.}
  \label{eq:frozen-explicit-action}
\end{equation}

Its parameter marginal is the Dirac mass at
$\theta_{\delta,\rho}^{+}$, so the output has finite first
parameter moment and belongs to $\ProbThetaOne(\Z)$.

The self-consistency condition consists in transporting the same measure that
determines the frozen operator:
\begin{equation}
  \boxed{
  \STO_\delta(\rho)
  :=\Lin_{\delta,\rho}\rho
  =P_{\bar\theta_\rho}\mu_\rho
   \otimes
   \delta_{\bar\theta_\rho-
   \delta\Gammahat G(\mu_\rho,\bar\theta_\rho)}.}
  \label{eq:extended-sto}
\end{equation}
Thus \(\STO_\delta\) has the abstract self-consistent form
\eqref{eq:abstract-sto}, although the ambient extension is deliberately chosen
to collapse the parameter marginal to a single controller after one step.

Define the pair map
\begin{equation}
  \Phi_\delta(\mu,\theta)
  :=\bigl(P_\theta\mu,
          \theta-\delta\Gammahat G(\mu,\theta)\bigr).
  \label{eq:pair-map}
\end{equation}
The preceding construction can then be summarized by the factorization
\begin{equation}
  \boxed{
  \STO_\delta
  =\iota\circ\Phi_\delta\circ\mathcal R.}
  \label{eq:sto-factorization}
\end{equation}

\begin{proposition}[One-step monokinetic reduction]
\label{prop:one-step-monokinetic}
For every \(\delta\geq0\),
\begin{equation}
  \boxed{
  \STO_\delta\bigl(\ProbThetaOne(\Z)\bigr)\subseteq\Mmono.}
  \label{eq:one-step-monokinetic}
\end{equation}
In particular, every invariant probability measure of \(\STO_\delta\) is
monokinetic.  On \(\Mmono\), the STO reproduces the original adaptive dynamics:
\begin{equation}
  \STO_\delta(\mu\otimes\delta_\theta)
  =P_\theta\mu\otimes
   \delta_{\theta-\delta\Gammahat G(\mu,\theta)}.
  \label{eq:monokinetic-action}
\end{equation}
\end{proposition}


Combining \eqref{eq:retraction-identity} and \eqref{eq:sto-factorization}
gives the exact conjugacy on the monokinetic set,
\begin{equation}
  \boxed{
  \STO_\delta\circ\iota
  =\iota\circ\Phi_\delta,}
  \label{eq:conjugacy}
\end{equation}
and, more generally, for every \(n\geq1\),
\begin{equation}
  \boxed{
  \STO_\delta^n
  =\iota\circ\Phi_\delta^n\circ\mathcal R.}
  \label{eq:iterate-factorization}
\end{equation}


\section{Adaptive equilibria and Action Error}\label{3}

The one-step reduction \eqref{eq:one-step-monokinetic} makes the stationary
problem particularly rigid.  Although the self-consistent operator is defined
on the full space \(\ProbThetaOne(\Z)\), any invariant probability measure must already
belong to the monokinetic set \(\Mmono\).  Thus there is no distinction between
``ambient'' invariant measures and invariant states carrying a single
controller parameter.

\begin{definition}[Adaptive equilibrium]
A probability measure \(\rho_*\in\ProbThetaOne(\Z)\) is an adaptive equilibrium if
\[
  \STO_\delta(\rho_*)=\rho_*.
\]
By Proposition~\ref{prop:one-step-monokinetic}, every adaptive equilibrium is
necessarily of the form
\(\rho_*=\mu_*\otimes\delta_{\theta_*}\).
\end{definition}

\begin{proposition}[Fixed-point characterization of all invariant probabilities]
\label{prop:all-fixed-points}
Assume \(\delta>0\) and that \(\Gammahat\) is injective.  A probability measure
\(\rho_*\in\ProbThetaOne(\Z)\) is invariant for \(\STO_\delta\) if and only if there
exist \(\mu_*\in\Prob(X)\) and \(\theta_*\in\ThetaSet\) such that
\begin{equation}
  \rho_*=\mu_*\otimes\delta_{\theta_*}
  \qquad\text{and}\qquad
  \boxed{
  P_{\theta_*}\mu_*=\mu_* ,
  \qquad
  G(\mu_*,\theta_*)=0.}
  \label{eq:adaptive-equilibrium}
\end{equation}
In particular, the barycentric extension creates no non-monokinetic invariant
probability measures.
\end{proposition}

\begin{proof}
Suppose first that \(\rho_*\) is invariant.  Since
\(\STO_\delta(\ProbThetaOne(\Z))\subseteq\Mmono\), invariance implies
\(\rho_*\in\Mmono\), and therefore
\[
  \rho_*=\mu_*\otimes\delta_{\theta_*}
\]
for some \((\mu_*,\theta_*)\).  On the monokinetic set the operator agrees
exactly with the original measure--parameter dynamics, so
\eqref{eq:monokinetic-action} gives
\[
  \STO_\delta(\rho_*)
  =P_{\theta_*}\mu_*\otimes
   \delta_{\theta_*-\delta\Gammahat G(\mu_*,\theta_*)}.
\]
Equality with \(\rho_*\) implies equality of the physical and parameter
marginals.  Hence
\[
  P_{\theta_*}\mu_*=\mu_*,
  \qquad
  \delta\Gammahat G(\mu_*,\theta_*)=0.
\]
Since \(\delta>0\) and \(\Gammahat\) is injective, the second equality is
 equivalent to \(G(\mu_*,\theta_*)=0\).

Conversely, if \eqref{eq:adaptive-equilibrium} holds, then
\eqref{eq:monokinetic-action} immediately yields
\(\STO_\delta(\rho_*)=\rho_*\).
\end{proof}


The two conditions in \eqref{eq:adaptive-equilibrium} represent physical stationarity of the plant and parameter stationarity of the learning rule, respectively. By Proposition~\ref{prop:all-fixed-points}, the ambient extension introduces no artificial equilibria: multiple invariant measures can arise only from non-uniqueness in the physical invariant measure or the learning root.
This is particularly natural for neural parameterizations, where symmetries and
redundant parameter directions generally preclude global parameter uniqueness.

\begin{corollary}[Zero-learning-rate equilibria]
\label{cor:zero-rate-fixed-points}
At \(\delta=0\), a probability measure \(\rho_*\in\ProbThetaOne(\Z)\) is invariant for
\(\STO_0\) if and only if
\begin{equation}
  \rho_*=\mu_*\otimes\delta_{\theta_*},
  \qquad
  P_{\theta_*}\mu_*=\mu_*.
  \label{eq:zero-rate-fixed-points}
\end{equation}
Thus the learning condition disappears at zero rate and the controller
parameter remains a neutral coordinate.
\end{corollary}

\begin{proof}

For $\rho_*=\mu_*\otimes\delta_{\theta_*}$, the formula \ref{eq:monokinetic-action}
for the action of the STO on monokinetic measures, evaluated
at $\delta=0$, gives
\[
  \STO_0(\rho_*)
  =P_{\theta_*}\mu_*\otimes\delta_{\theta_*},
\]
which equals \(\rho_*\) exactly when
\(P_{\theta_*}\mu_*=\mu_*\).
\end{proof}

\subsection{Action Error as failure of invariance}
\label{sec:action-error}

Let \((\mu_{\mathrm{off}},\theta_{\mathrm{off}})\) denote the distribution and
controller parameter produced by offline training, and set
\[
  \rho_{\mathrm{off}}
  :=\mu_{\mathrm{off}}\otimes\delta_{\theta_{\mathrm{off}}}.
\]
An offline stationarity condition may give
\begin{equation}
  G(\mu_{\mathrm{off}},\theta_{\mathrm{off}})=0,
  \label{eq:offline-learning-stationarity}
\end{equation}
but this does not imply
\(P_{\theta_{\mathrm{off}}}\mu_{\mathrm{off}}=\mu_{\mathrm{off}}\).

A natural operator-level diagnostic of the deployment-induced mismatch in
\cite{LippiActionError} is therefore the physical invariance residual, which we formalize as the \emph{Action Error}:

\begin{equation}
  \mathfrak R_A
  :=P_{\theta_{\mathrm{off}}}\mu_{\mathrm{off}}
    -\mu_{\mathrm{off}}.
  \label{eq:action-residual}
\end{equation}

\section{Reduced stationary problem: existence and uniqueness}
\label{4}
In this section we establish existence and uniqueness results for the equilibrium state of our system.

Let
\[
  B_s(X)\hookrightarrow B_w(X)\subseteq\SM(X)
\]
be Banach spaces with continuous embedding and continuous total-mass
functional, and let
\(\mathcal U\subseteq\operatorname{int}(\ThetaSet)\) be open.

\begin{assumption}[Unique invariant physical branch]
\label{ass:physical-branch}
For every \(\theta\in\mathcal U\), the Markov operator \(P_\theta:B_s(X) \to B_s(X) \) has a
unique invariant probability measure in \(B_s(X)\),  
denoted
\(\mu_\theta\). 
\end{assumption}

Define the reduced learning field
\begin{equation}
  \overline G(\theta):=G(\mu_\theta,\theta),
  \qquad \theta\in\mathcal U.
  \label{eq:effective-gradient}
\end{equation}

 Next statement gives a characterization of the equilibrium states of the system as fixed  points of the dynamics induced by the reduced learning field. We remark that once the field is computed, this is a finite dimensional condition.

\begin{theorem}[Reduced zeros and invariant probabilities]
\label{thm:equilibrium-correspondence}
Assume \(\delta>0\), that \(\Gammahat\) is injective, and
Assumption~\ref{ass:physical-branch}.  Then
\begin{equation}
  \theta\longmapsto\mu_\theta\otimes\delta_\theta
  \label{eq:equilibrium-correspondence-map}
\end{equation}
defines a bijection between
\[
  \{\theta\in\mathcal U:\overline G(\theta)=0\}
\]
and the invariant probability measures of \(\STO_\delta\) whose physical marginals belong to \(B_s(X)\) and whose controller
parameter belongs to \(\mathcal U\).
\end{theorem}

\begin{proof}
If \(\overline G(\theta)=0\), then
\[
  P_\theta\mu_\theta=\mu_\theta,
  \qquad
  G(\mu_\theta,\theta)=0,
\]
so Proposition~\ref{prop:all-fixed-points} gives the invariant probability
measure \(\mu_\theta\otimes\delta_\theta\).

Conversely, let $\rho$ be an invariant probability measure
whose physical marginal belongs to $B_s(X)$ and whose
controller parameter belongs to $U$.
  Proposition~\ref{prop:all-fixed-points}
implies that
\[
  \rho=\mu\otimes\delta_\theta,
  \qquad
  P_\theta\mu=\mu,
  \qquad
  G(\mu,\theta)=0.
\]
Uniqueness of the frozen physical invariant measure gives
\(\mu=\mu_\theta\), and therefore \(\overline G(\theta)=0\).  Distinct
parameter values give distinct Dirac parameter marginals, proving injectivity.
\end{proof}

Thus, within \(\mathcal U\), existence and uniqueness of invariant probability
measures are exactly equivalent to the reduced zero set being nonempty or a
singleton.  In particular, the ambient extension introduces no additional
stationary branches.

\begin{corollary}[Uniqueness from the reduced field]
\label{cor:reduced-global-uniqueness}
Assume the hypotheses of Theorem~\ref{thm:equilibrium-correspondence}.  If
\(\overline G\) has exactly one zero \(\theta_*\) in \(\mathcal U\), then $\rho_*:=\mu_{\theta_*}\otimes\delta_{\theta_*}$
is the unique invariant probability measure of $\STO_\delta$
whose physical marginal belongs to $B_s(X)$ and whose
controller parameter lies in $U$.
\end{corollary}

\begin{proof}
This is immediate from the bijection in
Theorem~\ref{thm:equilibrium-correspondence}.
\end{proof}

We now establish an existence result for equilibrium states of the system, based on topological argument. For this we need to assume a kind of stability for the invariant measure associated to a control parameter.

\begin{assumption}
\label{ass:stationary-continuity}
The maps
\[
  \mathcal U\ni\theta\longmapsto\mu_\theta\in B_w(X),
  \qquad
  G:B_w(X)\times\mathcal U\longrightarrow\R^p
\]
are continuous with respect to the norm topology of \(B_w(X)\) and the standard Euclidean topology on \(\mathcal{U}\) and \(\R^p\). 
\end{assumption}

Under Assumption~\ref{ass:stationary-continuity}, the reduced learning field
\(\overline G:\mathcal U\to\R^p\) is automatically continuous by composition, reflecting the assumed regularity of the underlying loss function gradient. 

\begin{theorem}[Brouwer-degree existence]
\label{thm:degree-existence}
Assume the hypotheses of Theorem~\ref{thm:equilibrium-correspondence} and
Assumption~\ref{ass:stationary-continuity}.  Let
\(\Omega\subset\R^p\) be bounded and open, with
\(\overline\Omega\subset\mathcal U\), and suppose
\footnote{Here $\deg(\overline G,\Omega,0)$ denotes the Brouwer degree of the
continuous map $\overline G:\overline\Omega\to\mathbb R^p$ with respect to
the value $0$, defined when $0\notin\overline G(\partial\Omega)$.  It is an
integer-valued topological invariant that, roughly speaking, counts the zeros
of $\overline G$ in $\Omega$ with orientation and multiplicity.  In
particular, $\deg(\overline G,\Omega,0)\neq0$ implies that there exists at
least one $\theta_*\in\Omega$ such that $\overline G(\theta_*)=0$.  If all
zeros are isolated and nondegenerate, then
$\deg(\overline G,\Omega,0)=
\sum_{\overline G(\theta_i)=0}\operatorname{sgn}\det
D\overline G(\theta_i)$.}
\[
  0\notin\overline G(\partial\Omega),
  \qquad
  \deg(\overline G,\Omega,0)\neq0.
  \label{eq:nonzero-degree}
\]
Then \(\STO_\delta\) has an invariant probability measure
\[
  \rho_*=\mu_{\theta_*}\otimes\delta_{\theta_*}
\]
with \(\theta_*\in\Omega\).
\end{theorem}

\begin{proof}
Because the Brouwer degree $  \deg(\overline G,\Omega,0)\neq0$, there is
\(\theta_*\in\Omega\) with \(\overline G(\theta_*)=0\).
Theorem~\ref{thm:equilibrium-correspondence} gives the corresponding invariant
probability measure.
\end{proof}

We show a first practical way to verify that the Brouwer degree is not null and apply the existence result.

\begin{corollary}[Outward-pointing criterion]
\label{cor:outward-pointing-existence}

Assume the hypotheses of Theorem~\ref{thm:equilibrium-correspondence}
and Assumption~\ref{ass:stationary-continuity}.
Let $\overline{B_R(\theta_0)}\subset U$.
  If
\begin{equation}
  \left\langle\overline G(\theta),\theta-\theta_0\right\rangle>0
  \qquad
  \text{for every }\theta\in\partial B_R(\theta_0),
  \label{eq:outward-pointing-condition}
\end{equation}
then \(\STO_\delta\) has an invariant probability measure whose controller
parameter belongs to \(B_R(\theta_0)\).
\end{corollary}

\begin{proof}
For
\(H_t(\theta)=(1-t)(\theta-\theta_0)+t\overline G(\theta)\),
condition \eqref{eq:outward-pointing-condition} implies
\[
  \left\langle H_t(\theta),\theta-\theta_0\right\rangle>0
  \qquad
  (\theta\in\partial B_R(\theta_0),\ 0\leq t\leq1).
\]
Hence the homotopy has no boundary zero and
\[
  \deg(\overline G,B_R(\theta_0),0)
  =\deg(\theta\mapsto\theta-\theta_0,B_R(\theta_0),0)=1.
\]
\end{proof}

\subsection{Verifying the degree condition}\label{4.1}

The degree condition is natural here because it constrains the reduced field
only on the boundary of a bounded parameter region.  It allows nonconvexity,
multiple critical points, and neural-network symmetries in the interior.
In applications one normally verifies a boundary or homotopy criterion rather
than computing the degree from its definition. We now outline some different ways to verify that the degree is not null.

\paragraph{Boundary coercivity.}
A sufficient condition for
\eqref{eq:outward-pointing-condition} is
\begin{equation}
  \left\langle\overline G(\theta),\theta-\theta_0\right\rangle
  \geq a\norm{\theta-\theta_0}^2-b,
  \qquad a>0,\quad b\geq0,
  \label{eq:degree-coercive-bound}
\end{equation}
for then every \(R>\sqrt{b/a}\) contained in \(\mathcal U\) has degree one.
For the unpreconditioned learning rule
$\widehat\Gamma=\mathrm{Id}$, this means that the stationary
learning update points inward on the boundary.
For a general injective preconditioner, the degree criterion
still establishes existence, since the equilibrium equation
remains $\overline G(\theta)=0$.

In practical machine learning applications, this coercivity condition can be  satisfied by applying \(L_2\) regularization (weight decay) \cite{goodfellow2016deep}. 
The regularization adds an outward contribution proportional to $\theta$
to the reduced learning field, ensuring the Brouwer degree
boundary condition whenever this contribution dominates
the remaining terms on the boundary.
This mechanism is made explicit in the example of
Section~\ref{sec:duffing-example}.

\paragraph{Homotopy to a nominal model.}
Let $G_0:\overline\Omega\to\mathbb R^p$ be continuous,
with $0\notin G_0(\partial\Omega)$ and
$\deg(G_0,\Omega,0)\ne0$. The boundary estimate
\begin{equation}
  \sup_{\theta\in\partial\Omega}
  \norm{\overline G(\theta)-\overline G_0(\theta)}
  <
  \inf_{\theta\in\partial\Omega}\norm{\overline G_0(\theta)}
  \label{eq:degree-homotopy-bound}
\end{equation}
prevents zeros along the straight-line homotopy, and therefore
\(\deg(\overline G,\Omega,0)=\deg(\overline G_0,\Omega,0)\).  The reference
field may come from a nominal plant, a simplified controller, or an offline
model.

\paragraph{Validated approximation.}
The same argument applies to a numerical field \(\overline G_N\).  If
\begin{equation}
  \sup_{\theta\in\partial\Omega}
  \norm{\overline G(\theta)-\overline G_N(\theta)}
  <
  \inf_{\theta\in\partial\Omega}\norm{\overline G_N(\theta)},
  \label{eq:degree-validated-bound}
\end{equation}
then \(\overline G\) and \(\overline G_N\) have the same degree.  Ulam-type
approximations, continuation in the controller parameters, and interval bounds
can therefore certify existence when their boundary errors are controlled.

\begin{remark}[Multiplicity and neural symmetries]
The degree argument gives existence, not generic global uniqueness.  By
Corollary~\ref{cor:reduced-global-uniqueness}, uniqueness within
\(\mathcal U\) is equivalent to uniqueness of the relevant zero of
\(\overline G\) once the physical stationary branch is unique; the statement
is global when every invariant controller parameter lies in \(\mathcal U\).  
\end{remark}

\section{Local stability and the zero-learning-rate center family}\label{5}

The factorization \eqref{eq:sto-factorization} has an important consequence for
stability.  In fact
\eqref{eq:retraction-identity} and \eqref{eq:sto-factorization} imply the exact
semiconjugacy
\begin{equation}
  \boxed{
  \mathcal R\circ\STO_\delta
  =\Phi_\delta\circ\mathcal R,}
  \label{eq:retraction-semiconjugacy}
\end{equation}
while \eqref{eq:iterate-factorization} gives
\(
\STO_\delta^n=\iota\circ\Phi_\delta^n\circ\mathcal R
\)
for every \(n\geq1\).  Thus the  asymptotic stability problem for $\STO_\delta$ is
equivalent to the stability problem for the pair dynamics \(\Phi_\delta\).

At \(\delta=0\),
\[
  \Phi_0(\mu,\theta)=(P_\theta\mu,\theta),
\]
so the controller coordinate is neutral.  Whenever
\(P_\theta\mu_\theta=\mu_\theta\), the pair \((\mu_\theta,\theta)\) and the
corresponding measure
\begin{equation}
  \rho_\theta:=\mu_\theta\otimes\delta_\theta
  \label{eq:equilibrium-family}
\end{equation}
are equilibria.  The zero-learning-rate system therefore possesses a finite-dimensional center family rather than a single isolated equilibrium.
For positive learning rate, we now develop below a criterion ensuring that a certain fixed point of the STO is attractive and hence also locally unique.

Let \(\rho_*=\mu_*\otimes\delta_{\theta_*}\) be a selected adaptive
equilibrium, for example one supplied by
Theorem~\ref{thm:degree-existence}, with \(\mu_*\in B_s(X)\) and
\(\theta_*\in\operatorname{int}(\ThetaSet)\).  Define
\begin{equation}
  V_s(X):=\{\eta\in B_s(X):\eta(X)=0\}.
  \label{eq:zero-mass-space}
\end{equation}
On the pair space we use the local product distance
\begin{equation}
  d_{\mathrm{pair},s}
  \bigl((\mu,\theta),(\mu_*,\theta_*)\bigr)
  :=\norm{\mu-\mu_*}_s+\norm{\theta-\theta_*}.
  \label{eq:pair-strong-distance}
\end{equation}
This formula also gives the monokinetic set the structure of a metric space.

A compatible ambient metric is obtained on the class of product probability measures
\[
  \mathcal A_s
  :=
  \left\{
    \mu\otimes m:
    \mu\in\mathcal P(X)\cap B_s(X),\
    m\in\mathcal P_1(\Theta)
  \right\}
\]
by setting
\[
  d_{\mathcal A_s}(\mu\otimes m,\nu\otimes n)
  :=
  \|\mu-\nu\|_s+W_1(m,n),
\]
where $W_1$ is the Wasserstein distance induced by the
Euclidean norm on $\Theta$.
\footnote{
Since
\[
  \left\|\int_\Theta\theta\,dm(\theta)
          -\int_\Theta\theta\,dn(\theta)\right\|
  \leq W_1(m,n),
\]
the reduction map satisfies
\[
  d_{\mathrm{pair},s}(R\rho,R\sigma)
  \leq d_{\mathcal A_s}(\rho,\sigma).
\]
Moreover, on monokinetic measures,
\[
  d_{\mathcal A_s}
  (\mu\otimes\delta_\theta,\nu\otimes\delta_\vartheta)
  =
  \|\mu-\nu\|_s+\|\theta-\vartheta\|.
\]
Thus the reduction map $\mathcal R$ is $1$-Lipschitz
from $(\mathcal A_s,d_{\mathcal A_s})$ to the pair space,
while the monokinetic embedding $\iota$ is an isometry.
These properties allow us to transfer the pair-space
stability estimate to $\mathcal A_s$ in
Corollary~\ref{cor:ambient-metric-stability}.}


\begin{assumption}[Local differentiability of the adaptive dynamics]
\label{ass:local-adaptive-C1}
Equip
\[
  V_s(X):=\{\eta\in B_s(X):\eta(X)=0\}
\]
with the norm inherited from \(B_s(X)\), and equip
\[
  E_s:=V_s(X)\times\R^p
\]
with the product norm
\[
  \norm{(\eta,v)}_{E_s}
  :=\norm{\eta}_s+\norm{v},
\]
where \(\norm{\cdot}\) is the Euclidean norm on \(\R^p\).

There exist open neighborhoods
\(\mathcal V\subset V_s(X)\) and
\(\mathcal W\subset\R^p\) of zero, with
\(\theta_*+\mathcal W\subseteq\ThetaSet\), such that the maps
initially defined on physically admissible pairs admit extensions
\[
  \mathcal F:\mathcal V\times\mathcal W\longrightarrow V_s(X),
  \qquad
  \mathcal F(\eta,v)
  :=
  P_{\theta_*+v}(\mu_*+\eta)-\mu_*,
\]
and
\[
  \mathcal G:\mathcal V\times\mathcal W\longrightarrow\R^p,
  \qquad
  \mathcal G(\eta,v)
  :=
  G(\mu_*+\eta,\theta_*+v),
\]
which are continuously Fr\'echet differentiable with respect to the
product Banach topology of \(E_s\).  
The extension $\mathcal F$ is given by the displayed formula
on all of $V\times W$, using the linear action of
$P_{\theta_*+v}$ on signed measures. In particular,
\[
  D_\eta\mathcal F(0,0)[\eta]=P_{\theta_*}\eta.
\]

Thus
\[
  D\mathcal F(\eta,v)\in
  \mathcal L(E_s,V_s(X)),
  \qquad
  D\mathcal G(\eta,v)\in
  \mathcal L(E_s,\R^p),
\]
and these differentials depend continuously on \((\eta,v)\) in the
corresponding operator norms.
\end{assumption}

Under Assumption~\ref{ass:local-adaptive-C1}, the pair map
\(\Phi_\delta\), expressed in translated coordinates
\((\eta,v)=(\mu-\mu_*,\theta-\theta_*)\), is continuously Fr\'echet
differentiable near \((0,0)\).  Its differential at the equilibrium is
\begin{equation}
 D\Phi_\delta(\mu_*,\theta_*)
 \begin{pmatrix}\eta\\ v\end{pmatrix}
 =
 \begin{pmatrix}
   P_*\eta+B_*v\\[1mm]
   v-\delta\Gammahat\bigl(C_*\eta+H_*v\bigr)
 \end{pmatrix},
 \label{eq:block-linearization}
\end{equation}
where \(P_*:=P_{\theta_*}\) and

\begin{equation}
  B_*v:=D_v\mathcal F(0,0)[v]
  =\left.\frac{d}{dt}\right|_{t=0}
    P_{\theta_*+tv}\mu_*.
  \label{eq:Bstar}
\end{equation}

\begin{align}
  C_*\eta&:=D_\mu G(\mu_*,\theta_*)[\eta],
  \label{eq:Cstar}\\
  H_*v&:=D_\theta G(\mu_*,\theta_*)[v].
  \label{eq:Hstar}
\end{align}
The definition of $B_*$ only requires differentiability
of $\theta\mapsto P_\theta\mu_*$.
Under the stronger operator-norm differentiability
assumption introduced in Section~6, it agrees with
\[
  B_*v=(D_\theta P_{\theta_*}[v])\mu_*.
\]

Equivalently,
\begin{equation}
  \mathcal A_\delta
  :=D\Phi_\delta(\mu_*,\theta_*)
  =
  \begin{pmatrix}
    P_* & B_*\\
    -\delta\Gammahat C_* &
    \Id-\delta\Gammahat H_*
  \end{pmatrix}.
  \label{eq:block-matrix}
\end{equation}

The upper-right block measures the first-order effect of a controller
perturbation on the physical distribution, whereas the lower-left block is the
feedback from the physical distribution to the learning update.

Assume that the physical dynamics has a spectral gap on the zero-mass space: \footnote{Spectral gaps provide a standard mechanism for exponential
relaxation in mixing hyperbolic dynamics and in stochastic
systems. Their realization depends on the functional
setting: hyperbolic systems often require anisotropic spaces of distributions, whereas a uniform Doeblin condition gives contraction on zero-mass measures in total variation. The latter fits the present framework
directly and is verified in
Section~\ref{sec:duffing-example}.}
\begin{equation}
  r\bigl(P_*|_{V_s(X)}\bigr)<1.
  \label{eq:fast-spectral-gap}
\end{equation}
Then
\[
  R_*:=(\Id-P_*)^{-1}:V_s(X)\to V_s(X)
\]
is bounded.  Define the effective slow Jacobian
\begin{equation}
  \boxed{
  \mathcal J_*:=H_*+C_*R_*B_*.}
  \label{eq:effective-slow-jacobian}
\end{equation}
The term \(R_*B_*v\) is the stationary first-order physical displacement
induced by a parameter perturbation \(v\), and \(C_*R_*B_*v\) is the
corresponding indirect contribution to the learning update. We remark that $\mathcal J_*$ is independent of $\delta$.

\begin{theorem}[Local exponential stability of the pair dynamics]
\label{thm:local-stability}
Let \(\rho_*\) be as above.  Assume:
\begin{enumerate}[label=(S\arabic*)]
  \item Assumption~\ref{ass:local-adaptive-C1} holds;
  \item the spectral-gap condition \eqref{eq:fast-spectral-gap} on the physical system 
         holds;
  \item the matrix
        \begin{equation}
          M_*:=\Gammahat\mathcal J_*
          \label{eq:effective-slow-matrix}
        \end{equation}
        is positive stable, that is,
        \(\operatorname{Re}\lambda>0\) for every
        \(\lambda\in\sigma(M_*)\).
\end{enumerate}
Then there exists \(\delta_0>0\) such that, for every
\(0<\delta<\delta_0\), the equilibrium \((\mu_*,\theta_*)\) is locally
exponentially stable for \(\Phi_\delta\).  Equivalently, \(\rho_*\) is locally
exponentially stable for the STO on the invariant monokinetic set.
More precisely, for each such \(\delta\) there exist a neighborhood
\(\mathcal N_\delta\) of \((0,0)\) in \(V_s(X)\times\R^p\), a constant
\(C_\delta\geq1\), and \(q_\delta\in(0,1)\) such that
\begin{equation}
  \norm{\mu_n-\mu_*}_s+\norm{\theta_n-\theta_*}
  \leq C_\delta q_\delta^n
  \left(
    \norm{\mu_0-\mu_*}_s+\norm{\theta_0-\theta_*}
  \right)
  \label{eq:local-exponential-stability}
\end{equation}
whenever \(\mu_0\) is a probability measure and
\((\mu_0-\mu_*,\theta_0-\theta_*)\in\mathcal N_\delta\).

The linearized operator
$A_\delta=D\Phi_\delta(\mu_*,\theta_*)$ has exactly $p$
eigenvalues near $1$, counted with algebraic multiplicity.
These eigenvalues  $\lambda_j^{\mathrm{slow}}(\delta)$  converge to $1$ as $\delta\to0$ and
describe the slow adaptive modes associated with the
controller parameters. More precisely

\begin{equation}
  \lambda_j^{\mathrm{slow}}(\delta)
  =1-\delta\lambda_j(M_*)+o(\delta),
  \qquad j=1,\ldots,p
  \label{eq:slow-multipliers}
\end{equation}
where $\lambda_1(M_*),\ldots,\lambda_p(M_*)$ are the
eigenvalues of $M_*=\widehat\Gamma\mathcal J_*$,
counted with algebraic multiplicity.
Consequently, if
\(m_*:=\min\{\operatorname{Re}\lambda:\lambda\in\sigma(M_*)\}>0\),
the slow contraction factor is
\(1-\delta m_*+o(\delta)\), and the corresponding relaxation time is of
order \(\delta^{-1}\).
\end{theorem}

\begin{proof}
At $\delta=0$,
\[
  \mathcal A_0=
  \begin{pmatrix}
    P_* & B_*\\
    0 & \Id
  \end{pmatrix},
\]
where throughout the block calculation $P_*$ acts on
$V_s(X)$.
By \emph{(S2)}, $1\notin\sigma(P_*)$, so
$R_*=(\Id-P_*)^{-1}$ is bounded and
\[
  P_*R_*B_*+B_*=R_*B_*.
\]
Define the bounded invertible transformation
\[
  S(\xi,v):=(\xi+R_*B_*v,v),
  \qquad
  S^{-1}(\eta,v)=(\eta-R_*B_*v,v).
\]
Then
\[
\begin{aligned}
  \mathcal A_0S(\xi,v)
  &=\bigl(P_*\xi+P_*R_*B_*v+B_*v,\ v\bigr)\\
  &=\bigl(P_*\xi+R_*B_*v,\ v\bigr),
\end{aligned}
\]
and hence
\[
  S^{-1}\mathcal A_0S(\xi,v)=(P_*\xi,v).
\]
Thus
\[
  S^{-1}\mathcal A_0S
  =\operatorname{diag}(P_*,\Id).
\]
In particular, the spectral value $1$ is isolated and
semisimple, with algebraic multiplicity $p$ and eigenspace
\[
  E_c=\ker(\mathcal A_0-\Id)
     =\{(R_*B_*v,v):v\in\R^p\}.
\]

For positive $\delta$, recalling
$\mathcal J_*=H_*+C_*R_*B_*$, we obtain
\[
  \mathcal A_\delta S(\xi,v)
  =
  \begin{pmatrix}
    P_*\xi+R_*B_*v\\
    -\delta\Gammahat C_*\xi
    +(\Id-\delta\Gammahat\mathcal J_*)v
  \end{pmatrix}.
\]
Applying $S^{-1}$ subtracts $R_*B_*$ times the second
component from the first, giving
\begin{equation}
  S^{-1}\mathcal A_\delta S
  =
  \begin{pmatrix}
    P_*+\delta R_*B_*\Gammahat C_*
      & \delta R_*B_*\Gammahat\mathcal J_*\\[1mm]
    -\delta\Gammahat C_*
      & \Id-\delta\Gammahat\mathcal J_*
  \end{pmatrix}.
  \label{eq:conjugated-linearization}
\end{equation}

At $\delta=0$, the eigenspace of
$S^{-1}\mathcal A_0S=\operatorname{diag}(P_*,\Id)$
associated with the eigenvalue $1$ is
$E_0:=\{0\}\times\R^p$, with spectral projection
$Q(\xi,v)=(0,v)$.
Write
\[
  S^{-1}\mathcal A_\delta S
  =
  \operatorname{diag}(P_*,\Id)+\delta K.
\]
Applying the first-order perturbation $K$ to a vector
$(0,v)\in E_0$ and projecting back onto $E_0$ gives,
by \eqref{eq:conjugated-linearization},
\[
  QK(0,v)
  =
  (0,-\Gammahat\mathcal J_*v)
  =
  (0,-M_*v).
\]
Thus the first-order perturbation on this eigenspace
is represented by $-M_*$.

Since $\mathcal A_\delta$ depends analytically on $\delta$,
perturbation theory for the isolated spectral cluster at $1$
gives a $p$-dimensional invariant spectral subspace.
After an analytic identification with the unperturbed
center space, the restriction of $\mathcal A_\delta$
to this subspace is represented by
\[
  T_\delta=\Id-\delta M_*+O(\delta^2).
\]
Continuity of the eigenvalues of
$(T_\delta-\Id)/\delta=-M_*+O(\delta)$, counted with
algebraic multiplicity, yields
\[
  \lambda_j^{\mathrm{slow}}(\delta)
  =1-\delta\lambda_j(M_*)+o(\delta),
  \qquad j=1,\ldots,p,
\]
which is \eqref{eq:slow-multipliers}.
This argument does not require $M_*$ to be diagonalizable.

By \emph{(S3)}, $\Re\lambda_j(M_*)>0$ for every $j$.
Consequently,
\[
  \bigl|\lambda_j^{\mathrm{slow}}(\delta)\bigr|^2
  =1-2\delta\Re\lambda_j(M_*)+o(\delta)<1
\]
for all sufficiently small positive $\delta$.
The remaining spectrum stays in a disk of radius strictly
less than one: indeed, choose $a$ with $r(P_*)<a<1$.
Spectral separation and upper semicontinuity of the spectrum
under operator-norm perturbations imply that, for sufficiently
small $\delta$, every spectral value outside the slow cluster
lies in $\{z:|z|<a\}$.
Hence, after choosing $\delta_0>0$ sufficiently small,
\[
  r(\mathcal A_\delta)<1
  \qquad (0<\delta<\delta_0).
\]

For each such fixed $\delta$, the principle of linearized
stability for continuously Fréchet differentiable maps on
Banach spaces gives
\eqref{eq:local-exponential-stability}.
Indeed, $r(\mathcal A_\delta)<1$ allows an equivalent norm
in which $\|\mathcal A_\delta\|<1$; continuity of the
derivative then makes the translated nonlinear map a
contraction on a sufficiently small neighborhood of zero.

Finally, with
\[
  m_*:=\min_{\lambda\in\sigma(M_*)}\Re\lambda>0,
\]
the spectral radius of the slow cluster satisfies
\[
  r_{\mathrm{slow}}(\delta)
  =1-\delta m_*+o(\delta).
\]
Thus its asymptotic linear relaxation time satisfies
\[
  \frac{1}{-\log r_{\mathrm{slow}}(\delta)}
  \sim\frac{1}{\delta m_*},
  \qquad \delta\to0^+,
\]
and is therefore of order $\delta^{-1}$.
\end{proof}
{

The factorization transfers pair-space stability to the
space $\mathcal A_s$ of product probability measures,
equipped with the metric $d_{\mathcal A_s}$ defined above.

\begin{corollary}[Stability in the product-measure metric]
\label{cor:ambient-metric-stability}
Assume the hypotheses of
Theorem~\ref{thm:local-stability}, and fix
$0<\delta<\delta_0$.
Let $\rho_*=\mu_*\otimes\delta_{\theta_*}$.
There exists $\varepsilon_\delta>0$ such that every
$\rho_0=\mu_0\otimes m_0\in\mathcal A_s$ satisfying
\[
  d_{\mathcal A_s}(\rho_0,\rho_*)<\varepsilon_\delta
\]
has iterates in $\mathcal A_s$ and satisfies
\begin{equation}
  d_{\mathcal A_s}(\STO_\delta^n\rho_0,\rho_*)
  \leq
  C_\delta q_\delta^n
  d_{\mathcal A_s}(\rho_0,\rho_*),
  \qquad n\geq0,
  \label{eq:ambient-metric-contraction}
\end{equation}
where $C_\delta\geq1$ and $q_\delta\in(0,1)$ are as in
Theorem~\ref{thm:local-stability}.
Moreover, $\STO_\delta^n\rho_0\in\Mmono$ for every $n\geq1$.
Thus $\rho_*$ is locally exponentially stable on
$(\mathcal A_s,d_{\mathcal A_s})$.
\end{corollary}

\begin{proof}
Let
\[
  \bar\theta_0:=\int_\Theta\theta\,dm_0(\theta).
\]
Since
\[
  \|\bar\theta_0-\theta_*\|
  \leq
  \int_\Theta\|\theta-\theta_*\|\,dm_0(\theta)
  =W_1(m_0,\delta_{\theta_*}),
\]
we have
\[
  \|\mu_0-\mu_*\|_s+\|\bar\theta_0-\theta_*\|
  \leq d_{\mathcal A_s}(\rho_0,\rho_*).
\]
Choose $\varepsilon_\delta>0$ sufficiently small that
this bound ensures
\[
  (\mu_0-\mu_*,\bar\theta_0-\theta_*)\in\mathcal N_\delta.
\]
Write
\[
  (\mu_n,\theta_n)
  :=\Phi_\delta^n(\mu_0,\bar\theta_0).
\]
By \eqref{eq:iterate-factorization}, for every $n\geq1$,
\[
  \STO_\delta^n\rho_0=\mu_n\otimes\delta_{\theta_n}.
\]
Consequently, Theorem~\ref{thm:local-stability} gives
\[
\begin{aligned}
  d_{\mathcal A_s}(\STO_\delta^n\rho_0,\rho_*)
  &=\|\mu_n-\mu_*\|_s+\|\theta_n-\theta_*\|\\
  &\leq
  C_\delta q_\delta^n
  \bigl(\|\mu_0-\mu_*\|_s+\|\bar\theta_0-\theta_*\|\bigr)\\
  &\leq
  C_\delta q_\delta^n d_{\mathcal A_s}(\rho_0,\rho_*).
\end{aligned}
\]
The case $n=0$ follows from $C_\delta\geq1$.
\end{proof}

\section{Linear response and the geometry of the reduced field}\label{LR}

The preceding sections describe the same adaptive equilibrium from two
complementary viewpoints.  The stationary problem is encoded by the reduced
field
\(
  \overline G(\theta)=G(\mu_\theta,\theta)
\),
whereas the local dynamics is governed by the block linearization of the pair
map \(\Phi_\delta\).  Linear response provides the bridge between these two
descriptions: it differentiates the stationary physical branch
\(\theta\mapsto\mu_\theta\) and shows that the Jacobian of
\(\overline G\) is exactly the effective slow Jacobian that appears in the
stability calculation.

This identification is entirely a statement on the physical measure--parameter
space.  In particular, no differentiation of the ambient barycentric STO is
needed.  The reduction map \(\mathcal R\) is relevant for arbitrary initial
measures through the factorization \eqref{eq:iterate-factorization}, but once a
state has been reduced to a single controller the stationary sensitivity and
the local stability mechanism are those of the pair dynamics.

\begin{assumption}[Strong regularity of the physical family]
\label{ass:physical-regularity}
In addition to Assumption~\ref{ass:physical-branch}, suppose that:
\begin{enumerate}[label=(\roman*)]
  \item the map
        \[
          \mathcal U\ni\theta\longmapsto
          P_\theta\in\mathcal L(B_s(X),B_s(X))
        \]
        is continuously Fr\'echet differentiable;
  \item for every \(\theta\in\mathcal U\), the restriction
        \((\Id-P_\theta)|_{V_s(X)}\) is boundedly invertible.  We write
        \begin{equation}
          R_\theta:=(\Id-P_\theta)^{-1}:V_s(X)\to V_s(X).
          \label{eq:physical-resolvent}
        \end{equation}
\end{enumerate}
\end{assumption}

The term $(D_\theta P_\theta[v])\mu_\theta$ describes the
first-order change produced by perturbing the controller
while keeping the input distribution fixed.
The inverse $R_\theta$ converts this perturbation into the
corresponding stationary response.
When $r(P_\theta|_{V_s(X)})<1$, this response also admits
the convergent representation
\[
  D_\theta\mu_\theta[v]
  =\sum_{n=0}^{\infty}
    P_\theta^n(D_\theta P_\theta[v])\mu_\theta.
\]

The following proposition is a classical linear-response
result, stated here in the present functional framework.
We include its short proof for completeness.
The response holds in the strong norm and in the
Fr\'echet sense, so the first-order approximation is
uniform with respect to the direction of the parameter
perturbation.

\begin{proposition}
\label{prop:physical-LR}
Under Assumptions~\ref{ass:physical-branch} and
\ref{ass:physical-regularity}, the map
\[
  \mathcal U\ni\theta\longmapsto\mu_\theta\in B_s(X)
\]
is continuously Fr\'echet differentiable.  For every \(v\in\R^p\),
\begin{equation}
  \boxed{
  D_\theta\mu_\theta[v]
  =R_\theta
  \bigl(D_\theta P_\theta[v]\bigr)\mu_\theta.}
  \label{eq:LR-formula}
\end{equation}
\end{proposition}

\begin{proof}
Fix \(\theta_0\in\mathcal U\) and work on the affine space of unit-mass
elements \(\mu_{\theta_0}+V_s(X)\).  Define
\[
  \mathcal F(\theta,\eta)
  :=P_\theta(\mu_{\theta_0}+\eta)-(\mu_{\theta_0}+\eta),
  \qquad \eta\in V_s(X).
\]
Because every \(P_\theta\) preserves total mass,
\(\mathcal F(\theta,\eta)\in V_s(X)\).  Moreover,
\[
  \mathcal F(\theta_0,0)=0,
  \qquad
  D_\eta\mathcal F(\theta_0,0)=P_{\theta_0}-\Id.
\]
The last operator is boundedly invertible on \(V_s(X)\) by
Assumption~\ref{ass:physical-regularity}.  The implicit-function theorem therefore gives a unique local continuously
differentiable mass-one invariant branch in the strong topology.  For nearby
\(\theta\), let \(\widetilde\mu_\theta\) denote this branch.  The prescribed
invariant probability \(\mu_\theta\) also belongs to \(B_s(X)\), and
\(\widetilde\mu_\theta-\mu_\theta\in V_s(X)\) satisfies
\[
  (\Id-P_\theta)(\widetilde\mu_\theta-\mu_\theta)=0.
\]
Invertibility of \(\Id-P_\theta\) on \(V_s(X)\) therefore gives
\(\widetilde\mu_\theta=\mu_\theta\).  Since \(\theta_0\) is arbitrary,
the stationary branch is continuously differentiable on \(\mathcal U\).

Differentiating \(P_\theta\mu_\theta=\mu_\theta\) in the direction \(v\)
gives
\[
  \bigl(D_\theta P_\theta[v]\bigr)\mu_\theta
  +P_\theta D_\theta\mu_\theta[v]
  =D_\theta\mu_\theta[v].
\]
Since \(\mu_\theta(X)=1\), one has
\(D_\theta\mu_\theta[v]\in V_s(X)\).  Differentiating mass preservation of
\(P_\theta\) also shows
\( (D_\theta P_\theta[v])\mu_\theta\in V_s(X)\), so the resolvent is applied
to an element of its stated domain.  Applying \(R_\theta\) proves
\eqref{eq:LR-formula}.
\end{proof}

The term
\(
  (D_\theta P_\theta[v])\mu_\theta
\)
is the instantaneous first-order change in the physical distribution produced
by a controller perturbation.  The resolvent \(R_\theta\) accumulates this
perturbation through the stationary physical dynamics.  Thus
\eqref{eq:LR-formula} is the stationary statistical sensitivity associated with
the controlled plant.

\begin{assumption}[Differentiability of the learning field]
\label{ass:learning-C1}
In addition to Assumption~\ref{ass:stationary-continuity}, suppose that the
restriction
\[
  G:B_s(X)\times\mathcal U\to\R^p
\]
is continuously Fr\'echet differentiable in both variables.
\end{assumption}

\begin{proposition}[Derivative of the reduced stationary field]
\label{prop:reduced-field-derivative}
Under Assumptions~\ref{ass:physical-regularity} and
\ref{ass:learning-C1}, the reduced field
\(\overline G:\mathcal U\to\R^p\) is continuously Fr\'echet differentiable
and, for every \(v\in\R^p\),
\begin{equation}
  D\overline G(\theta)[v]
  =D_\mu G(\mu_\theta,\theta)
     \bigl[D_\theta\mu_\theta[v]\bigr]
   +D_\theta G(\mu_\theta,\theta)[v].
  \label{eq:effective-gradient-derivative}
\end{equation}
Equivalently, using linear response,
\begin{equation}
\boxed{
\begin{aligned}
  D\overline G(\theta)[v]
  &=D_\theta G(\mu_\theta,\theta)[v]\\
  &\quad+
  D_\mu G(\mu_\theta,\theta)
  \left[
    R_\theta
    \bigl(D_\theta P_\theta[v]\bigr)\mu_\theta
  \right].
\end{aligned}}
  \label{eq:effective-gradient-LR}
\end{equation}
\end{proposition}

\begin{proof}
The differentiability follows from
Proposition~\ref{prop:physical-LR} and the chain rule applied to
\(\overline G(\theta)=G(\mu_\theta,\theta)\).  Formula
\eqref{eq:effective-gradient-derivative} is the chain rule itself, and
\eqref{eq:effective-gradient-LR} follows from \eqref{eq:LR-formula}.
\end{proof}

Formula \eqref{eq:effective-gradient-LR} separates the direct parameter
sensitivity of the learning rule from the indirect sensitivity caused by the
controller-induced displacement of the stationary state distribution.  The
second contribution is precisely the feedback effect that is absent if the
data distribution is treated as fixed during optimization.

\begin{proposition}[Stationary--dynamical Jacobian identity]
\label{prop:effective-jacobian-identification}
Let \(\rho_*=\mu_*\otimes\delta_{\theta_*}\) be an adaptive equilibrium at
which the assumptions of this section and
Assumption~\ref{ass:local-adaptive-C1} hold.  Then
\begin{equation}
  \boxed{
  D\overline G(\theta_*)
  =H_*+C_*R_*B_*
  =\mathcal J_*.}
  \label{eq:reduced-jacobian}
\end{equation}
\end{proposition}

\begin{proof}
At \(\theta_*\), Proposition~\ref{prop:physical-LR} gives
\[
  D_\theta\mu_{\theta_*}[v]
  =R_*B_*v.
\]
Substituting this identity into
\eqref{eq:effective-gradient-derivative}, and using the definitions
\eqref{eq:Cstar}--\eqref{eq:Hstar}, gives
\[
  D\overline G(\theta_*)[v]
  =H_*v+C_*R_*B_*v
\]
for every \(v\in\R^p\).
\end{proof}

\begin{remark}\label{cor:LR-stability} The identity \eqref{eq:reduced-jacobian} is the main bridge between the
stationary and dynamical parts of the paper.  The same matrix controls both the
local geometry of the reduced equilibrium equation
\(\overline G(\theta)=0\) and the first-order displacement of the slow
multipliers from the neutral value \(1\).  In particular, the stability
criterion can be stated entirely in terms of the reduced stationary field.

By this identity, the stability condition in
Theorem~\ref{thm:local-stability} can be written as
positive stability of $\Gammahat D\overline G(\theta_*)$.
Under the remaining hypotheses of that theorem, this
gives local exponential stability for sufficiently
small positive learning rates.
\end{remark}
\begin{remark}\label{cor:local-uniqueness} If $D\overline G(\theta_*)$ is invertible, the
inverse-function theorem implies that $\theta_*$ is
an isolated zero of $\overline G$.
For positive learning rate and injective $\Gammahat$,
Theorem~\ref{thm:equilibrium-correspondence} then gives
local uniqueness of the adaptive equilibrium among
those with physical marginal in $B_s(X)$.
We remark that this conclusion does not require the stability condition
of Theorem~\ref{thm:local-stability}.
\end{remark}

\paragraph{A convenient positivity condition.}
Although $G(\mu,\theta)$ is a loss gradient for fixed $\mu$,
the reduced field $\overline G$ need not itself be a gradient:
the dependence of $\mu_\theta$ on $\theta$ can make
$D\overline G(\theta)$ nonsymmetric.

If \(\Gammahat\) is symmetric positive definite and the symmetric part of
\(D\overline G(\theta_*)\) is positive definite, then
\(\Gammahat D\overline G(\theta_*)\) is positive stable.  Indeed,
\(\Gammahat D\overline G(\theta_*)\) is similar to
\[
  \Gammahat^{1/2}
  D\overline G(\theta_*)
  \Gammahat^{1/2},
\]
whose symmetric part is positive definite; hence all of its eigenvalues have
strictly positive real part.

\paragraph{Degree and local index.}
If the zeros of \(\overline G\) in \(\Omega\) are isolated and
nondegenerate, then

\[
  \deg(\overline G,\Omega,0)
  =
  \sum_{\substack{\theta_i\in\Omega\\
                   \overline G(\theta_i)=0}}
  \operatorname{sgn}\det D\overline G(\theta_i).
\]

By Proposition~\ref{prop:effective-jacobian-identification},
the local Brouwer index at a nondegenerate adaptive
equilibrium is $\operatorname{sgn}\det\mathcal J_*$.
Thus the effective slow Jacobian determines both
the local index and, through the spectrum of
$\Gammahat\mathcal J_*$, the stability criterion
for small positive learning rates.
These local conclusions complement the global
existence criterion of Theorem~\ref{thm:degree-existence}.

\section{Concrete reflected noisy cubic (Duffing-type) control system}
\label{sec:duffing-example}

We apply the theory to a stochastic cubic control model inspired by
\cite{LippiActionError}. We retain its physical coefficients and baseline
controller, replace hard saturation by reflection on a compact interval,
and use the distributional learning field of \eqref{eq:mean-gradient}.
These choices allow us to verify the mixing and differentiability hypotheses,
prove existence of adaptive equilibria, and obtain explicit sufficient
conditions for uniqueness and local stability.

\subsection{Model and learning rule}
\label{sec:duffing-loss}

Let $X=[-2,2]$ and $\ThetaSet=\R^p$. The input-driven plant is
\begin{equation}
  x_{n+1}=\mathfrak r_2\!\left(1.05x_n+0.15x_n^3+u_n+\xi_n\right),
  \label{eq:duffing-stochastic-plant}
\end{equation}
where $\mathfrak r_2:\R\to X$ is the continuous $8$-periodic reflection
map, equal to the identity on $X$, with alternating slopes $+1$ and $-1$
on consecutive intervals of length $4$. The independent noise variables
$\xi_n$ have a common $C^2$ density satisfying
\begin{equation}
  q_\sigma(z)>0\quad(z\in\R),
  \qquad q_\sigma,q_\sigma',q_\sigma''\in L^1(\R).
  \label{eq:duffing-noise-assumption}
\end{equation}
Gaussian noise of variance $\sigma^2>0$ is the principal example.

We use the control $u_n=-0.25x_n-g(x_n,\theta_n)$: the baseline changes
the linear coefficient to $0.8$, while the compensator aims to cancel
the deterministic residual
\begin{align}
  d(x)&:=0.15x^3,
  \label{eq:duffing-residual}\\
  \norm{d}_{L^\infty(X)}&=1.2.
  \label{eq:duffing-residual-bound}
\end{align}
The compensator has the linear readout form
\begin{equation}
  g(x,\theta)=\theta^{\mathsf T}\psi(x),
  \qquad \psi=(\psi_1,\ldots,\psi_p)^{\mathsf T}:X\to\R^p,
  \label{eq:duffing-neural-readout}
\end{equation}
where the features $\psi_j$ are fixed $C^1$ functions. Learning updates
their coefficients, so $g$ may be nonlinear in $x$ but is linear in
$\theta$. The features may be outputs of fixed hidden layers of a neural
network; alternatively, $\psi(x)=(1,x,x^2,x^3)^{\mathsf T}$ gives a
polynomial compensator that can represent $d$ exactly. Set
\begin{equation}
  S_\psi:=\sup_{x\in X}\norm{\psi(x)}<\infty.
  \label{eq:duffing-Spsi}
\end{equation}
The closed-loop dynamics is therefore
\begin{equation}
  x_{n+1}=\mathfrak r_2\!\left(
    0.8x_n+d(x_n)-g(x_n,\theta_n)+\xi_n\right).
  \label{eq:duffing-adaptive-plant}
\end{equation}
Under exact compensation, $g(\cdot,\theta)=d$, its update is
$0.8$-Lipschitz in $x_n$ for each fixed noise realization. Since noise
persists, the stationary object is a probability distribution, not the
point $0$. Learning targets the residual $d$, rather than the process noise.

For a frozen controller, write
\begin{align}
  F_\theta(x)&:=0.8x+d(x)-\theta^{\mathsf T}\psi(x),
  \label{eq:duffing-Ftheta}\\
  x_{n+1}&=\mathfrak r_2\!\left(F_\theta(x_n)+\xi_n\right).
  \label{eq:reflected-duffing}
\end{align}
Let $P_\theta$ be its transfer operator and $\mu_\theta$ its unique
invariant probability, whose existence is established below.

The learning objective penalizes the compensation error and, optionally,
departure from a nominal parameter $\theta_c\in\R^p$:
\begin{align}
  \ell_0(x,\theta)&:=\frac12\bigl(g(x,\theta)-d(x)\bigr)^2,
  \label{eq:duffing-compensation-loss}\\
  \ell(x,\theta)&:=\ell_0(x,\theta)
    +\frac\lambda2\norm{\theta-\theta_c}^2,\qquad \lambda\ge0.
  \label{eq:duffing-regularized-loss}
\end{align}
Since $\psi$ is independent of $\theta$,
\begin{align}
  \nabla_\theta\ell_0(x,\theta)
    &=\bigl(g(x,\theta)-d(x)\bigr)\psi(x),
  \label{eq:duffing-grad-loss}\\
  D_\theta^2\ell_0(x,\theta)&=\psi(x)\psi(x)^{\mathsf T}\succeq0.
  \label{eq:duffing-hessian-loss}
\end{align}
Here $A\succeq B$ means that $A-B$ is positive semidefinite;
the quadratic form of the last matrix is $(v^{\mathsf T}\psi(x))^2$.
Thus the distributional and reduced learning fields are
\begin{align}
  G(\mu,\theta)
    &=\int_X\bigl(g(x,\theta)-d(x)\bigr)\psi(x)\,\dd\mu(x)
      +\lambda(\theta-\theta_c),
  \label{eq:duffing-learning-field-explicit}\\
  \overline G(\theta)&=h(\theta)+\lambda(\theta-\theta_c),
  \qquad h(\theta):=\int_X\nabla_\theta\ell_0(x,\theta)\,\dd\mu_\theta(x).
  \label{eq:duffing-reduced-field}
\end{align}
The autonomous pair dynamics is
\begin{equation}
  \mu_{n+1}=P_{\theta_n}\mu_n,
  \qquad \theta_{n+1}=\theta_n-\delta\Gammahat G(\mu_n,\theta_n),
  \qquad \delta>0.
  \label{eq:duffing-parameter-update}
\end{equation}

\subsection{Mixing and linear response}

The frozen transition kernel and its density on $(-2,2)$ are
\begin{align}
  K_\theta(x,A)
    &:=\int_\R\mathbf1_A\!\left(\mathfrak r_2(F_\theta(x)+z)\right)
      q_\sigma(z)\,\dd z,
  \label{eq:duffing-kernel}\\
  k_\theta(x,y)
    &=\sum_{m\in\mathbb Z}\left[
      q_\sigma(y+8m-F_\theta(x))
      +q_\sigma(4-y+8m-F_\theta(x))\right].
  \label{eq:duffing-folded-density}
\end{align}
The noise is absolutely continuous, so reflection creates no boundary atoms.
Fix a bounded open parameter region $\mathcal U\Subset\R^p$, and set
\begin{align}
  M_{\mathcal U}
    &:=\sup_{\theta\in\overline{\mathcal U},\,x\in X}|F_\theta(x)|,
  \label{eq:duffing-MU}\\
  c_{\mathcal U}
    &:=\min_{|z|\le2+M_{\mathcal U}}q_\sigma(z)>0,
  \qquad \varepsilon_{\mathcal U}:=4c_{\mathcal U}\in(0,1).
  \label{eq:duffing-doeblin-constant}
\end{align}
We work on $B_s(X)=B_w(X)=\SM(X)$ with the total-variation norm,
and write $V(X)=\{\eta\in\SM(X):\eta(X)=0\}$ and $m_X$ for normalized
Lebesgue measure on $X$.

\begin{proposition}[Uniform frozen mixing]
\label{prop:duffing-doeblin}
For every $\theta\in\overline{\mathcal U}$,
\begin{equation}
  K_\theta(x,A)\ge\varepsilon_{\mathcal U}m_X(A)
  \qquad(x\in X,\ A\in\mathcal B(X)).
  \label{eq:duffing-doeblin}
\end{equation}
The operator $P_\theta$ has a unique invariant probability $\mu_\theta$,
and
\begin{align}
  \norm{P_\theta\eta}_{\rm TV}
    &\le(1-\varepsilon_{\mathcal U})\norm{\eta}_{\rm TV},
    \qquad \eta\in V(X),
  \label{eq:duffing-TV-contraction}\\
  r(P_\theta|_{V(X)})&\le1-\varepsilon_{\mathcal U}<1,
  \label{eq:duffing-spectral-gap}\\
  \norm{R_\theta}_{V(X)\to V(X)}&\le\varepsilon_{\mathcal U}^{-1},
    \qquad R_\theta:=(\Id-P_\theta)^{-1}|_{V(X)}.
  \label{eq:duffing-resolvent-bound}
\end{align}
\end{proposition}

\begin{proof}
The first $m=0$ term in \eqref{eq:duffing-folded-density} gives
$k_\theta(x,y)\ge q_\sigma(y-F_\theta(x))\ge c_{\mathcal U}$.
Integration yields \eqref{eq:duffing-doeblin}. The Doeblin decomposition
then gives the contraction on zero-mass measures and, by the contraction
principle on probabilities, the unique invariant law. The spectral-gap
and resolvent bounds follow from the contraction and the Neumann series.
\end{proof}

Positive noise supplies mixing uniformly on bounded parameter regions,
without any hyperbolicity assumption on the cubic drift. It also yields
strong differentiability. Indeed, set
\begin{align}
  M_F&:=\sup_{x\in X}\norm{D_\theta F_\theta(x)}=S_\psi,
  \label{eq:duffing-MF}\\
  C_P&:=\norm{q_\sigma'}_{L^1}S_\psi.
  \label{eq:duffing-CP}
\end{align}

\begin{proposition}[Differentiability of the frozen family]
\label{prop:duffing-P-differentiability}
The map $\theta\mapsto P_\theta\in\mathcal L(\SM(X),\SM(X))$
is continuously Fr\'echet differentiable, with
\begin{equation}
  \norm{D_\theta P_\theta[v]}_{\rm TV\to TV}\le C_P\norm{v}.
  \label{eq:duffing-DP-bound}
\end{equation}
The stationary branch is continuously Fr\'echet differentiable and satisfies
\begin{align}
  D_\theta\mu_\theta[v]
    &=R_\theta\bigl(D_\theta P_\theta[v]\bigr)\mu_\theta,
  \label{eq:duffing-LR}\\
  \norm{D_\theta\mu_\theta[v]}_{\rm TV}
    &\le\frac{C_P}{\varepsilon_{\mathcal U}}\norm{v},
    \qquad \theta\in\mathcal U.
  \label{eq:duffing-Dmu-bound}
\end{align}
\end{proposition}

\begin{proof}
Write $K_\theta(x,\cdot)=(\mathfrak r_2)_\#\nu_{\theta,x}$, where
$\dd\nu_{\theta,x}(z)=q_\sigma(z-F_\theta(x))\,\dd z$. Differentiation gives
\[
  D_\theta\nu_{\theta,x}[v]
  =-q_\sigma'(z-F_\theta(x))D_\theta F_\theta(x)[v]\,\dd z.
\]
The $L^1$ Taylor remainder is bounded uniformly in $x$ by
$\frac12\norm{q_\sigma''}_{L^1}S_\psi^2\norm{v}^2$.
Since push-forward contracts total variation, this proves operator-norm
differentiability and \eqref{eq:duffing-DP-bound}; the same translation
estimate gives continuity of the derivative. Proposition~\ref{prop:physical-LR},
together with Proposition~\ref{prop:duffing-doeblin}, gives the remaining claims.
\end{proof}

For Gaussian noise, $C_P=\sqrt{2/\pi}\,S_\psi/\sigma$; these strong-norm
estimates are therefore not uniform as $\sigma\to0$.

\subsection{Existence, uniqueness, and stability}

Existence can be checked directly from the model, without computing
the stationary measures.

\begin{proposition}[Constructive existence of an adaptive equilibrium]
\label{prop:duffing-constructive-existence}
Let $\lambda>0$ and
\begin{equation}
  r>r_{\rm ex}:=\frac{(\norm{\theta_c}S_\psi+1.2)S_\psi}{\lambda}.
  \label{eq:duffing-explicit-existence-radius}
\end{equation}
Then
\begin{align}
  \langle\overline G(\theta),\theta-\theta_c\rangle&>0,
    \qquad \theta\in\partial B_r(\theta_c),
  \label{eq:duffing-explicit-outward}\\
  \deg(\overline G,B_r(\theta_c),0)&=1.
  \label{eq:duffing-explicit-degree}
\end{align}
There exists a zero $\theta_*\in B_r(\theta_c)$, and every zero lies in
$\overline B_{r_{\rm ex}}(\theta_c)$. For every $\delta>0$ and injective
$\Gammahat$, this gives an invariant probability of the ambient STO:
\begin{equation}
  \rho_*:=\mu_{\theta_*}\otimes\delta_{\theta_*}.
  \label{eq:duffing-equilibrium}
\end{equation}
\end{proposition}

\begin{proof}
For $w=\theta-\theta_c$, positivity of the squared term gives
\begin{equation}
\begin{aligned}
  \langle\overline G(\theta),w\rangle
  &=\lambda\norm{w}^2+\int_X(w^{\mathsf T}\psi)^2\,\dd\mu_\theta
    +\int_X(\theta_c^{\mathsf T}\psi-d)w^{\mathsf T}\psi\,\dd\mu_\theta\\
  &\ge\lambda\norm{w}^2-(\norm{\theta_c}S_\psi+1.2)S_\psi\norm{w}.
\end{aligned}
  \label{eq:duffing-radial-coercivity}
\end{equation}
This is strictly positive when $\norm{w}>r_{\rm ex}$, proving both
the localization of all zeros and the outward condition. The straight-line
homotopy to $\theta\mapsto\theta-\theta_c$ is nonzero on the boundary,
so the degree is $1$. The existence property of Brouwer degree gives a zero,
and the equilibrium correspondence yields \eqref{eq:duffing-equilibrium}.
\end{proof}

In particular, for $\theta_c=0$ the radius condition is simply
\begin{equation}
  r>1.2S_\psi/\lambda.
  \label{eq:duffing-zero-anchor-radius}
\end{equation}

To obtain uniqueness and stability, we compare the regularization with
the statistical feedback. Choose
\begin{align}
  \mathcal U_T&:=B_T(0),\qquad T>\norm{\theta_c}+r,
  \label{eq:duffing-UT}\\
  |F_\theta(x)|&\le\widehat M_T:=2.8+TS_\psi,
    \qquad \theta\in\overline{\mathcal U_T},\ x\in X.
  \label{eq:duffing-MT-explicit}
\end{align}
A valid Doeblin constant is
$\underline\varepsilon_T:=4\min_{|z|\le4.8+TS_\psi}q_\sigma(z)>0$.
For Gaussian noise this and the derivative bound become
\begin{align}
  \underline\varepsilon_T
    &=\frac{4}{\sqrt{2\pi}\sigma}
      \exp\!\left(-\frac{(4.8+TS_\psi)^2}{2\sigma^2}\right),
  \label{eq:duffing-epsilonT-explicit}\\
  \norm{q_\sigma'}_{L^1}&=\frac{\sqrt{2/\pi}}{\sigma},
  \qquad C_P=\frac{\sqrt{2/\pi}}{\sigma}S_\psi.
  \label{eq:duffing-gaussian-CP}
\end{align}
Set
\begin{align}
  M_{1,T}&:=(TS_\psi+1.2)S_\psi
    \ge\sup_{x\in X,\,\theta\in\overline{\mathcal U_T}}
      \norm{\nabla_\theta\ell_0(x,\theta)},
  \label{eq:duffing-M1T-explicit}\\
  L_{\rm fb}(T)&:=M_{1,T}\frac{C_P}{\underline\varepsilon_T}.
  \label{eq:duffing-Lfb}
\end{align}
Differentiating the reduced field gives
\begin{equation}
  Dh(\theta)[v]
  =\int_X\psi\psi^{\mathsf T}v\,\dd\mu_\theta
   +\int_X(g(\cdot,\theta)-d)\psi\,\dd(D_\theta\mu_\theta[v]).
  \label{eq:duffing-Dh}
\end{equation}
The first term is positive semidefinite; the second has operator norm
at most $L_{\rm fb}(T)$. Writing $\operatorname{sym}A=(A+A^{\mathsf T})/2$,
we obtain
\begin{equation}
  \operatorname{sym}D\overline G(\theta)
  \succeq\bigl(\lambda-L_{\rm fb}(T)\bigr)\Id,
  \qquad\theta\in\mathcal U_T.
  \label{eq:duffing-symmetric-lower-bound}
\end{equation}
Thus sufficiently strong regularization dominates the statistical feedback,
with a criterion involving only known model bounds.

\begin{theorem}[Uniqueness in $\mathcal U_T$ and local exponential stability]
\label{thm:duffing-main-results}
Choose $r$ as in \eqref{eq:duffing-explicit-existence-radius} and
$T>\norm{\theta_c}+r$. If
\begin{equation}
  \lambda>L_{\rm fb}(T),
  \label{eq:duffing-feedback-dominance}
\end{equation}
then $\theta_*$ is the unique zero of $\overline G$ in $\mathcal U_T$.
For every $\delta>0$ and injective $\Gammahat$, $\rho_*$ is consequently
the unique invariant probability of $\STO_\delta$ whose controller
parameter belongs to $\mathcal U_T$.

If $\Gammahat$ is symmetric positive definite, there is $\delta_0>0$
such that, for each $0<\delta<\delta_0$, the pair equilibrium is locally
exponentially stable: there are a neighborhood $\mathcal N_\delta$,
$C_\delta\ge1$ and $q_\delta\in(0,1)$ such that
\begin{equation}
\begin{aligned}
  \norm{\mu_n-\mu_{\theta_*}}_{\rm TV}+\norm{\theta_n-\theta_*}
  \le C_\delta q_\delta^n
  \bigl(\norm{\mu_0-\mu_{\theta_*}}_{\rm TV}+\norm{\theta_0-\theta_*}\bigr)
\end{aligned}
  \label{eq:duffing-local-exp-convergence}
\end{equation}
for initial pairs in $\mathcal N_\delta$. In particular,
\begin{equation}
  \sup_{x\in X}|g(x,\theta_n)-g(x,\theta_*)|
  \le S_\psi\norm{\theta_n-\theta_*}\longrightarrow0.
  \label{eq:duffing-controller-convergence}
\end{equation}
The slow multipliers satisfy
\begin{equation}
  \lambda_j^{\rm slow}(\delta)
  =1-\delta\lambda_j\!\left(\Gammahat D\overline G(\theta_*)\right)
    +o(\delta),\qquad j=1,\ldots,p,
  \label{eq:duffing-slow-multipliers}
\end{equation}
with asymptotic slow relaxation time of order $\delta^{-1}$.
\end{theorem}

\begin{proof}
Integrating \eqref{eq:duffing-symmetric-lower-bound} along segments in
the convex set $\mathcal U_T$ gives
\begin{equation}
  \langle\overline G(\theta_1)-\overline G(\theta_2),\theta_1-\theta_2\rangle
  \ge(\lambda-L_{\rm fb}(T))\norm{\theta_1-\theta_2}^2.
  \label{eq:duffing-strong-monotonicity}
\end{equation}
This proves uniqueness, since the existence ball lies in $\mathcal U_T$.
At the zero, $\operatorname{sym}D\overline G(\theta_*)\succ0$ and
$\Gammahat\succ0$ imply positive stability of
$\Gammahat D\overline G(\theta_*)$.
Propositions~\ref{prop:duffing-doeblin} and
\ref{prop:duffing-P-differentiability} provide the spectral gap and joint
$C^1$ regularity of $(\mu,\theta)\mapsto P_\theta\mu$; the learning field
is affine in $\mu$ and $C^1$ in $\theta$. Theorem~\ref{thm:local-stability}
and Proposition~\ref{prop:effective-jacobian-identification} therefore yield
the stability estimate and slow multipliers. The controller estimate follows
from the linear readout and \eqref{eq:duffing-Spsi}.
\end{proof}

For fixed $T>\norm{\theta_c}$, the conditions on $r$ and $\lambda$ hold
for all sufficiently large $\lambda$. The bound may nevertheless be very
conservative, especially for small noise: its failure does not imply instability.

\begin{remark}[Role of the learning rate]
\label{rem:duffing-learning-rate}
For $\delta>0$ and injective $\Gammahat$, equilibria are determined by
$\overline G(\theta_*)=0$, independently of $\delta$. The learning rate
affects stability and convergence speed. Convergence of the adaptive pair
implies neither convergence of noisy sample paths to $0$ nor exact compensation.
\end{remark}

\subsection{Exact compensation and small regularization}

The response formula makes the feedback contribution explicit:
\begin{equation}
\begin{aligned}
  D\overline G(\theta)[v]
  &=\lambda v+\int_X\psi\psi^{\mathsf T}v\,\dd\mu_\theta\\
  &\quad+\int_X\nabla_\theta\ell_0(x,\theta)\,
    \dd\!\left[R_\theta(D_\theta P_\theta[v])\mu_\theta\right](x).
\end{aligned}
  \label{eq:duffing-explicit-reduced-jacobian}
\end{equation}
At equilibrium, this is precisely the effective slow Jacobian:
\begin{equation}
  D\overline G(\theta_*)=H_*+C_*R_*B_*=\mathcal J_*.
  \label{eq:duffing-jacobian-identification}
\end{equation}
For any bounded neighborhood $\mathcal U$ of $\theta_*$, the sharper estimate
\[
  \operatorname{sym}D\overline G(\theta_*)\succeq
  \left(\lambda+
    \lambda_{\min}\!\left(\int_X\psi\psi^{\mathsf T}\,\dd\mu_{\theta_*}\right)
    -\norm{g(\cdot,\theta_*)-d}_{L^\infty(X)}
      S_\psi\frac{C_P}{\varepsilon_{\mathcal U}}\right)\Id
\]
shows that statistical feedback is controlled by the equilibrium misfit,
rather than the bound $M_{1,T}$ over an entire parameter region.

\begin{remark}[The exactly realizable case ($\lambda=0$)]
\label{rem:duffing-realizable-case}
If $g(\cdot,\theta_*)=d$ and $\lambda=0$, then
\[
  C_*=0,\qquad
  \mathcal J_*=H_*=\int_X\psi(x)\psi(x)^{\mathsf T}\,\dd\mu_{\theta_*}(x).
\]
Positive definiteness of this Gram matrix and of $\Gammahat$ therefore
gives local exponential stability for sufficiently small positive learning
rates, without regularization.
\end{remark}

The next result shows that this stable regime persists under a small
quadratic penalty, with an equilibrium shift of order $\lambda$.

\begin{proposition}[Persistence under small regularization]
\label{prop:duffing-small-regularization}
Assume that $\psi_1,\ldots,\psi_p$ are linearly independent as functions
on $X$ and that $g(\cdot,\theta_{\mathrm{opt}})=d$ for some
$\theta_{\mathrm{opt}}\in\R^p$. Define
\[
  Q_{\mathrm{opt}}:=\int_X\psi(x)\psi(x)^{\mathsf T}
    \,\dd\mu_{\theta_{\mathrm{opt}}}(x),
  \qquad \overline G_\lambda(\theta):=h(\theta)+\lambda(\theta-\theta_c).
\]
For fixed $\theta_c\in\R^p$ and symmetric positive definite $\Gammahat$,
there are $\lambda_0>0$, a neighborhood $U_{\mathrm{opt}}$ of
$\theta_{\mathrm{opt}}$, and a $C^1$ branch $\theta_*^\lambda$,
$0\le\lambda<\lambda_0$, giving the unique zero of $\overline G_\lambda$
in $U_{\mathrm{opt}}$. Moreover,
\[
  \theta_*^0=\theta_{\mathrm{opt}},\qquad
  \theta_*^\lambda=\theta_{\mathrm{opt}}
    -\lambda Q_{\mathrm{opt}}^{-1}(\theta_{\mathrm{opt}}-\theta_c)+o(\lambda),
  \qquad
  \norm{g(\cdot,\theta_*^\lambda)-d}_{L^\infty(X)}=O(\lambda).
\]
For each such $\lambda$, $(\mu_{\theta_*^\lambda},\theta_*^\lambda)$ is
locally exponentially stable for all sufficiently small positive learning rates.
\end{proposition}

\begin{proof}
Invariance and the Doeblin minorization give
$\mu_{\theta_{\mathrm{opt}}}\ge\varepsilon_{\mathcal U}m_X$ on a bounded
parameter neighborhood. For $v\ne0$, continuity and linear independence
of the features imply
\[
  v^{\mathsf T}Q_{\mathrm{opt}}v
  \ge\varepsilon_{\mathcal U}\int_X|v^{\mathsf T}\psi(x)|^2\,\dd m_X(x)>0.
\]
Thus $Q_{\mathrm{opt}}\succ0$. Exact compensation gives
$h(\theta_{\mathrm{opt}})=0$ and
$Dh(\theta_{\mathrm{opt}})=Q_{\mathrm{opt}}$, since the response term in
\eqref{eq:duffing-Dh} vanishes. The implicit-function theorem gives the
locally unique branch, and differentiating its defining equation yields
\[
  Q_{\mathrm{opt}}
    \left.\frac{d\theta_*^\lambda}{d\lambda}\right|_{\lambda=0}
    +\theta_{\mathrm{opt}}-\theta_c=0.
\]
This proves the expansion; the misfit bound follows from
$\norm{g(\cdot,\theta_*^\lambda)-d}_{L^\infty(X)}
\le S_\psi\norm{\theta_*^\lambda-\theta_{\mathrm{opt}}}$.
Finally, $D\overline G_\lambda(\theta_*^\lambda)\to Q_{\mathrm{opt}}$
as $\lambda\to0$, so its symmetric part stays positive definite after
reducing $\lambda_0$. Theorem~\ref{thm:local-stability} and
Proposition~\ref{prop:effective-jacobian-identification}, with the frozen
mixing and differentiability proved above, give local exponential stability.
\end{proof}

\section{Numerical Experiments}
\label{sec:empirical-validation}

We now present numerical experiments on the reflected noisy cubic control system of Section~\ref{sec:duffing-example} to illustrate the analytical results of Sections~\ref{5}--\ref{sec:duffing-example} and examine regimes beyond our formal hypotheses.

Section~\ref{sec:numerical-implementation} summarizes the numerical implementation of the coupled dynamics, contrasting the exact distributional operator $\Phi_\delta$ (discretized on an Ulam grid) with finite-sample particle simulations ($N<\infty$) driven by a one-step tracking-error gradient estimator. We quantify two effects introduced by finite-sample online learning: boundary reflection bias on the Gaussian tails and the $O(\sigma\sqrt{\delta/N})$ stochastic fluctuation floor. The experiments are organized into two parts:

\begin{itemize}
    \item \emph{Quantitative verification of the local theory (Section~\ref{sec:part1-illustration}):} In the autonomous adaptation protocol, $(\mu_n,\theta_n)$ jointly evolve under $\Phi_\delta$ from $(\mu_0,\theta_0)$ at $\sigma=0.18$. These experiments test the predictions of Sections~\ref{5}--\ref{sec:duffing-example}: existence and the $O(\lambda)$ equilibrium shift under small regularization (Propositions~\ref{prop:duffing-constructive-existence} and~\ref{prop:duffing-small-regularization}), local exponential stability for small learning rates (Theorem~\ref{thm:local-stability}), the statistical linear-response contribution $C_*R_*B_*$ to the effective Jacobian $\mathcal J_* = D\overline G(\theta_*)$ (Proposition~\ref{prop:effective-jacobian-identification}), and the slow relaxation rates
    \[
      \lambda_j^{\rm slow}(\delta)
      =
      1-\delta\lambda_j\!\left(
        \Gammahat D\overline G(\theta_*)
      \right)+o(\delta),
    \]
    which yield an asymptotic relaxation time of order $\delta^{-1}$.
    \item \emph{Extension to practical control regimes (Sections~\ref{sec:part2-practice} and~\ref{sec:nonconvex_extension}):} In robotic and control applications \cite{lippi2018prediction,lippi2021icinco,DulacArnoldEtAl2021,BrunkeEtAl2022}, controllers are often pre-trained offline under conservative conditions and then deployed online over finite horizons, under weak persistent excitation \cite{AnnaswamyFradkov2021,VernadeEtAl2026}, with relatively large learning rates, or with nonlinear network parameterizations. Section~\ref{sec:part2-practice} examines weak excitation, a two-phase deployment protocol in which a parameter $\theta_{\rm off}$ identified under a conservative baseline gain ($-0.80$) is released online under a relaxed gain ($-0.25$) \cite{LippiActionError}, and large learning rates pushed across and beyond the linear stability threshold
    \[
      \delta_c=\frac{2}{\lambda_{\max}(H_*)}
    \]
    in the unregularized, realizable setting ($\Gammahat=\Id$, $C_*=0$). Section~\ref{sec:nonconvex_extension} tests a multilayer perceptron compensator.
\end{itemize}

Figure~\ref{fig:block_diagram} shows the closed-loop control architecture.

\begin{figure}[htpb]
\centering
\begin{tikzpicture}[
    auto,
    block/.style={rectangle, draw, thick, align=center, minimum height=1.2cm},
    sum/.style={circle, draw, thick, inner sep=0pt, minimum size=6mm},
    line/.style={draw, thick, -{Stealth[length=2.5mm]}},
    optional/.style={rectangle, draw, thick, dashed, align=center, minimum height=1.2cm, fill=white}
]

\node [sum] (sum_u) at (0,0) {$\Sigma$};

\node [block] (duffing) at (3.5, 0) {\textbf{Cubic Plant Dynamics} \\ $1.05x_n + 0.15x_n^3 + u_n$};

\node [sum] (sum_noise) at (6.8, 0) {$\Sigma$};

\node [block] (reflect) at (9.5, 0) {\textbf{Reflection} \\ $\mathfrak{r}_2(\cdot)$};

\node (output) at (13, 0) {};

\node [align=center] (proc_noise) at (6.8, 1.8) {Process \\ Noise $\xi_n$};

\begin{scope}[on background layer]
    \node [draw=black!70, thick, dashed, fill=green!5, inner sep=4mm, fit=(duffing) (sum_noise) (reflect), label={[font=\bfseries, yshift=-1mm]below:Physical Plant}] (plant_box) {};
\end{scope}

\node [block] (baseline) at (5, -2.5) {\textbf{Baseline Controller} \\ $u_{\rm base}(x_n) = -0.25x_n$};
\node [optional] (nn) at (5, -4.5) {\textbf{Neural Compensator} \\ $g(x_n, \theta_n)$ \\ {\footnotesize (Polynomial or MLP)}};

\draw [line] (sum_u) -- node [above] {$u_n$} (duffing);
\draw [line] (duffing) -- (sum_noise);
\draw [line] (proc_noise) -- (sum_noise);
\draw [line] (sum_noise) -- (reflect);
\draw [thick] (reflect) -- (12,0) coordinate (y_coord);
\draw [line] (12,0) -- node [above] {$x_n$} (output);

\draw [thick] (y_coord) |- (baseline.east);
\draw [line] (y_coord) |- (nn.east);

\draw [line] (baseline.west) -| node[pos=0.9, right] {$+$} (sum_u.south);
\draw [line] (nn.west) -- (-1.5, -4.5) |- node[pos=0.9, above] {$-$} (sum_u.west);

\end{tikzpicture}
\caption{Block diagram of the closed-loop system. The plant has cubic drift $1.05x_n + 0.15x_n^3$, stabilized locally around $x_{\rm ref}=0$ by $u_{\rm base}(x_n) = -0.25x_n$, perturbed by $\xi_n \sim \mathcal{N}(0,\sigma^2)$, and folded onto $X=[-2,2]$ by $\mathfrak{r}_2(\cdot)$. In Part~1 (Figures~\ref{fig:part1_nominal_convergence} and~\ref{fig:part1_timescale_and_reg}), $\sigma=0.18$. In Part~2, $\sigma$ is reduced to $0.10$ and $0.05$ to study the invariance residual $\mathfrak{R}_A$, ill-conditioning of the Gram matrix, and loss of linear stability. The compensator $g(x_n,\theta_n)$ acts with opposite sign to cancel $d(x_n)$ and is parameterized by a polynomial basis (Figures~\ref{fig:part1_nominal_convergence}--\ref{fig:part2_density_shattering}) or a multilayer perceptron (Figure~\ref{fig:mlp_convergence}).}
\label{fig:block_diagram}
\end{figure}
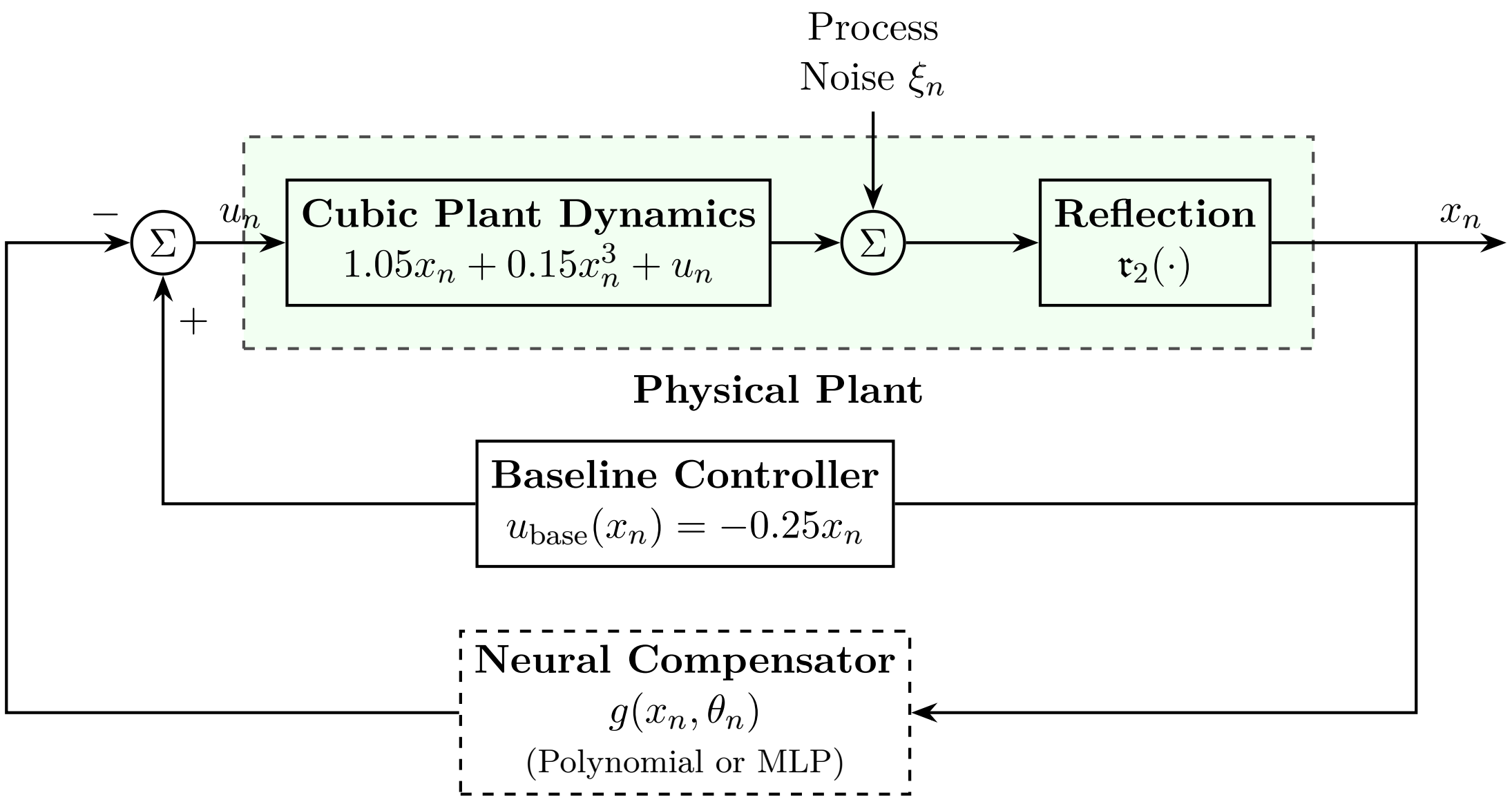

\subsection{Numerical Implementation}
\label{sec:numerical-implementation}

Consider the reflected stochastic plant of Section~\ref{sec:duffing-example} with baseline controller $u_{\rm base}(x) = -0.25x$, target cubic residual $d(x) = 0.15x^3$, and Gaussian process noise $\xi_n \sim \mathcal{N}(0,\sigma^2)$.

The bound $\norm{d}_{L^\infty(X)} = 1.2$ on $X=[-2,2]$ sets the existence radius $r_{\rm ex} = (\norm{\theta_c}S_\psi + 1.2)S_\psi/\lambda$ in Proposition~\ref{prop:duffing-constructive-existence}. We parameterize the compensator $g(x,\theta) = \theta^{\mathsf T}\psi(x)$ with the polynomial feature map
\begin{equation}
    \psi(x) = [1, x, x^2, x^3]^{\mathsf T}.
\end{equation}
Since $\psi(x)$ spans $d(x)$, the unique parameter satisfying $g(\cdot,\theta_{\rm opt}) \equiv d(\cdot)$ is
\begin{equation}
    \theta_{\rm opt} = [0, 0, 0, 0.15]^{\mathsf T}.
\end{equation}
Reflection symmetry of $\mu_{\theta_{\rm opt}}$ decouples the unregularized Gram matrix $H_* = Q_{\rm opt} = \int_X \psi\psi^{\mathsf T}\dd\mu_{\theta_{\rm opt}}$ into even and odd blocks, where the odd block is
\[
  H_{\rm odd} = \begin{pmatrix} m_2 & m_4 \\ m_4 & m_6 \end{pmatrix}, \qquad m_{2k} = \int_X x^{2k}\dd\mu_{\theta_{\rm opt}}(x).
\]
Approximating $\mu_{\theta_{\rm opt}}$ by the unreflected stationary Gaussian law of variance $\sigma_x^2 = \sigma^2/0.36$ gives $\det(H_{\rm odd}) = 6\sigma_x^8$ and
\begin{equation}
    \lambda_{\min}(H_*) \approx 6\sigma_x^6 \approx 128.6\,\sigma^6.
    \label{eq:duffing-sigma6-law}
\end{equation}
This estimate yields $\lambda_{\min}(H_*) \approx 4.37\times 10^{-3}$ at $\sigma=0.18$ (Section~\ref{sec:part1-illustration}) and $\lambda_{\min}(H_*) \approx 1.29\times 10^{-4}$ at $\sigma=0.10$ (Section~\ref{sec:part2-practice}).

To monitor uniform convergence as in \eqref{eq:duffing-controller-convergence}, we evaluate on a uniform grid of $K=500$ points $\{x_k\}_{k=1}^K \subset [-2, 2]$ the quantity
\begin{equation}
    \mathcal{E}_\infty^{\rm grid}(\theta_n) := \sqrt{\frac{1}{K}\sum_{k=1}^K \bigl|g(x_k, \theta_n) - g(x_k, \theta_{\rm opt})\bigr|^2}.
\end{equation}
Setting $\widehat H_K:=\frac{1}{K}\sum_{k=1}^K\psi(x_k)\psi(x_k)^{\mathsf T} \succ 0$, $c_K:=\lambda_{\min}(\widehat H_K)^{1/2}$, and $S_{\psi,K}:=\max_{1\le k\le K}\norm{\psi(x_k)}$ gives the norm equivalence
\begin{equation}
    c_K\norm{\theta-\theta_{\rm opt}} \le \mathcal{E}_\infty^{\rm grid}(\theta) \le S_{\psi,K}\norm{\theta-\theta_{\rm opt}}.
\end{equation}
Thus $\mathcal{E}_\infty^{\rm grid}(\theta_n)$ measures uniform compensation error across $X=[-2, 2]$ (visualized in the deterministic drift return maps), whereas time-series trajectories record the empirical $L^2(\mu_n)$ residual error
\begin{equation}
    \mathcal{E}_A(\theta_n) := \norm{g(\cdot, \theta_n) - d(\cdot)}_{L^2(\mu_n)} = \left(\frac{1}{N}\sum_{i=1}^N \bigl|g(x_n^{(i)}, \theta_n) - d(x_n^{(i)})\bigr|^2\right)^{1/2}.
    \label{eq:empirical-compensation-error}
\end{equation}

In the particle simulations, an empirical ensemble $\mu_n^N = \frac{1}{N}\sum_{i=1}^N \delta_{x_n^{(i)}}$ of $N$ trajectories evolving under $x_{n+1}^{(i)} = \mathfrak{r}_2(F_{\theta_n}(x_n^{(i)}) + \xi_n^{(i)})$ approximates $\mu_n$, and the parameter update is driven by the empirical one-step tracking-error estimator
\begin{equation}
    \widehat G_{\lambda,N}^{\rm tr}(\mu_n^N,\theta_n)
    := -\frac{1}{N}\sum_{i=1}^N \bigl(x_{n+1}^{(i)} - 0.8 x_n^{(i)}\bigr)\psi\bigl(x_n^{(i)}\bigr)
       + \lambda(\theta_n - \theta_c).
    \label{eq:empirical-tracking-gradient}
\end{equation}
This update differs from the deterministic distributional field $G(\mu_n,\theta_n)$ in \eqref{eq:duffing-learning-field-explicit} in two respects:
\begin{enumerate}[label=(\alph*)]
   \item \emph{Boundary reflection discrepancy:} Let $X \sim \mu$ and $\xi \sim \mathcal{N}(0,\sigma^2)$ be independent, and set $Z := F_\theta(X) + \xi$. When $Z \in [-2,2]$, we have $-(\mathfrak{r}_2(Z) - 0.8X) = g(X,\theta) - d(X) - \xi$, whose conditional expectation over $\xi$ coincides with $g(X,\theta) - d(X)$. When $\vert Z\vert > 2$, the folding map $\mathfrak{r}_2$ introduces a boundary bias supported on the Gaussian tail $\{\vert Z\vert > 2\}$:
   \[
       \Delta G_{\rm refl}(\mu,\theta)
       := \mathbb{E}\bigl[\widehat G_{\lambda,N}^{\rm tr}(\mu,\theta)\bigr] - G(\mu,\theta)
       = \mathbb{E}\bigl[(Z - \mathfrak{r}_2(Z))\psi(X)\bigr],
   \]
   which satisfies the analytical upper bound
   \[
       \norm{\Delta G_{\rm refl}(\mu,\theta)}
       \le B_{\rm refl}(\mu,\theta)
       := S_\psi\,\mathbb{E}\bigl[\vert Z - \mathfrak{r}_2(Z)\vert\,\mathbf{1}_{\{\vert Z\vert > 2\}}\bigr].
   \]
   For the cubic basis ($S_\psi = \sqrt{85} \approx 9.22$) at $\sigma = 0.18$, Gaussian-tail integration against the $M=1{,}200$ Ulam stationary measure at $(\mu_{\theta_{\rm opt}},\theta_{\rm opt})$ yields $\mathbb{P}(\vert Z\vert > 2) \approx 2.58\times 10^{-11}$, $B_{\rm refl} \approx 2.03\times 10^{-11}$, and $\norm{\Delta G_{\rm refl}} \approx 1.09\times 10^{-11}$. At the regularized equilibrium $(\mu_{\theta_*^\lambda},\theta_*^\lambda)$ ($\lambda=10^{-3}$), these values are $\mathbb{P}(\vert Z\vert > 2) \approx 1.78\times 10^{-9}$, $B_{\rm refl} \approx 2.07\times 10^{-9}$, and $\norm{\Delta G_{\rm refl}} \approx 1.45\times 10^{-9}$ (even-mode components vanish by parity). 
   
   Along the adaptation trajectory from $\theta_0 = [-0.25, 0.45, 0, 0]^{\mathsf T}$ ($\delta=0.10$, $\sigma=0.18$), transient state broadening briefly raises the Ulam values to peaks of $\mathbb{P}(\vert Z\vert > 2) \approx 4.50\times 10^{-4}$ and $\norm{\Delta G_{\rm refl}(\mu_n,\theta_n)} \approx 2.52\times 10^{-3}$ ($B_{\rm refl} \approx 3.30\times 10^{-3}$). In the $N=10^5$ ensemble simulation at the transient peak ($n=482$, reflected fraction $1.14\times 10^{-3}$), the samplewise bound $S_\psi\frac{1}{N}\sum_{i=1}^N \vert Z_n^{(i)} - \mathfrak{r}_2(Z_n^{(i)})\vert \approx 9.14\times 10^{-3}$ bounds the empirical reflection bias $\norm{\Delta G_{\rm emp}} \approx 7.07\times 10^{-3}$ (ratio $1.29$). For $n \ge 1{,}500$, $\mathbb{P}(\vert Z\vert > 2)$ drops below $10^{-6}$ and $\norm{\Delta G_{\rm refl}(\mu_n,\theta_n)} \le 1.8\times 10^{-5}$, with zero reflected particles over the final $5{,}000$ steps for $N=10^5$. 
   
   Comparing the $N=10^5$ trajectory under \eqref{eq:empirical-tracking-gradient} against an identical-noise trajectory driven by the exact unreflected residual $g(x_n^{(i)},\theta_n)-d(x_n^{(i)})$ shows that the discrepancy $\norm{\theta_n^{\rm trk} - \theta_n^{\rm theo}}$ decays from $2.53\times 10^{-2}$ ($n=500$) to $2.69\times 10^{-4}$ at $n=15{,}000$. By contrast, in the higher-noise misspecified experiment of Section~\ref{sec:linear-response-experiment} ($\sigma=0.24$, $\psi_{\rm mis}=[1,x,x^2]^{\mathsf T}$, $S_{\psi_{\rm mis}}=\sqrt{21}\approx 4.58$), using \eqref{eq:empirical-tracking-gradient} at $(\mu_{\theta_*},\theta_*)$ would incur $\mathbb{P}(\vert Z\vert > 2) \approx 2.48\times 10^{-3}$, $B_{\rm refl} \approx 7.39\times 10^{-3}$, and $\norm{\Delta G_{\rm refl}} \approx 3.04\times 10^{-3}$; to isolate $C_*R_*B_*$ without reflection bias, that simulation evaluates the exact residual field directly.
   
   \item \emph{Finite-ensemble fluctuations:} Averaging over $N < \infty$ particles and noise realizations $\xi_n^{(i)}$ introduces a zero-mean perturbation of order $O(N^{-1/2})$ at each step. Whereas deterministic discretizations of $\Phi_\delta$ decay to machine precision, finite-$N$ trajectories settle at a stochastic fluctuation floor of order $O(\sigma\sqrt{\delta/N})$ ($\approx 1.8\times 10^{-4}$ for $\sigma=0.18$, $\delta=0.10$, and $N=10^5$). The asymptotic prefactor depends on the stationary update covariance and the damping rates of individual linearized modes of $\mathcal{A}_\delta$. In the well-specified unregularized limit, the one-step noise covariance is $\frac{\sigma^2}{N}Q_{\rm opt}$, and the small-$\delta$ discrete Lyapunov balance gives the root-mean-square estimate $\sigma\sqrt{p\delta/(2N)} \approx 2.55\times 10^{-4}$ ($p=4$), consistent with the observed stationary plateau ($\approx 3\times 10^{-4}$).
\end{enumerate}

We compare two positive regularization weights with $\theta_c = 0$ under Proposition~\ref{prop:duffing-small-regularization}. The global dominance condition $\lambda > L_{\rm fb}(T)$ of Theorem~\ref{thm:duffing-main-results} is conservative due to the Gaussian tail in the uniform Doeblin bound \eqref{eq:duffing-epsilonT-explicit} and is not required for local stability:
\begin{itemize}
    \item $\lambda = 10^{-6}$: Evaluating $Q_{\rm opt} = \int_X \psi\psi^{\mathsf T}\dd\mu_{\theta_{\rm opt}}$ at $\theta_{\rm opt}$ on the $M=1{,}200$ Ulam grid gives $\lambda_{\min}(Q_{\rm opt}) \approx 4.06\times 10^{-3}$ (matching the stationary $N=10^5$ ensemble estimate $4.06\times 10^{-3}$, whereas $4.86\times 10^{-3}$ corresponds to the finite-horizon empirical decay ratio $-\text{slope}/\delta$ at $\delta=0.05$). The first-order shift $\norm{\lambda Q_{\rm opt}^{-1}\theta_{\rm opt}} \approx 3.55\times 10^{-5}$ lies below the finite-$N$ fluctuation floor ($\sim 3\times 10^{-4}$).
    \item $\lambda = 10^{-3}$: The equilibrium shift is resolved above the ensemble floor, with empirical bias $\norm{\hat{\theta}_*^\lambda - \theta_{\rm opt}} \approx 2.66\times 10^{-2}$, matching the Ulam root $\theta_*^\lambda \approx [0, 0.00699, 0, 0.12429]^{\mathsf T}$ of $\overline G_\lambda(\theta_*^\lambda)=0$ ($\norm{\theta_*^\lambda - \theta_{\rm opt}} \approx 2.66\times 10^{-2}$).
\end{itemize}

For $\lambda = 10^{-3}$, integrating the particle dynamics at $\delta=0.10$ for $25{,}000$ steps and averaging over the final $5{,}000$ steps yields $\hat{\theta}_*^\lambda \approx [-0.00001, 0.00699, 0.00004, 0.12435]^{\mathsf T}$. Using the Ulam Gram matrix $Q_{\rm opt}$, the first-order expansion of Proposition~\ref{prop:duffing-small-regularization} gives $\theta_{\rm pred} = \theta_{\rm opt} - \lambda Q_{\rm opt}^{-1}\theta_{\rm opt} \approx [0, 0.00926, 0, 0.11571]^{\mathsf T}$ ($\norm{\theta_{\rm pred} - \theta_{\rm opt}} \approx 3.55\times 10^{-2}$), leaving a nonlinear remainder $\norm{\theta_*^\lambda - \theta_{\rm pred}} \approx 8.9\times 10^{-3}$ consistent with $\lambda/\lambda_{\min}(Q_{\rm opt}) \approx 0.25$ (using the finite-horizon empirical estimate gives $[0, 0.00612, 0, 0.12015]^{\mathsf T}$, norm $3.05\times 10^{-2}$).

The uncompensated cubic residual at $\theta_*^\lambda$ broadens the stationary variance relative to $\mu_{\theta_{\rm opt}}$, increasing the minimum Gram eigenvalue to $\lambda_{\min}(H(\theta_*^\lambda)) \approx 4.65\times 10^{-3}$ on the Ulam grid ($6.97\times 10^{-3}$ empirically). Numerical evaluation of the reduced Jacobian $D\overline G_\lambda(\theta_*^\lambda) = H(\theta_*^\lambda) + \lambda\Id + C_*R_*B_*$ confirms that both $D\overline G_\lambda(\theta_*^\lambda)$ and its symmetric part are positive stable, with $\lambda_{\min}(D\overline G_\lambda(\theta_*^\lambda)) \approx 5.81\times 10^{-3} > 0$ on the Ulam grid and $\lambda_{\min}(\operatorname{sym} D\overline G_\lambda(\hat{\theta}_*^\lambda)) \approx 7.82\times 10^{-3}$ empirically, verifying the spectral stability hypothesis of Theorem~\ref{thm:local-stability} and Proposition~\ref{prop:duffing-small-regularization}. Throughout this section, evaluations on the finite Ulam grid represent standard numerical discretizations rather than rigorous computer-assisted enclosures.

The simulations use the parameters in Table~\ref{tab:parameters} (with initial conditions, parameter vectors, random seeds, and pseudocode in Section~\ref{sec:reproducibility-protocols}):
\begin{itemize}
    \item \emph{Ensemble approximation:} $\mu_n \in \Prob(X)$ is represented by an empirical measure of $N$ particles sharing $\theta_n$ (hence in $\Mmono$), and $G(\mu_n,\theta_n)$ is evaluated as the ensemble average.
    \item \emph{Burn-in:} In Phase~1 (Figures~\ref{fig:part2_action_error_recovery} and~\ref{fig:part2_escaping_trap_offline}), $\theta_{\rm off}$ is trained for $3{,}000$ steps under $u_{\rm base}(x) = -0.80x$, followed by $1{,}000$ steps with $\theta_{\rm off}$ frozen under the gain $-0.25$, during which $\mu_n$ relaxes to $\mu_{\rm deploy}$. In Figure~\ref{fig:part2_bifurcation_diagram}, each $\delta \in [1.50, 2.25]$ is integrated for $1{,}600$ steps, discarding $N_{\rm settle}=1{,}500$ transient steps and recording $N_{\rm record}=100$ states.
\end{itemize}

\begin{table}[htpb]
\centering
\footnotesize
\setlength{\tabcolsep}{3pt}
\begin{tabular}{lccccccc}
\toprule
\textbf{Experiment (Figure)} & $N$ & \textbf{Steps} & $\delta$ & $\sigma$ & $\lambda$ & \textbf{Baseline $u$} & \textbf{Burn-in} \\
\midrule
Fig.~\ref{fig:part1_nominal_convergence} (Autonomous Adapt.) & $10^5$ & 40k & 0.10 & 0.18 & $0, 10^{-6}, 10^{-3}$ & $-0.25x$ & None \\
Fig.~\ref{fig:part1_timescale_and_reg} (Regularized, Prop.~\ref{prop:duffing-small-regularization}) & $10^5$ & 15k & $0.005$ to $2.02$ & 0.18 & $10^{-3}, 5\times 10^{-3}$ & $-0.25x$ & None \\
Fig.~\ref{fig:linear_response_contribution} (Linear Response $C_*R_*B_*$) & $10^5$ & 800 & 0.05 & 0.24 & $0$ & $-0.25x$ & None \\
Fig.~\ref{fig:part2_identifiability_trap_pe} (Weak Excitation) & $10^5$ & 40k & 0.10 & $0.10, 0.18$ & $0$ & $-0.25x$ & None \\
Fig.~\ref{fig:part2_action_error_recovery} (Deployment Mismatch) & $10^5$ & 6k & 0.04 & $0.10 \to 0.18$ & $0$ & $-0.80x \to -0.25x$ & 3k + 1k \\
Fig.~\ref{fig:part2_escaping_trap_offline} (Matched Baseline) & $10^5$ & 5k & 0.015 & 0.10 & $0$ & $-0.25x$ & 4k + 1k \\
Fig.~\ref{fig:part2_bifurcation_diagram} (Stability Loss) & $10^2$ & 1.6k & $1.50$ to $2.25$ & 0.05 & $0$ & $-0.25x$ & 1.5k \\
Fig.~\ref{fig:part2_density_shattering} (Density Oscillations) & $10^4$ & 150 & $1.85, 2.22$ & 0.05 & $0$ & $-0.25x$ & None \\
Fig.~\ref{fig:mlp_convergence} (MLP Compensator) & $10^4$ & 1.5k & 0.05 & 0.10 & $10^{-5}$ & $-0.25x$ & None \\
\bottomrule
\end{tabular}
\caption{Simulation and system parameters for the numerical experiments.}
\label{tab:parameters}
\end{table}

\subsection{Part 1: Autonomous Online Adaptation}
\label{sec:part1-illustration}

Figures~\ref{fig:part1_nominal_convergence} and \ref{fig:part1_timescale_and_reg} use $\sigma=0.18$ and $\theta_0=[-0.25,0.45,0,0]^{\mathsf T}$. By Proposition~\ref{prop:duffing-doeblin}, $P_\theta$ contracts zero-mass measures uniformly in total variation. At $\sigma=0.18$, the stationary measure $\mu_*$ spreads across $[-2,2]$ (Figure~\ref{fig:part2_identifiability_trap_pe}, left), so that $H_*$ is well-conditioned with $\lambda_{\min}(H_*) \approx 4.37\times 10^{-3}$. By \eqref{eq:duffing-slow-multipliers}, the slowest mode relaxes on the time scale $1/(\delta\lambda_{\min}(H_*)) \approx 2.3\times 10^{3}$ steps at $\delta=0.10$.

In the unregularized realizable case ($\lambda=0$), $g(\cdot,\theta_{\rm opt}) = d(\cdot)$ implies $\nabla_\theta\ell_0(\cdot,\theta_{\rm opt}) \equiv 0$ and hence $C_*=0$. Since $D_\theta\mu_\theta$ is bounded (Proposition~\ref{prop:duffing-P-differentiability}), the linear-response contribution vanishes at equilibrium, and
\begin{equation}
    D\overline G(\theta_{\rm opt}) = \mathcal{J}_* = H_* = Q_{\rm opt} \succ 0.
\end{equation}
With $\Gammahat=\Id$, Remark~\ref{cor:LR-stability} gives local exponential stability. In Figure~\ref{fig:part1_nominal_convergence} (left), $\norm{\theta_n - \theta_{\rm opt}}$ for $\lambda=0$ and $\lambda=10^{-6}$ decays exponentially to the numerical floor ($\approx 3\times 10^{-4}$). For $\lambda=10^{-3}$, $\norm{\theta_n - \theta_{\rm opt}}$ settles at the bias plateau $\norm{\theta_*^\lambda - \theta_{\rm opt}} \approx 2.66\times 10^{-2}$, whereas $\norm{\theta_n - \hat{\theta}_*^\lambda}$ decays exponentially to the numerical floor, in agreement with Proposition~\ref{prop:duffing-small-regularization} and Theorem~\ref{thm:local-stability}.

\begin{figure}[htpb]
\centering
\includegraphics[width=\textwidth]{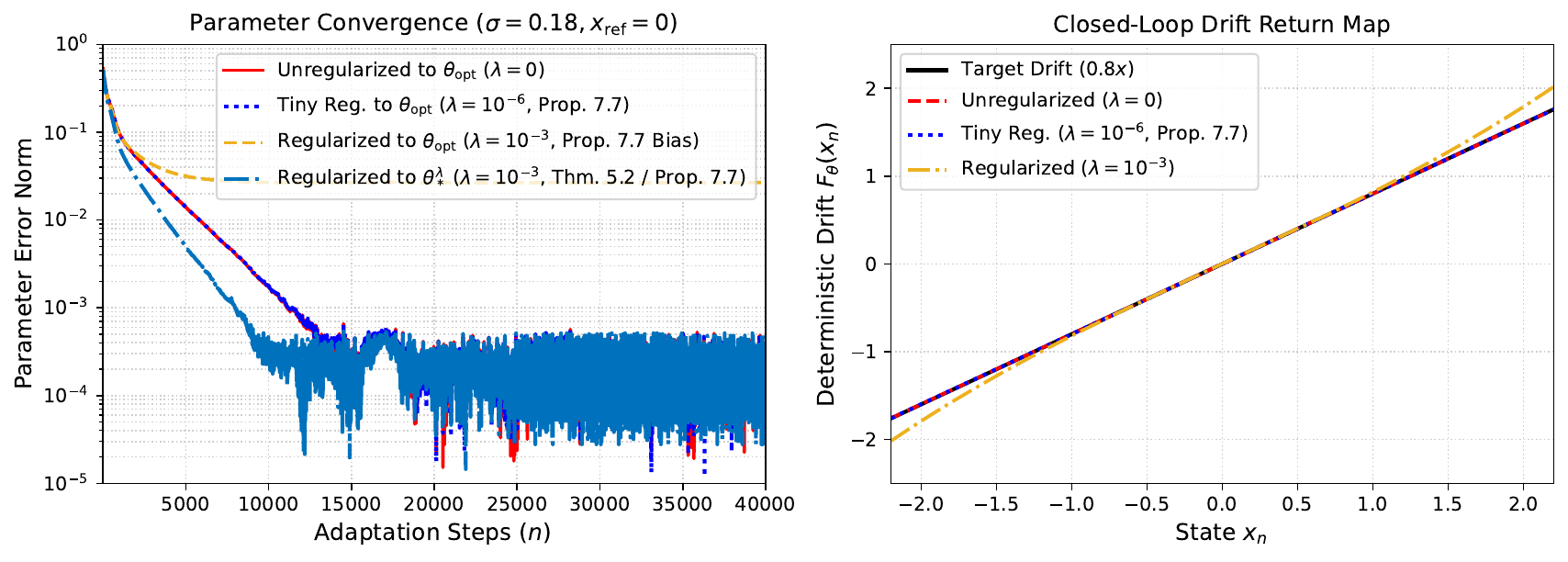}
\caption{Convergence of autonomous adaptation ($\sigma=0.18$, $\delta=0.10$). (Left) Parameter error. The trajectories for $\lambda=0$ (solid red) and $\lambda=10^{-6}$ (dotted blue) decay to the numerical floor. For $\lambda=10^{-3}$, the distance to $\theta_{\rm opt}$ (dashed yellow) settles at $\norm{\theta_*^\lambda - \theta_{\rm opt}} \approx 2.66\times 10^{-2}$, whereas the distance to $\theta_*^\lambda$ (dash-dotted blue) decays to the numerical floor. (Right) Return map of the deterministic drift. The unregularized (dashed red) and $\lambda=10^{-6}$ (dotted blue) fits match $x_{\rm target}(x) = 0.8x$ (solid black), while the $\lambda=10^{-3}$ fit (dash-dotted gold) retains a cubic residual.}
\label{fig:part1_nominal_convergence}
\end{figure}

For $\lambda = 10^{-3}$ ($\Gammahat=\Id$), Theorem~\ref{thm:local-stability} and Proposition~\ref{prop:effective-jacobian-identification} yield local exponential convergence to $(\mu_{\theta^\lambda_*},\theta^\lambda_*)$ for $0<\delta<\delta_0$, with slow multipliers governed by \eqref{eq:duffing-slow-multipliers} and relaxation time of order $\delta^{-1}$. Both the shifted stationary measure $\mu_{\theta_*^\lambda}$ and the non-vanishing linear-response term $C_*R_*B_*$ enter the reduced Jacobian $D\overline G_\lambda(\theta_*^\lambda) = H(\theta_*^\lambda) + \lambda\Id + C_*R_*B_*$, whose eigenvalues on the $M=1{,}200$ Ulam grid are $\{0.00581,\, 0.01860,\, 0.10227,\, 1.01054\}$. Because $C_* \neq 0$, the coupled linearization $\mathcal{A}_\delta$ in \eqref{eq:block-matrix} is not upper block-triangular; thus the small-$\delta$ expansion \eqref{eq:duffing-slow-multipliers} and $\lambda_{\max}(D\overline G_\lambda(\theta_*^\lambda)) \approx 1.0105$ do not alone determine the stability threshold at finite $\delta$ (unlike the $C_*=0$ case of Section~\ref{sec:loss-of-stability}, where $\delta_c = 2/\lambda_{\max}(H_*)$ is exact).

To find the finite-$\delta$ stability boundary of the discretized pair dynamics, we restrict $P_*$ and $B_*$ to the $(M-1)$-dimensional zero-mass subspace $V(X)$ on the $M=1{,}200$ Ulam grid (where $r(P_*\vert_{V(X)}) \approx 0.8002$) and compute the spectrum of the $(M+3)\times(M+3)$ block matrix
\[
    \mathcal{A}_\delta = \begin{pmatrix} P_*\vert_{V(X)} & B_* \\ -\delta C_* & \Id - \delta\bigl(H(\theta_*^\lambda) + \lambda\Id\bigr) \end{pmatrix}.
\]
Its spectral radius satisfies $r(\mathcal{A}_\delta) \in \{0.99997,\, 0.99988,\, 0.99971,\, 0.99954\} < 1$ across the stable learning rates $\delta \in \{0.005, 0.02, 0.05, 0.08\}$. At $\delta=2.02$, $r(\mathcal{A}_{2.02}) \approx 1.0385 > 1$, crossing unity through a real eigenvalue $-1.00000006$ at the critical threshold $\delta_c(\mathcal{A}_\delta) \approx 1.981814$. Because $\det(\mathcal{A}_\delta + \Id) = 0$ at $z=-1$ reduces via the Schur complement on $V(X)$ to $\det\bigl(2\Id - \delta M_{\rm flip}\bigr) = 0$ with
\[
    M_{\rm flip} := K - C_*\vert_{V(X)}\bigl(\Id + P_*\vert_{V(X)}\bigr)^{-1} B_*\vert_{V(X)},
    \qquad
    K := H(\theta_*^\lambda) + \lambda\Id,
\]
where $\lambda_{\max}(M_{\rm flip}) \approx 1.009176$, this threshold agrees to seven significant digits with the exact flip formula $\delta_{c,{\rm flip}} = 2/\lambda_{\max}(M_{\rm flip}) \approx 1.981814$, lying within $0.01\%$ of $2/\lambda_{\max}(K) \approx 1.9817$ and $0.14\%$ of the adiabatic estimate $2/\lambda_{\max}(D\overline G_\lambda(\theta_*^\lambda)) \approx 1.9791$.

Figure~\ref{fig:part1_timescale_and_reg} (left) plots the empirical trajectories across $\delta \in \{0.005, 0.02, 0.05, 0.08, 2.02\}$, and Figure~\ref{fig:part1_timescale_and_reg} (right) compares convergence at $\delta=0.05$. The $\lambda=10^{-3}$ trajectory converges to $\theta_*^\lambda$ faster than the unregularized trajectory ($\lambda=0$) converges to $\theta_{\rm opt}$. This reflects the larger minimum eigenvalue of $D\overline G_\lambda(\theta_*^\lambda)$ relative to $D\overline G(\theta_{\rm opt}) = Q_{\rm opt}$: on the Ulam grid, $\lambda_{\min}$ increases from $4.06\times 10^{-3}$ to $5.81\times 10^{-3}$ (per-step decay rate $\delta\lambda_{\min}$ increasing from $2.03\times 10^{-4}$ to $2.90\times 10^{-4}$), and in the finite-horizon empirical estimate from $4.86\times 10^{-3}$ to $7.82\times 10^{-3}$. This acceleration is driven jointly by the broadened stationary distribution $\mu_{\theta_*^\lambda}$ in $H(\theta_*^\lambda)$, the $+\lambda\Id$ shift, and the linear-response term $C_*R_*B_*$.

\begin{figure}[htpb]
\centering
\includegraphics[width=\textwidth]{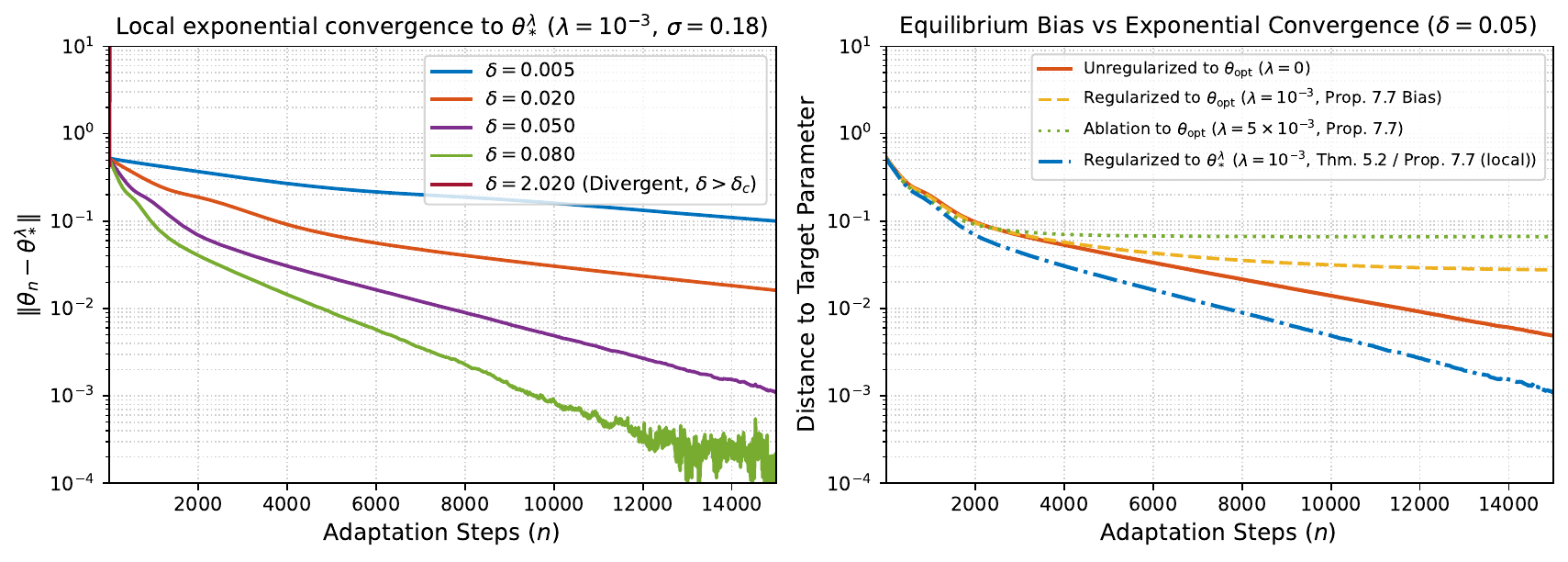}
\caption{Local exponential convergence to $\theta_*^\lambda$ ($\sigma=0.18$). (Left) Regularized error $\norm{\theta_n - \theta_*^\lambda}$ for $\lambda = 10^{-3}$. For $\delta < \delta_c$, trajectories converge with decay slopes proportional to $\delta$; for $\delta=2.02 > \delta_c$, the trajectory diverges. (Right) Convergence at $\delta=0.05$. Measured relative to $\theta_{\rm opt}$, trajectories settle at $O(\lambda)$ bias plateaus ($\approx 2.66\times 10^{-2}$ for $\lambda=10^{-3}$, dashed yellow; $\approx 6.64\times 10^{-2}$ for $\lambda=5\times 10^{-3}$, dotted green). Measured relative to $\theta_*^\lambda$, the $\lambda=10^{-3}$ trajectory (dash-dotted blue) converges exponentially.}
\label{fig:part1_timescale_and_reg}
\end{figure}

\subsubsection{Isolating the Statistical Linear-Response Contribution \texorpdfstring{$C_*R_*B_*$}{C*R*B*}}
\label{sec:linear-response-experiment}

In the realizable cubic setting with $\lambda=0$, exact cancellation $g(\cdot,\theta_{\rm opt})=d$ gives $\nabla_\theta\ell_0(\cdot,\theta_{\rm opt})\equiv 0$ and $C_*=0$. To isolate the statistical feedback term $C_*R_*B_*$ in Proposition~\ref{prop:effective-jacobian-identification}, we consider an under-parameterized quadratic compensator
\begin{equation}
    \psi_{\rm mis}(x) = [1,\, x,\, x^2]^{\mathsf T} \in \R^3, \qquad \lambda = 0,
\end{equation}
controlling the same reflected cubic plant ($d(x)=0.15x^3$) at $\sigma=0.24$. Because $x^3 \notin \operatorname{span}(\psi_{\rm mis})$, the equilibrium misfit $g(\cdot,\theta_*)-d$ is nonzero at the self-consistent root $\theta_* \approx [0,\, 0.1121,\, 0]^{\mathsf T}$ of $\overline G(\theta_*)=0$. Parity symmetry of $\mu_{\theta_*}$ decouples the even and odd coordinates, giving $\theta_{*,1} = 0.15\,m_4(\theta_*)/m_2(\theta_*)$.

Along the odd parameter coordinate $\theta_1$ (the linear feedback gain in $\theta=[\theta_0,\theta_1,\theta_2]^{\mathsf T}$, corresponding to the $(2,2)$ matrix entry), the frozen-measure Hessian and the full reduced Jacobian satisfy
\begin{equation}
    (H_*)_{22} = m_2(\theta_*),
    \qquad
    \bigl(D\overline G(\theta_*)\bigr)_{22}
    = (H_*)_{22} + (C_*R_*B_*)_{22}
    = m_2(\theta_*) + \theta_{*,1}m_2'(\theta_*) - 0.15\,m_4'(\theta_*),
\end{equation}
where primes denote derivatives with respect to $\theta_1$ at $\theta_*$. Substituting $\theta_{*,1} = 0.15\,m_4(\theta_*)/m_2(\theta_*)$ expresses the linear-response contribution in terms of relative moment sensitivities:
\begin{equation}
    (C_*R_*B_*)_{22}
    = 0.15\,m_4(\theta_*)
      \left(
        \frac{m_2'(\theta_*)}{m_2(\theta_*)}
        - \frac{m_4'(\theta_*)}{m_4(\theta_*)}
      \right).
    \label{eq:CRB22-relative-moments}
\end{equation}

Increasing $\theta_1$ strengthens the linear restoring drift $0.8-\theta_1$, contracting $\mu_\theta$ toward the origin so that $m_2'(\theta_*) < 0$ and $m_4'(\theta_*) < 0$. By \eqref{eq:CRB22-relative-moments}, $(C_*R_*B_*)_{22} > 0$ holds if and only if the fourth moment contracts faster than the second moment in relative terms ($\lvert m_4'(\theta_*)\rvert/m_4(\theta_*) > \lvert m_2'(\theta_*)\rvert/m_2(\theta_*)$). Under a Gaussian closure $m_4 = 3 m_2^2$, logarithmic differentiation gives $m_4'/m_4 = 2 m_2'/m_2$ and hence $(C_*R_*B_*)_{22} = -0.15\,m_4\,(m_2'/m_2) > 0$. Although the true invariant measure deviates from Gaussianity due to the cubic drift and boundary reflections—exhibiting an equilibrium kurtosis of $m_4/m_2^2 \approx 4.98$ versus $3.0$—it preserves the positive sign of the feedback term. Moreover, the misspecified deterministic drift $F_{\theta_*}(x) = (0.8 - \theta_{*,1})x + 0.15 x^3$ has unstable fixed points ($F_{\theta_*}'(x_u) \approx 1.62 > 1$) at $x_u = \pm\sqrt{(0.20 + \theta_{*,1})/0.15} \approx \pm 1.44$, beyond which deterministic trajectories are repelled toward the reflecting boundaries at $\pm 2$.

To evaluate these operators numerically, we partition $X=[-2,2]$ into $M$ uniform cells of width $\Delta x = 4/M$ centered at midpoints $\{x_i\}_{i=1}^M$. The discretized transfer matrix $P_\theta \in \R^{M\times M}$ is assembled by evaluating the folded Gaussian transition density between cell midpoints across the reflection branches of $\mathfrak{r}_2$ and normalizing columns to unit mass. The stationary distribution $\mu_\theta \in \R^M$ is obtained by solving $(\Id - P_\theta)\mu_\theta = 0$ with $\mathbf{1}^{\mathsf T}\mu_\theta = 1$, and $\theta_*$ is computed by Newton iteration on $\overline G(\theta)=0$. At $\theta_*$, the sensitivity columns $B_* e_j = (D_\theta P_{\theta_*}[e_j])\mu_*$ ($j=1,2,3$) are evaluated by differentiating the folded Gaussian kernel analytically along $\partial_{\theta_j} F_\theta = -\psi_{\rm mis,j}$ and subtracting $(\mathbf{1}^{\mathsf T} B_* e_j)\mu_*$ so that $B_* e_j \in V(X) = \{\eta \in \R^M : \mathbf{1}^{\mathsf T}\eta = 0\}$. The resolvent action $R_* B_* = (\Id - P_*)^{-1}\vert_{V(X)} B_*$ is computed from the nonsingular rank-one stabilized system
\[
    \bigl(\Id - P_* + \mu_* \mathbf{1}^{\mathsf T}\bigr)(R_* B_*) = B_*
\]
followed by projection onto $V(X)$, and $C_* R_* B_*$ is obtained by integrating $(g(x_i,\theta_*) - d(x_i))\psi_{\rm mis}(x_i)$ against $R_* B_*$.

At $M=1{,}200$, the equilibrium moments and their derivatives along $\theta_1$ (via $D_{\theta_1}\mu_{\theta_*} = R_* B_* e_2$) are $m_2(\theta_*) \approx 0.1500$, $m_4(\theta_*) \approx 0.1121$, $m_2'(\theta_*) \approx -1.1313$, and $m_4'(\theta_*) \approx -2.2809$. This gives relative sensitivities
\[
    \frac{m_2'(\theta_*)}{m_2(\theta_*)} \approx -7.5398,
    \qquad
    \frac{m_4'(\theta_*)}{m_4(\theta_*)} \approx -20.3384,
\]
confirming that $m_4$ contracts roughly $2.7$ times faster in relative terms than $m_2$ and yielding
\[
    (H_*)_{22} \approx 0.1500, \qquad (C_*R_*B_*)_{22} \approx +0.2153 \quad (+143.5\%), \qquad \bigl(D\overline G(\theta_*)\bigr)_{22} \approx 0.3653.
\]
On the even subspace $\operatorname{span}\{1,x^2\}$, the non-symmetric feedback block $(C_*R_*B_*)_{\rm even} \approx \bigl(\begin{smallmatrix} -0.0052 & 0.1086 \\ 0.5413 & 0.3585 \end{smallmatrix}\bigr)$ has real eigenvalues $\{-0.1265,\, +0.4797\}$, increasing the smallest eigenvalue of the even block from $\lambda_{\min}(H_{\rm even}) \approx 0.0875$ to $\lambda_{\min}(H_{\rm even} + (C_*R_*B_*)_{\rm even}) \approx 0.2352$ ($+168.9\%$). Across grids with $M \in \{600,\, 1{,}200,\, 2{,}400\}$, all reported equilibrium coordinates, moments, and eigenvalues vary by less than $10^{-5}$. 

As an independent check at $M=1{,}200$, approximating $D\overline G(\theta_*)$ by central finite differences $[\overline G(\theta_*+\varepsilon e_j)-\overline G(\theta_*-\varepsilon e_j)]/(2\varepsilon)$ with step $\varepsilon=10^{-6}$ (re-solving $(\Id - P_{\theta_*\pm\varepsilon e_j})\mu_{\theta_*\pm\varepsilon e_j}=0$ at each perturbed parameter) matches the analytical linear-response sum $H_* + C_*R_*B_*$ to $\norm{D_{\rm FD}\overline G(\theta_*) - (H_*+C_*R_*B_*)} \approx 1.19\times 10^{-8}$ (relative error $9.1\times 10^{-9}$), whereas omitting $C_*R_*B_*$ leaves $\norm{D_{\rm FD}\overline G(\theta_*) - H_*} \approx 0.6517$ (for the regularized cubic system of Section~\ref{sec:part1-illustration} at $\lambda=10^{-3}$ and $\sigma=0.18$, the corresponding errors are $2.04\times 10^{-9}$ with $C_*R_*B_*$ versus $7.28\times 10^{-3}$ without it). Figure~\ref{fig:linear_response_contribution} shows that the local exponential decay of $\norm{\theta_n - \theta_*}$ under $\Phi_\delta$ ($\delta=0.05$, with the finite-$N$ ensemble driven directly by $\frac{1}{N}\sum_{i=1}^N (g(x_n^{(i)},\theta_n)-d(x_n^{(i)}))\psi_{\rm mis}(x_n^{(i)})$ to exclude boundary-reflection bias at $\sigma=0.24$) follows the full linear-response rate $\lambda_j(H_* + C_*R_*B_*)$, with log-linear fitted slopes matching $\log(1-\delta\lambda_j(H_*+C_*R_*B_*))$ ($\approx -0.0184$ odd, $-0.0118$ even), whereas the frozen-distribution approximation $\lambda_j(H_*)$ ($\approx -0.0075$ odd, $-0.0044$ even) systematically underestimates the relaxation rate.

\begin{figure}[htpb]
\centering
\includegraphics[width=\textwidth]{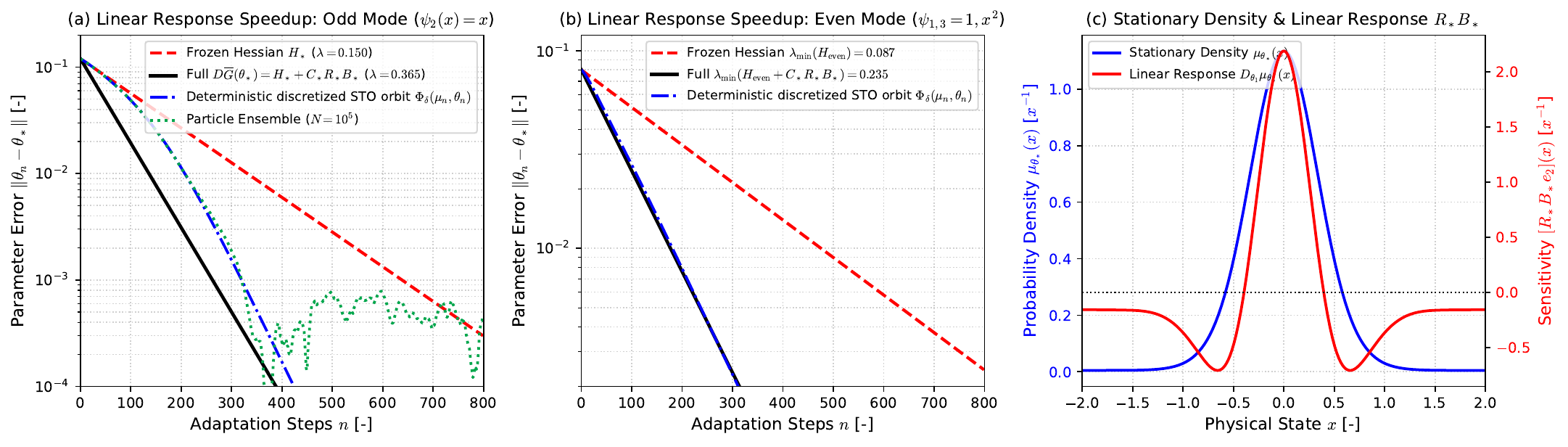}
\caption{Effect of the statistical linear-response term $C_*R_*B_*$ for the misspecified quadratic compensator $\psi_{\rm mis}(x)=[1,x,x^2]^{\mathsf T}$ ($\sigma=0.24$, $\delta=0.05$, $\lambda=0$; all quantities dimensionless). \textbf{(a)} Relaxation of $\norm{\theta_n - \theta_*}$ versus step $n$ for an odd-mode perturbation ($\theta_1$). The empirical finite-particle trajectory ($N=10^5$, dotted green, evaluated with the exact residual $g-d$) and deterministic discretized STO orbit (dash-dotted blue) match the full Jacobian decay parameter $(H_*+C_*R_*B_*)_{22} \approx 0.365$ (solid black; $0.3653$ unrounded), whereas the frozen-measure parameter $(H_*)_{22} \approx 0.150$ (dashed red; $0.1500$ unrounded) underpredicts the decay rate. \textbf{(b)} Relaxation of $\norm{\theta_n - \theta_*}$ along the slowest even mode in $\mathrm{span}\{1,x^2\}$, comparing the deterministic discretized STO orbit (dash-dotted blue) against $\lambda_{\min}(H_{\rm even}+(C_*R_*B_*)_{\rm even}) \approx 0.235$ (solid black; $0.2352$ unrounded) and $\lambda_{\min}(H_{\rm even}) \approx 0.087$ (dashed red; $0.0875$ unrounded). \textbf{(c)} Stationary density of $\mu_{\theta_*}$ (left axis, blue) and its linear response $D_{\theta_1}\mu_{\theta_*}=R_*B_*e_2$ (right axis, red), where $e_2=(0,1,0)^{\mathsf T}$. Perturbing $\theta_1\mapsto\theta_1+\varepsilon$ changes the physical drift by $-\varepsilon x$; relative to the dotted zero-response line, the response is positive near the origin and negative farther out, reflecting stronger concentration for $\varepsilon>0$.}
\label{fig:linear_response_contribution}
\end{figure}

\subsection{Part 2: Ill-Conditioning, Deployment Mismatch, and Large Learning Rates}
\label{sec:part2-practice}

In robotic posture and exoskeleton control \cite{lippi2018prediction,lippi2021icinco}, a compensator is often pre-trained offline under a high-gain stabilizing baseline controller and subsequently deployed online under a more compliant gain, over finite horizons and with finite learning rates \cite{LippiActionError}. This subsection examines three phenomena arising in such regimes:
\begin{enumerate}[label=(\roman*)]
    \item \emph{Finite-horizon ill-conditioning vs.\ asymptotic stability (Section~\ref{sec:weak-excitation}):} While Proposition~\ref{prop:duffing-doeblin} and Remark~\ref{rem:duffing-realizable-case} guarantee $H_*\succ 0$ and asymptotic exponential convergence for any $\sigma>0$, small noise ($\sigma=0.10$) causes $\lambda_{\min}(H_*)\propto \sigma^6$ to nearly vanish, illustrating the classical trade-off between regulation and persistent excitation \cite{AnnaswamyFradkov2021,Annaswamy2023,LippiActionError}.
    \item \emph{Transient recovery from deployment mismatch (Section~\ref{sec:deployment-invariance-failure}):} Whereas Theorem~\ref{thm:local-stability} applies in a local neighborhood $\mathcal N_\delta$ of $(\mu_*,\theta_*)$, deploying an offline-trained controller $\theta_{\rm off}$ with a nonzero invariance residual $\mathfrak R_A = P_{\theta_{\rm off}}\mu_{\rm off}-\mu_{\rm off}\neq 0$ induces an immediate distributional shift ($\mu_{\rm off}\to\mu_{\rm deploy}$) \cite{LippiActionError,sinha2023systemlevelviewoutofdistributiondata,kang2022lyapunovdensitymodelsconstraining} before online adaptation restores convergence to $(\mu_*,\theta_*)$.
    \item \emph{Large learning rates and post-instability dynamics (Section~\ref{sec:loss-of-stability}):} Increasing $\delta$ across and beyond the linear stability threshold $\delta_c = 2/\lambda_{\max}(H_*)$ produces damped alternating density oscillations, a nonlinear period-$2$ saturation regime, and a secondary loss of physical contractivity.
\end{enumerate}

\subsubsection{Weak Excitation and Ill-Conditioning}
\label{sec:weak-excitation}
Gaussian noise of intensity $\sigma > 0$ ensures that $\mu_*$ has full support on $X=[-2,2]$ and that the polynomial features remain linearly independent ($H_* \succ 0$). Reducing $\sigma$ to $0.10$, however, lowers $\lambda_{\min}(H_*)$ to $\approx 1.29\times 10^{-4}$ via \eqref{eq:duffing-sigma6-law}. Because $\mu_*$ concentrates tightly near the origin (Figure~\ref{fig:part2_identifiability_trap_pe}, left), higher-order monomials become nearly collinear in $L^2(\mu_*)$, rendering $H_*$ ill-conditioned. The slowest multiplier approaches $1$ (with relaxation time $1/(\delta\lambda_{\min}(H_*)) \approx 7.8\times 10^4$ steps at $\delta=0.10$), slowing parameter convergence over finite horizons. At $n=40{,}000$ steps (Figure~\ref{fig:part2_identifiability_trap_pe}, right), the compensator matches the target drift near the origin but retains large errors in regions of low stationary probability.

\begin{figure}[htpb]
\centering
\includegraphics[width=\textwidth]{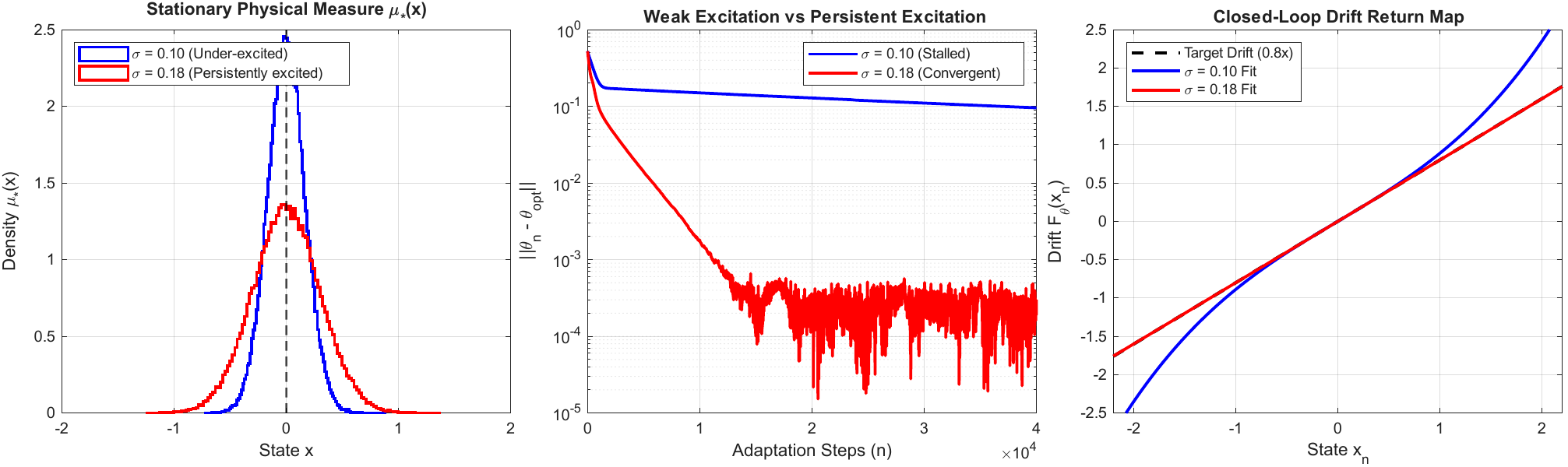}
\caption{Weak excitation and ill-conditioning ($\delta=0.10$). (Left) Stationary distribution $\mu_*(x)$. At $\sigma=0.10$, the density concentrates near $x=0$. (Middle) Parameter error $\norm{\theta_n - \theta_{\rm opt}}$. At $\sigma=0.10$, $\lambda_{\min}(H_*) \approx 1.29\times 10^{-4}$, stalling convergence. (Right) Return map of the deterministic drift. For $\sigma=0.10$, the fit matches $0.8x$ near the origin but diverges outside $[-0.5, 0.5]$.}
\label{fig:part2_identifiability_trap_pe}
\end{figure}

\subsubsection{Failure of Invariance at Deployment}
\label{sec:deployment-invariance-failure}
In the two-phase deployment protocol ($\lambda=0$), offline pre-training is performed under the conservative baseline controller $u_{\rm base}^{\rm off}(x) = -0.80x$ at noise intensity $\sigma_{\rm off}=0.10$, while the parameter update evaluates the tracking-error residual $-(x_{n+1}-0.8x_n)$ of \eqref{eq:empirical-tracking-gradient} relative to the nominal reference $0.8x_n$. Since the pre-reflection dynamics under $u_{\rm base}^{\rm off}(x)=-0.80x$ is
\[
    x_{n+1} = 0.25x_n + d(x_n) - g(x_n,\theta) + \xi_n,
\]
neglecting boundary reflection and using $\mathbb{E}[\xi_n]=0$ gives $\mathbb{E}[-(x_{n+1}-0.8x_n)\mid x_n=x] = g(x,\theta) - d(x) + 0.55x$. Up to the negligible reflection tail, the offline update follows the gradient of the shifted quadratic loss
\[
    \ell_{\rm off}(x,\theta) = \frac{1}{2}\bigl(g(x,\theta) - d(x) + 0.55x\bigr)^2,
\]
with offline learning field
\begin{equation}
    G_{\rm off}(\mu,\theta) := \int_X \bigl(g(x,\theta) - d(x) + 0.55x\bigr)\psi(x)\,\dd\mu(x).
    \label{eq:offline-learning-field}
\end{equation}
Offline identification under $P_{\rm off,\theta}$ thus seeks $\theta_{\rm off}$ satisfying
\begin{equation}
    G_{\rm off}(\mu_{\rm off}, \theta_{\rm off}) \approx 0.
\end{equation}
By Proposition~\ref{prop:all-fixed-points}, equilibrium after deployment under the relaxed baseline $u_{\rm base}^{\rm dep}(x) = -0.25x$ requires instead
\begin{equation}
    P_{\rm dep,\theta_*}\mu_* = \mu_* \qquad \text{and} \qquad G_{\rm dep}(\mu_*,\theta_*) = 0,
\end{equation}
where $P_{\rm dep,\theta} \equiv P_\theta$ and $G_{\rm dep}(\mu,\theta) \equiv G(\mu,\theta) = \int_X (g(x,\theta)-d(x))\psi(x)\,\dd\mu(x)$. Relaxing the baseline gain from $-0.80x$ to $-0.25x$ therefore modifies both the physical kernel ($P_{\rm off,\theta}\mapsto P_{\rm dep,\theta}$) and the effective learning target ($d(x)-0.55x \mapsto d(x)$). In the unregularized realizable cubic case, exact offline identification ($G_{\rm off}(\mu_{\rm off},\theta_{\rm off}^*)=0$, neglecting reflection) yields the ideal offline compensator
\[
    g_{\rm off}^*(x) = 0.15x^3 - 0.55x, \qquad \theta_{\rm off}^* = [0,\, -0.55,\, 0,\, 0.15]^{\mathsf T}.
\]
Upon deployment under $u_{\rm base}^{\rm dep}(x)=-0.25x$, the deterministic map before reflection becomes
\begin{equation}
    F_{\rm dep}(x) := 0.8x + 0.15x^3 - g(x,\theta),
    \label{eq:F-dep-general}
\end{equation}
which for $g_{\rm off}^*$ reduces to $F_{\rm dep}(x) = 1.35x$. Thus, even exact offline identification against the fixed reference residual $x_{n+1}-0.8x_n$ destabilizes the origin once the baseline gain is relaxed.

When the deployed operator $P_{\rm dep,\theta_{\rm off}} \equiv P_{\theta_{\rm off}}$ does not preserve $\mu_{\rm off}$, physical invariance fails with residual \cite{LippiActionError}
\begin{equation}
    \mathfrak{R}_A := P_{\rm dep,\theta_{\rm off}}\mu_{\rm off} - \mu_{\rm off} \neq 0.
    \label{eq:action-residual-ref}
\end{equation}
Note the distinction between the signed invariance residual $\mathfrak{R}_A \in V(X)$, which vanishes if and only if $\mu_{\rm off}$ is invariant under $P_{\rm dep,\theta_{\rm off}}$ and drives the transient shift $\mu_{\rm off}\to\mu_{\rm deploy}$, and the scalar compensation error $\mathcal{E}_A(\theta_n)=\norm{g(\cdot,\theta_n)-d}_{L^2(\mu_n)}$ in \eqref{eq:empirical-compensation-error} plotted in Figure~\ref{fig:part2_action_error_recovery} (middle).

In the finite-horizon simulation ($3{,}000$ pre-training steps from $\theta_{\rm off,0}=[0,-0.45,0,0.28]^{\mathsf T}$ and $x_0^{(i)}=0$ at $\delta_{\rm off}=0.002$, $\sigma_{\rm off}=0.10$, and $\lambda=0$), pre-training under $u_{\rm base}^{\rm off}(x) = -0.80x$ confines $\mu_{\rm off}$ near the origin ($\operatorname{std}_{\mu_{\rm off}}(X)\approx 0.141$, $\lambda_{\min}(H(\mu_{\rm off}))\approx 4.49\times 10^{-5}$; Figure~\ref{fig:part2_action_error_recovery}, left). The linear coefficient adapts partially toward $-0.55$, whereas the weakly excited cubic coefficient remains near its initialization, yielding
\begin{equation}
    \theta_{\rm off} = [\theta_{0,\rm off},\, \theta_{1,\rm off},\, \theta_{2,\rm off},\, \theta_{3,\rm off}]^{\mathsf T} \approx [1.03\times 10^{-5},\, -0.46188,\, 4.62\times 10^{-7},\, 0.27928]^{\mathsf T}.
    \label{eq:theta-off-numerical}
\end{equation}
Substituting \eqref{eq:theta-off-numerical} into \eqref{eq:F-dep-general} gives the deployed pre-reflection map
\[
\begin{aligned}
    F_{\rm dep}(x) &= -\theta_{0,\rm off} + (0.8 - \theta_{1,\rm off})x - \theta_{2,\rm off}x^2 + (0.15 - \theta_{3,\rm off})x^3 \\
    &\approx -1.03\times 10^{-5} + 1.26188\,x - 4.62\times 10^{-7}\,x^2 - 0.12928\,x^3,
\end{aligned}
\]
with slope at the origin
\begin{equation}
    F_{\rm dep}'(0) = 0.8 - \theta_{1,\rm off} \approx 1.26188 > 1
    \label{eq:F-dep-slope}
\end{equation}
(compared to $1.35$ for $\theta_{\rm off}^*$). Since $\theta_{0,\rm off}\approx 1.03\times 10^{-5}$ is negligible, the central fixed point of $F_{\rm dep}(x)=x$ lies at $x_{\rm fp,0} \approx +3.92\times 10^{-5} \approx 0$ with multiplier $F_{\rm dep}'(x_{\rm fp,0}) \approx 1.26188 > 1$ (Figure~\ref{fig:part2_action_error_recovery}, right). Meanwhile, the cubic over-compensation $(0.15-\theta_{3,\rm off}) \approx -0.12928 < 0$ creates two locally stable outer fixed points at $x_{{\rm fp},\pm} \approx \pm 1.4233$ ($F_{\rm dep}'(x_{{\rm fp},\pm}) \approx 0.4762 \in (-1,1)$), splitting $\mu_{\rm off}$ into the bimodal distribution $\mu_{\rm deploy}$ peaked near $\pm 1.42$ (Figure~\ref{fig:part2_action_error_recovery}, left) and reflecting $\mathfrak{R}_A \neq 0$.

To quantify $\mathfrak{R}_A$ without comparing an atomic empirical measure directly against a continuous density in total variation, we evaluate the discretized invariance residual on the $M$-cell Ulam grid ($h = 4/M$):
\begin{equation}
    r_{A,h} := P_{{\rm dep},h}(\theta_{\rm off})\mu_{{\rm off},h} - \mu_{{\rm off},h},
    \qquad
    \norm{r_{A,h}}_{\ell^1} = \sum_{i=1}^M \bigl|(r_{A,h})_i\bigr|,
    \label{eq:discretized-rA}
\end{equation}
where $\mu_{{\rm off},h} \in \R^M$ solves $(\Id - P_{{\rm off},h}(\theta_{\rm off}))\mu_{{\rm off},h}=0$, $\mathbf{1}^{\mathsf T}\mu_{{\rm off},h}=1$ under the offline kernel ($u_{\rm base}=-0.80x$, $\sigma_{\rm off}=0.10$), cross-checked against the cell-binned histogram $\mu_{{\rm off},h}^{\rm emp}$ of the $N=10^5$ offline ensemble ($\norm{\mu_{{\rm off},h}-\mu_{{\rm off},h}^{\rm emp}}_{\ell^1} \approx 0.0355$ at $M=1{,}200$). Holding $\theta_{\rm off}$ fixed at \eqref{eq:theta-off-numerical}, Table~\ref{tab:deployment_ablation} separates the effect of relaxing the baseline gain ($-0.80x \to -0.25x$) from that of increasing the noise ($\sigma: 0.10 \to 0.18$). Across $M \in \{600,\, 1{,}200,\, 2{,}400\}$, all Ulam values of $\norm{r_{A,h}}_{\ell^1}$ and $\norm{\mu_{{\rm dep},h}-\mu_{{\rm off},h}}_{\ell^1}$ vary by less than $6\times 10^{-5}$ (e.g., $\norm{r_{A,h}}_{\ell^1} \in \{0.35216,\, 0.35212,\, 0.35213\}$ for the gain-only change and $\{0.54875,\, 0.54880,\, 0.54879\}$ for the simultaneous change). The ablation confirms that relaxing the baseline gain alone destabilizes the origin ($F_{\rm dep}'(0)\approx 1.2619>1$), triggers the bimodal split ($\norm{\mu_{{\rm dep},h}-\mu_{{\rm off},h}}_{\ell^1}\approx 2.0000$, $\operatorname{std}_{\mu_{\rm deploy}}(X)\approx 1.4139$), and raises $\mathcal{E}_A(\theta_{\rm off})$ from $0.0640$ to $0.2842$, whereas increasing $\sigma$ alone under $-0.80x$ preserves unimodality ($F_{\rm off}'(0) = 0.25-\theta_{1,\rm off} \approx 0.7119 < 1$, $\operatorname{std}_{\mu_{\rm deploy}}(X)\approx 0.2483$, $\mathcal{E}_A(\theta_{\rm off})\approx 0.1090$) and primarily supplies persistent excitation during online recovery.

\begin{table}[htpb]
\centering
\footnotesize
\setlength{\tabcolsep}{2.5pt}
\begin{tabular}{lccccccc}
\toprule
\textbf{Transition ($\theta_{\rm off}$ fixed)} & $u_{\rm base}$ & $\sigma$ & $\norm{r_{A,h}}_{\ell^1}$ \textbf{(Ulam/Emp.)} & $\norm{\mu_{{\rm dep},h}-\mu_{{\rm off},h}}_{\ell^1}$ & $\operatorname{std}(X)$ & $\mathcal{E}_A(\theta_{\rm off})$ & $\lambda_{\min}(H)$ \\
\midrule
Offline ref.\ ($\mu_{{\rm off},h}$) & $-0.80x$ & $0.10$ & $0$ & $0$ & $0.1409$ & $0.0640$ & $4.49\times 10^{-5}$ \\
Gain change only & $-0.25x$ & $0.10$ & $0.3521$ / $0.3528$ & $2.0000$ (bimodal) & $1.4139$ & $0.2842$ & $2.08\times 10^{-2}$ \\
Noise change only & $-0.80x$ & $0.18$ & $0.3603$ / $0.3624$ & $0.5398$ (unimodal) & $0.2483$ & $0.1090$ & $1.24\times 10^{-3}$ \\
Simultaneous change & $-0.25x$ & $0.18$ & $0.5488$ / $0.5501$ & $1.9972$ (bimodal) & $1.3898$ & $0.2825$ & $6.97\times 10^{-2}$ \\
\bottomrule
\end{tabular}
\caption{Separation of baseline gain relaxation and noise intensity change at deployment on the $M=1{,}200$ Ulam grid with $\theta_{\rm off}$ fixed at \eqref{eq:theta-off-numerical}. The column $\norm{r_{A,h}}_{\ell^1}$ reports the one-step invariance residual \eqref{eq:discretized-rA} evaluated on the stationary Ulam vector $\mu_{{\rm off},h}$ and on the cell-binned $N=10^5$ empirical histogram $\mu_{{\rm off},h}^{\rm emp}$.}
\label{tab:deployment_ablation}
\end{table}

Under $\Phi_\delta$ ($\delta=0.04$, $\sigma=0.18$), the trajectory recovers to $(\mu_*, \theta_*)$: the spread of $\mu_{\rm deploy}$ increases $\lambda_{\min}(H)$, restoring a contractive drift slope and returning the state to the unimodal measure $\mu_*$. As the misfit diminishes, the linear-response term $C_*R_*B_*$ decays and $\norm{\theta_n - \theta_{\rm opt}}$ relaxes along the slow manifold on the time scale $O(1/(\delta\lambda_{\min}(H_*)))$.

\begin{figure}[htpb]
\centering
\includegraphics[width=\textwidth]{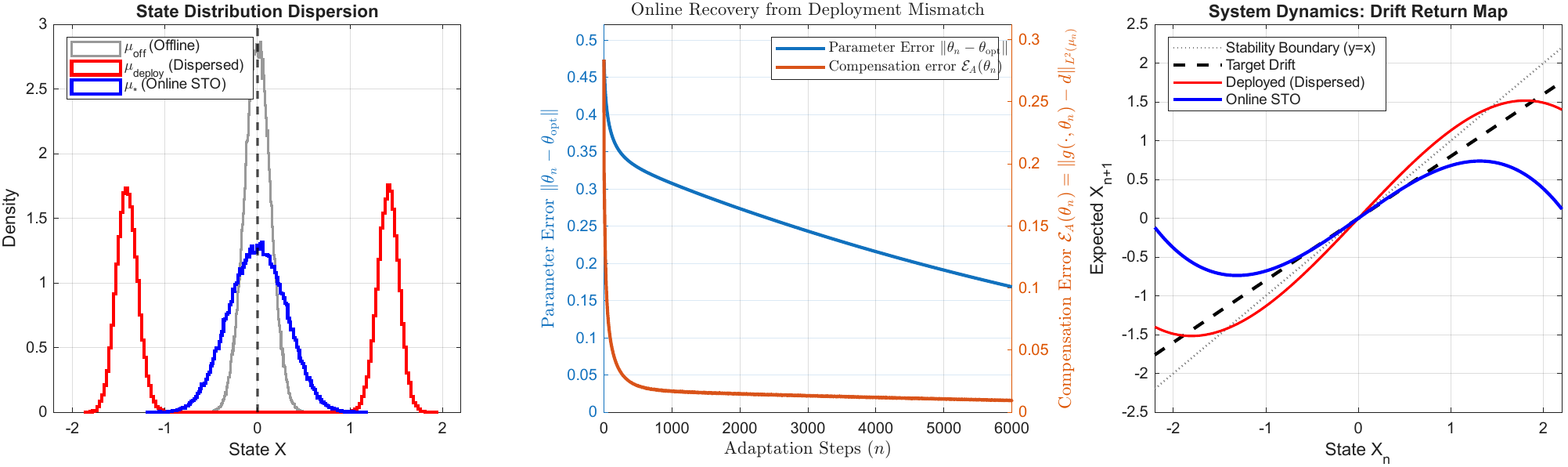}
\caption{Failure of invariance at deployment and online recovery ($\sigma_{\rm off}=0.10 \to \sigma_{\rm online}=0.18$, $\delta=0.04$). (Left) State distribution. High-gain pre-training ($-0.80$) produces $\mu_{\rm off}$ (gray). Relaxing the gain to $-0.25$ and increasing the noise $\sigma = 0.10 \to 0.18$ (Table~\ref{tab:parameters}) yield $\mathfrak R_A \neq 0$ and drive the state to $\mu_{\rm deploy}$ (red); adaptation under $\Phi_\delta$ restores contractivity and returns the state to $\mu_*$ (blue). (Middle) Error trajectories. Online adaptation reduces the compensation error $\mathcal{E}_A(\theta_n)$ (orange), after which both $\mathcal{E}_A$ and $\norm{\theta_n - \theta_{\rm opt}}$ (blue) decay along the slow manifold. (Right) Return map of the deterministic drift. The deployed drift $F_{\rm dep}$ (red) has slope $F_{\rm dep}'(0)>1$ at the origin; adaptation (blue) restores a contractive slope matching $x_{\rm target}$ (dashed black).}
\label{fig:part2_action_error_recovery}
\end{figure}

Conversely, when offline identification uses the deployment gain ($-0.25$) and the same noise law $\mathcal{N}(0,\sigma^2)$ so that $G_{\rm off} \equiv G_{\rm dep}$, let $\mu_{\rm off}$ denote the stationary distribution reached during burn-in with the resulting controller $\theta_{\rm off}^{\rm match} \approx [-1.14\times 10^{-5},\, 0.01079,\, 3.54\times 10^{-5},\, 0.00270]^{\mathsf T}$ held frozen ($F_{\rm dep}'(0) = 0.8 - \theta_{1,\rm off}^{\rm match} \approx 0.7892 < 1$). The offline and deployed kernels then coincide:
\begin{equation}
    P_{\rm dep,\theta_{\rm off}^{\rm match}}\mu_{\rm off} = \mu_{\rm off}, \qquad \mathfrak{R}_A = 0 \quad (\norm{r_{A,h}}_{\ell^1} = 0).
\end{equation}
Figure~\ref{fig:part2_escaping_trap_offline} confirms that under this matched protocol, the deployed distribution remains unimodal and avoids the transient bimodal split.

\begin{figure}[htpb]
\centering
\includegraphics[width=0.7\textwidth]{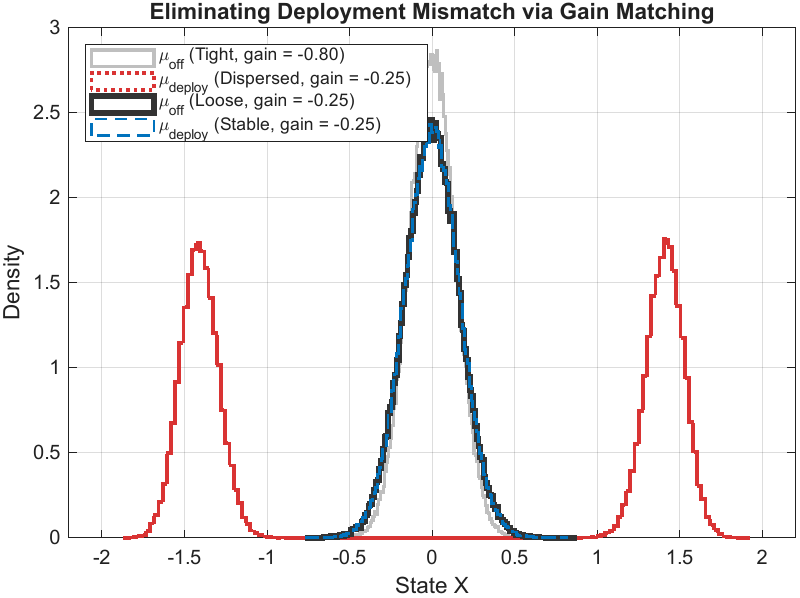}
\caption{Elimination of the invariance residual \eqref{eq:action-residual-ref} at $\sigma=0.10$. Pre-training under $u_{\rm base}(x) = -0.80x$ (solid light gray) induces a kernel mismatch that splits the deployed distribution into $\mu_{\rm deploy}$ (dotted red). Pre-training under the matched baseline $u_{\rm base}(x) = -0.25x$ (solid dark charcoal) yields $\mathfrak{R}_A = 0$, so the deployed distribution (dashed blue) remains unimodal.}
\label{fig:part2_escaping_trap_offline}
\end{figure}

\subsubsection{Loss of Linear Stability and Oscillatory Behavior}
\label{sec:loss-of-stability}

In the unregularized, realizable case ($\lambda=0$, $g(\cdot,\theta_{\rm opt})=d$, $\Gammahat=\Id$), the feedback block $C_*$ vanishes identically and the linearized pair operator \eqref{eq:block-matrix} is upper block-triangular:
\begin{equation}
    \mathcal{A}_\delta = \begin{pmatrix} P_* & B_* \\ 0 & \Id - \delta H_* \end{pmatrix}.
\end{equation}
The multipliers of the parameter block are therefore $\lambda_j^{\rm slow}(\delta) = 1 - \delta \lambda_j(H_*)$ for all $\delta > 0$, and the linear stability threshold is
\begin{equation}
    \delta_c = \frac{2}{\lambda_{\max}(H_*)},
\end{equation}
where the only approximation lies in the numerical evaluation of $H_*$. As $\delta$ increases, the leading multiplier $1 - \delta\lambda_{\max}(H_*)$ crosses zero at $\delta_{\rm osc} := 1/\lambda_{\max}(H_*)$, producing an alternating linearized mode, and exits the unit disk through $-1$ at $\delta_c$ (establishing a rigorous flip bifurcation would require additional nonlinear nondegeneracy conditions).

Figure~\ref{fig:part2_bifurcation_diagram} shows the empirical distribution of $X_n$ across $\delta \in [1.50, 2.25]$ at $\sigma=0.05$ ($\sigma_x \approx 0.0833$). At $\theta_{\rm opt}$, the eigenvalues of $H_*$ are $\lambda_1 \approx 2.01\times 10^{-6}$ ($x^3$ mode), $\lambda_2 \approx 9.64\times 10^{-5}$ ($x^2$ mode), $\lambda_3 \approx 6.95\times 10^{-3}$ ($x$ mode), and $\lambda_4 = \lambda_{\max}(H_*) \approx 1.00005$ (aligned with the constant bias feature $\psi_1(x)=1$, eigenvector $v_4 \approx [1,0,0.0069,0]^{\mathsf T}$), giving $\delta_c = 2/\lambda_{\max}(H_*) \approx 1.9999 \approx 2.000$ ($1.9998$ for the $N=100$ ensemble). The initial perturbation $\theta_0 - \theta_{\rm opt} = [10^{-2}, -10^{-2}, 0, 0]^{\mathsf T}$ projects equally ($70.7\%$ of its Euclidean norm) onto $v_4$ and $v_3$. Above $\delta_c$, one must distinguish the local loss of equilibrium stability from the macroscopic state broadening visible in Figure~\ref{fig:part2_bifurcation_diagram}. On a grid with step $\Delta\delta \approx 1.5\times 10^{-3}$ ($500$ points in $[1.50, 2.25]$), for $\delta \in (2.00, 2.19)$, the unstable multiplier $1 - \delta\lambda_{\max}(H_*) \in (-1.19, -1.00)$ amplifies sign-alternating perturbations along the constant offset $\theta_0$ (at $\delta=2.18$, linear growth amplifies a perturbation by $10^2$ in $\log 100/\log 1.18 \approx 28$ steps). Since $H(\mu)_{11} = \int_X 1^2\dd\mu = 1$ implies $\lambda_{\max}(H(\mu)) \ge 1$ and $1 - \delta\lambda_{\max}(H(\mu)) < -1$ for all $\delta > 2$, this instability cannot saturate via the Gram matrix alone. Instead, as $\theta_{0,n} \approx (-1)^n a_0$ alternates in sign, the mean state $m_{1,n} = \mathbb{E}_{\mu_n}[X_n] \approx \frac{a_0}{1.8 - \bar{\theta}_1}(-1)^n$ oscillates in phase with $\theta_{0,n}$, coupling the odd coefficients into the bias update $G_0(\mu_n,\theta_n) = \theta_{0,n} + \theta_{1,n}m_{1,n} + \theta_{2,n}m_{2,n} + (\theta_{3,n}-0.15)m_{3,n} \approx \bigl(1 + \frac{\bar{\theta}_1}{1.8 - \bar{\theta}_1}\bigr)\theta_{0,n}$. A negative steady shift of the linear gain ($\bar{\theta}_1 \approx -0.056$ at $\delta=2.05$, $-0.135$ at $\delta=2.15$, and $-0.166$ at $\delta=2.18$, with $\bar{\theta}_3 \approx 0.175$) reduces this effective loop gain to $2/\delta$, saturating the primary instability into a bounded period-$2$ orbit ($\norm{\theta_n - \theta_{\rm opt}} \approx 0.10$--$0.31$, with state excursions $\vert X_n(1)\vert \lesssim 0.08$ over $n \in [1{,}501, 1{,}600]$ that remain near zero on the $[-2, 2]$ scale). Because the closed-loop linear slope at the origin is $0.8 - \theta_1$, however, shifting $\theta_1$ negatively weakens physical contractivity: when $\delta$ exceeds the secondary threshold $\delta \approx 2.19$ ($2.1895 \pm 0.0028$ for $N=10^2$ and $2.1950 \pm 0.0000$ for $N=10^4$ across $10$ seeds), $\theta_{1,n}$ crosses $-0.20$ ($0.8 - \theta_{1,n} > 1$), the physical state loses contractivity, and boundary reflections drive trajectories across $[-2, 2]$. Figure~\ref{fig:part2_density_shattering} compares the evolution of $\mu_n(x)$ over $n \in [1, 150]$ for $\delta=1.85 \in (\delta_{\rm osc}, \delta_c)$ (damped sign-alternating density oscillations, multiplier $\approx -0.85$) and $\delta=2.22 > \delta_c$ (supercritical dispersion across $[-2, 2]$, multiplier $\approx -1.22$).

\begin{figure}[htpb]
\centering
\includegraphics[width=0.7\textwidth]{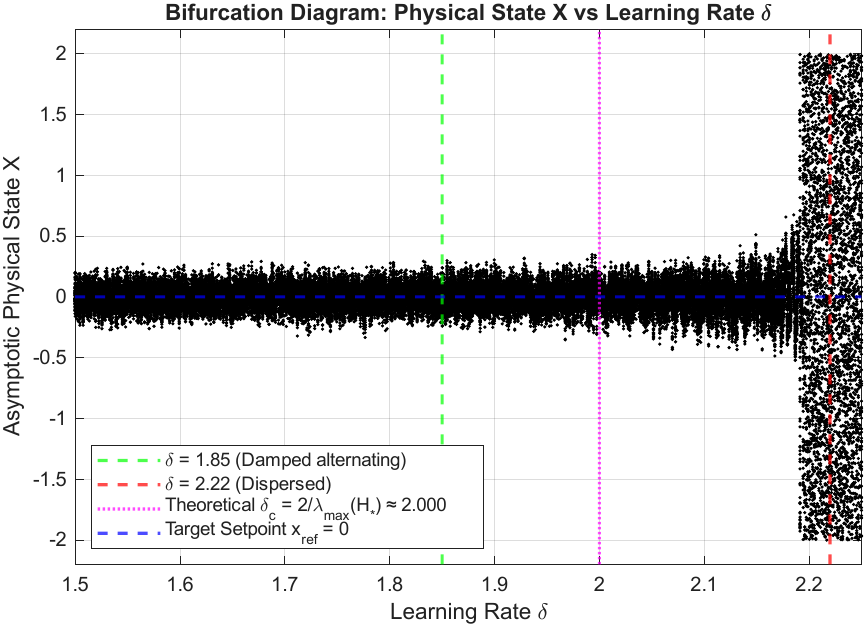}
\caption{Empirical bifurcation diagram of $X_n \in [-2, 2]$ versus $\delta \in [1.50, 2.25]$ at $\sigma=0.05$. Vertical dashed lines mark $\delta=1.85$ (green) and $\delta=2.22$ (red); the dotted magenta line marks $\delta_c = 2/\lambda_{\max}(H_*) \approx 2.000$. Below $\delta_c$, the equilibrium is stable and the state remains near $x=0$ (dashed blue); for $\delta \in (\delta_c, 2.19)$, the primary instability along $\theta_0$ saturates into a bounded period-2 regime around $x=0$, whereas for $\delta \gtrsim 2.19$ (grid step $\Delta\delta \approx 1.5\times 10^{-3}$) a secondary nonlinear instability spreads trajectories across $[-2, 2]$.}
\label{fig:part2_bifurcation_diagram}
\end{figure}

\begin{figure}[htpb]
\centering
\includegraphics[width=\textwidth]{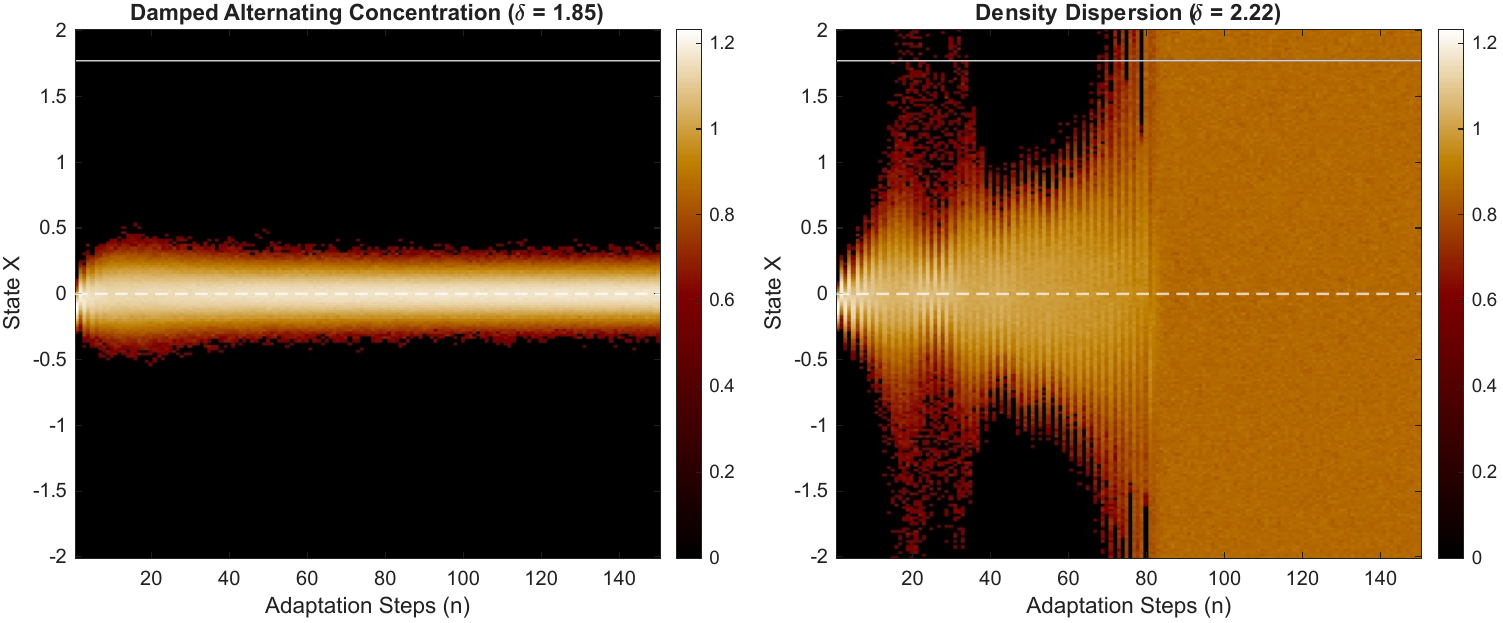}
\caption{Evolution of the state density $\mu_n(x)$ over $n \in [1, 150]$ from $\theta_0=[0.10, -0.10, 0, 0]^{\mathsf T}$ ($N=10^4, \sigma=0.05$). Color intensity shows $(\mathrm{PDF})^{0.1}$. (Left) Damped alternating regime ($\delta=1.85$, multiplier $\approx -0.85 \in (-1,0)$). (Right) Supercritical regime ($\delta=2.22 > \delta_c$, multiplier $\approx -1.22 < -1$), where sign-alternating oscillations amplify and disperse the density across $[-2, 2]$.}
\label{fig:part2_density_shattering}
\end{figure}

\subsection{Numerical Exploration of Non-Convex Neural Parameterizations}
\label{sec:nonconvex_extension}

While the explicit verification in Section~\ref{sec:duffing-example} uses the fixed-feature linear-readout model $g(x,\theta)=\theta^{\mathsf T}\psi(x)$, the stationary reduction and Brouwer degree criterion of Section~\ref{4} do not require convexity of the loss landscape: if the frozen dynamics admits a unique continuous stationary branch $\theta\mapsto\mu_\theta$ and $\overline G(\theta) = G(\mu_\theta,\theta)$ satisfies $0 \notin \overline G(\partial\Omega)$ and $\deg(\overline G,\Omega,0)\neq 0$ on a bounded open set $\Omega\subset\R^p$ (Theorem~\ref{thm:degree-existence}), then $\STO_\delta$ admits an adaptive equilibrium $\rho_* = \mu_{\theta_*}\otimes\delta_{\theta_*}$. When internal network weights are adapted, however, feature sensitivities couple to the readout magnitude, and $\lambda\norm{\theta}^2$ alone no longer controls the radial growth of $\langle\overline G(\theta),\theta\rangle$ on large spheres because the features $\tanh(W_1x+b_1)$ depend on $\theta$. Without establishing analytical boundary estimates for this case, we test the coupled dynamics numerically by replacing the polynomial basis with a fully adaptive multilayer perceptron ($8$ $\tanh$ hidden units, $p=25$ parameters):
\begin{equation}
    g(x, \theta) = W_2 \tanh(W_1 x + b_1) + b_2, \qquad \theta = \mathrm{vec}(W_1, b_1, W_2, b_2) \in \R^{25}.
    \label{eq:mlp-compensator}
\end{equation}
For initial weights $\theta_0 \sim \mathcal{N}(0, 0.1^2 I_{25})$ and parameters $\delta=0.05$, $\sigma=0.10$, and $\lambda=10^{-5}$, the empirical tracking loss converges smoothly to the Gaussian noise floor $\sigma^2/2 = 0.005$ (Figure~\ref{fig:mlp_convergence}).

\begin{figure}[htpb]
\centering
\includegraphics[width=0.7\textwidth]{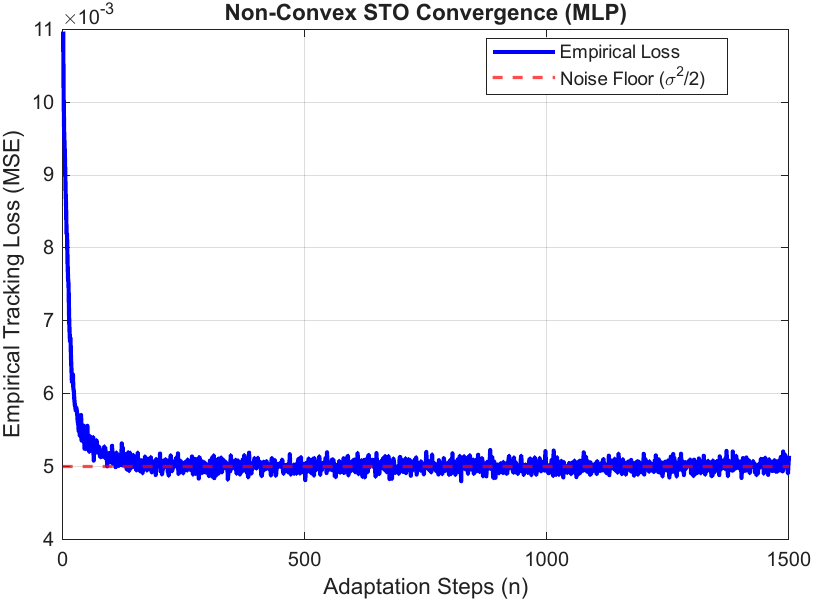}
\caption{Autonomous adaptation of an MLP compensator ($p=25$, $\delta=0.05$, $\sigma=0.10$, $\lambda=10^{-5}$). The empirical tracking loss converges to the Gaussian noise floor $\sigma^2/2 = 0.005$ (dashed red).}
\label{fig:mlp_convergence}
\end{figure}

\subsection{Reproducibility Details and Algorithmic Protocols}
\label{sec:reproducibility-protocols}

Table~\ref{tab:reproducibility_vectors} lists the initial conditions $(x_0^{(i)},\theta_0)$, random seeds, ensemble sizes $N$, realization counts $R$, and resulting parameter vectors for all experiments in Table~\ref{tab:parameters}. All particle simulations use the MATLAB Mersenne Twister generator initialized at seed $42$ (\texttt{rng(42)}), averaging over $N \in \{10^2, 10^4, 10^5\}$ coupled particles driven by independent Gaussian draws $\xi_n^{(i)} \sim \mathcal{N}(0,\sigma^2)$. The secondary instability threshold in Section~\ref{sec:loss-of-stability} is evaluated across $R=10$ independent realizations with seeds $s \in \{1,\ldots,10\}$.

\begin{table}[htpb]
\centering
\scriptsize
\setlength{\tabcolsep}{1.8pt}
\begin{tabular}{llll}
\toprule
\textbf{Experiment} & \textbf{Initial $(x_0^{(i)},\, \theta_0^{\mathsf T})$} & \textbf{Seed ($R$)} & \textbf{Equilibrium / Final Parameter Vector $\theta^{\mathsf T}$} \\
\midrule
Figs.~\ref{fig:part1_nominal_convergence}--\ref{fig:part1_timescale_and_reg} ($\lambda=0$) & $0$, $[-0.25, 0.45, 0, 0]$ & $42$ ($1$) & $\theta_{\rm opt}=[0,0,0,0.15]$; $\hat\theta_{40{\rm k}}\approx[2.1\cdot 10^{-4},4.4\cdot 10^{-5},-1.1\cdot 10^{-4},0.1502]$ \\
Figs.~\ref{fig:part1_nominal_convergence}--\ref{fig:part1_timescale_and_reg} ($\lambda=10^{-3}$) & $0$, $[-0.25, 0.45, 0, 0]$ & $42$ ($1$) & $\theta_*^\lambda\approx[0, 0.00699, 0, 0.12429]$; $\hat\theta_*^\lambda\approx[-10^{-5}, 0.00699, 4\cdot 10^{-5}, 0.12435]$ \\
Fig.~\ref{fig:linear_response_contribution} (Odd / Even) & $\mu_{\theta_*}$, $\theta_*+\Delta\theta_0$ & $42$ ($1$) & $\theta_*\approx[0, 0.11212, 0]$ ($\Delta\theta_{0,\rm odd}=[0,0.12,0]$, $\norm{\Delta\theta_{0,\rm even}}=0.08$) \\
Fig.~\ref{fig:part2_identifiability_trap_pe} ($\sigma=0.10$) & $0$, $[-0.25, 0.45, 0, 0]$ & $42$ ($1$) & $\hat\theta_{40{\rm k}}\approx[1.1\cdot 10^{-4}, 0.00825, -2.0\cdot 10^{-5}, 0.05427]$ \\
Fig.~\ref{fig:part2_action_error_recovery} (Offline) & $0$, $[0, -0.45, 0, 0.28]$ & $42$ ($1$) & $\theta_{\rm off}\approx[1.0\cdot 10^{-5}, -0.46188, 4.6\cdot 10^{-7}, 0.27928]$ \\
Fig.~\ref{fig:part2_action_error_recovery} (Online) & $\mu_{\rm deploy}^N$, $\theta_{\rm off}$ & $42$ ($1$) & $\hat\theta_{6{\rm k}}\approx[-1.1\cdot 10^{-4}, -0.04287, 2.2\cdot 10^{-5}, 0.31305]$ \\
Fig.~\ref{fig:part2_escaping_trap_offline} (Matched) & $0$, $[0, 0, 0, 0]$ & $42$ ($1$) & $\theta_{\rm off}^{\rm match}\approx[-1.1\cdot 10^{-5}, 0.01079, 3.5\cdot 10^{-5}, 0.00270]$ \\
Fig.~\ref{fig:part2_bifurcation_diagram} (Sweep $\delta$) & $0$, $[0.01, -0.01, 0, 0.15]$ & $42$, $1$--$10$ ($10$) & Period-$2$ saturation for $\delta\in(2.00,2.19)$; runaway for $\delta\gtrsim 2.19$ \\
Fig.~\ref{fig:part2_density_shattering} (Density) & $0$, $[0.10, -0.10, 0, 0]$ & $42$ ($1$) & Damped alternating ($\delta=1.85$) vs.\ dispersed ($\delta=2.22$) \\
Fig.~\ref{fig:mlp_convergence} (MLP, $p=25$) & $0$, $\mathcal{N}(0,0.01 I_{25})$ & $42$ ($1$) & Stationary tracking MSE $\approx \sigma^2/2 = 0.005$ \\
\bottomrule
\end{tabular}
\caption{Initial conditions, random seeds, realization counts $R$, and learned/equilibrium parameter vectors for the numerical experiments.}
\label{tab:reproducibility_vectors}
\end{table}

\paragraph{Algorithmic protocol for offline pre-training, deployment burn-in, and online STO adaptation.}
Given ensemble size $N$, feature map $\psi(x)$, regularization weight $\lambda \ge 0$ (with $\theta_c=0$), and periodic reflection $\mathfrak{r}_2(z) = \lvert ((z-2)\bmod 8) - 4\rvert - 2$, the simulations execute the following steps:
\begin{enumerate}[label=\textbf{Phase~\arabic*:}, leftmargin=4.2em]
    \item \emph{Offline pre-training ($n=0,\ldots,N_{\rm off}-1$, used in Figures~\ref{fig:part2_action_error_recovery}--\ref{fig:part2_escaping_trap_offline}):} Initialize $x_0^{(i)}=0$ ($i=1,\ldots,N$) and $\theta_0 = \theta_{\rm off,0}$. At each step $n$, draw i.i.d.\ $\xi_n^{(i)}\sim\mathcal{N}(0,\sigma_{\rm off}^2)$ and update
    \[
        x_{n+1}^{(i)} = \mathfrak{r}_2\!\bigl(1.05x_n^{(i)} + 0.15(x_n^{(i)})^3 + u_{\rm base}^{\rm off}(x_n^{(i)}) - \theta_n^{\mathsf T}\psi(x_n^{(i)}) + \xi_n^{(i)}\bigr),
    \]
    \[
        \widehat G_{n}^{\rm off} = -\frac{1}{N}\sum_{i=1}^N \bigl(x_{n+1}^{(i)} - 0.8x_n^{(i)}\bigr)\psi(x_n^{(i)}), \qquad \theta_{n+1} = \theta_n - \delta_{\rm off}\widehat G_{n}^{\rm off}.
    \]
    Set $\theta_{\rm off} := \theta_{N_{\rm off}}$ and record the empirical offline distribution $\mu_{\rm off}^N = \frac{1}{N}\sum_{i=1}^N \delta_{x_{N_{\rm off}}^{(i)}}$.
    \item \emph{Frozen deployment burn-in ($n=0,\ldots,N_{\rm burn}-1$, used in Figures~\ref{fig:part2_action_error_recovery}--\ref{fig:part2_escaping_trap_offline}):} Hold $\theta \equiv \theta_{\rm off}$ frozen, switch the baseline controller to $u_{\rm base}^{\rm dep}(x) = -0.25x$, and evolve $x_{n+1}^{(i)} = \mathfrak{r}_2\bigl(0.8x_n^{(i)} + 0.15(x_n^{(i)})^3 - \theta_{\rm off}^{\mathsf T}\psi(x_n^{(i)}) + \xi_n^{(i)}\bigr)$ for $N_{\rm burn}=1{,}000$ steps to obtain the deployed distribution $\mu_{\rm deploy}^N$.
    \item \emph{Online STO adaptation ($n=0,\ldots,N_{\rm steps}-1$, used in all experiments; autonomous runs skip Phases~1--2):} From $(x_0^{(i)},\theta_0)$ (set to $(\mu_{\rm deploy}^N,\theta_{\rm off})$ after Phase~2, or to the autonomous initial condition in Table~\ref{tab:reproducibility_vectors}), draw $\xi_n^{(i)}\sim\mathcal{N}(0,\sigma^2)$ and update the coupled state--parameter pair under $u_{\rm base}^{\rm dep}(x)=-0.25x$:
    \[
        x_{n+1}^{(i)} = \mathfrak{r}_2\!\bigl(0.8x_n^{(i)} + 0.15(x_n^{(i)})^3 - g(x_n^{(i)},\theta_n) + \xi_n^{(i)}\bigr),
    \]
    \[
        \theta_{n+1} = \theta_n - \delta\left[-\frac{1}{N}\sum_{i=1}^N \bigl(x_{n+1}^{(i)} - 0.8x_n^{(i)}\bigr)\nabla_\theta g(x_n^{(i)},\theta_n) + \lambda\theta_n\right]
    \]
    (with $-(x_{n+1}^{(i)}-0.8x_n^{(i)})$ replaced by the exact unreflected residual $g(x_n^{(i)},\theta_n)-d(x_n^{(i)})$ in the $\sigma=0.24$ linear-response experiment of Figure~\ref{fig:linear_response_contribution}).
\end{enumerate}

\section{Conclusion} \label{9}

We have formulated autonomous neural adaptive control as a self-consistent transfer operator on the extended state--parameter space $\Z=X\times\ThetaSet$. 

This formulation makes three analytical tools available for coupled plant--learner systems:
\begin{enumerate}[label=(\roman*)]
    \item \emph{Finite-dimensional stationary reduction:} When each frozen controller admits a unique invariant physical distribution $\mu_\theta \in B_s(X)$, adaptive equilibria $(\mu_*,\theta_*)$ correspond bijectively to zeros of the reduced learning field $\overline G(\theta) = G(\mu_\theta,\theta)$ (Theorem~\ref{thm:equilibrium-correspondence}). In this setting, offline-to-online deployment mismatch (\emph{Action Error}) corresponds to a nonzero physical invariance residual $\mathfrak R_A = P_{\theta_{\rm off}}\mu_{\rm off} - \mu_{\rm off} \neq 0$, and Brouwer degree certifies equilibrium existence from boundary coercivity without requiring global convexity (Theorem~\ref{thm:degree-existence}).
    \item 
\emph{Local stability for slow adaptation:}
Under the differentiability assumptions of
Theorem~\ref{thm:local-stability}, a spectral gap for
the frozen physical operator and positive real parts
of all eigenvalues of $\Gammahat\mathcal J_*$
ensure local exponential convergence of both the
physical distribution and the controller parameters
to the adaptive equilibrium, for every sufficiently
small positive learning rate $\delta$.
The slow relaxation time is of order $\delta^{-1}$,
so reducing the learning rate also slows convergence.
    \item
\emph{Statistical feedback and adaptive stability:}
Linear response quantifies how the stationary physical
distribution changes with the controller. It yields
\[
  D\overline G(\theta_*)
  =H_*+C_*R_*B_*=\mathcal J_*
\]
(Proposition~\ref{prop:effective-jacobian-identification}).
Here $H_*$ describes the learning response at a fixed
physical distribution, while $C_*R_*B_*$ accounts for
the change in that distribution induced by the controller.
Thus the derivative of the reduced equilibrium equation
also determines the stability of slow adaptation.
\end{enumerate}

Several extensions remain open. Applying these criteria to higher--dimensional robotic systems, such as compliant aerial interaction \cite{AuconePallottino2024,TognonPallottino2018} or exoskeleton gait adaptation \cite{lippi2021icinco}, requires verifying uniform mixing, stationary-branch differentiability, and degree boundary conditions for the corresponding physical kernels. Extending the explicit estimates from linear readouts to deep nonlinear parameterizations (motivated empirically by Section~\ref{sec:nonconvex_extension}) requires accounting for parameter symmetries, potentially on a quotient or policy space, and controlling the radial growth of internal layer gradients. Finally, incorporating finite-batch stochastic approximation noise, actor--critic states, replay buffers, or time--dependent references will require enlarging the autonomous state space or weakening the strong product--pace differentiability assumptions via two--norm transfer--operator techniques \cite{CastorriniGalatoloTanziJNS}.

\section*{Code and Data Availability}
The simulation code, Ulam transfer-operator discretization solvers, precomputed datasets, and figure-generation scripts reproducing all numerical experiments in Section~\ref{sec:duffing-example} and Section~\ref{sec:empirical-validation} are openly available at
\begin{center}
  \url{https://github.com/VittorioFreiburg/sto-adaptive-control}
\end{center}
and permanently archived on Zenodo under DOI \href{https://doi.org/10.5281/zenodo.24112629}{\texttt{10.5281/zenodo.24112629}}.

\end{document}